\documentclass[11pt]{article}
\usepackage[utf8]{inputenc}
\usepackage[top=1in,bottom =1in,left=1in,right=1in]{geometry}

\usepackage{amsmath,amssymb,mathtools,bm,isomath,amsthm}
\usepackage{graphicx}
\usepackage{epstopdf}
\usepackage{color}
\usepackage[dvipsnames]{xcolor}
\usepackage{soul}
\usepackage{enumerate}
\usepackage{nicefrac}
\usepackage{enumitem}
\usepackage{mathrsfs}
\usepackage{scalerel,stackengine}
\usepackage{afterpage}
\usepackage{url}
\usepackage[affil-it]{authblk} 
\usepackage[toc,page]{appendix}
\usepackage{lineno}
\usepackage{physics}
\usepackage{tikz}
\usetikzlibrary{arrows}

\stackMath
\newcommand\reallywidehat[1]{%
	\savestack{\tmpbox}{\stretchto{%
			\scaleto{%
				\scalerel*[\widthof{\ensuremath{#1}}]{\kern-.6pt\bigwedge\kern-.6pt}%
				{\rule[-\textheight/2]{1ex}{\textheight}}%WIDTH-LIMITED BIG WEDGE
			}{\textheight}% 
		}{0.5ex}}%
	\stackon[1pt]{#1}{\tmpbox}%
}

\newcommand{\R}{\mathbb{R}}

\newcommand{\bx}{{\boldsymbol x}}

\newcommand{\bW}{{\boldsymbol W}}
\newcommand{\bb}{{\boldsymbol b}}

\newcommand{\bz}{{\boldsymbol z}}

\newtheorem{example}{Example}

\begin{document}
\setlength{\baselineskip}{20pt}

\title{Multi-level Neural Networks for Accurate Solutions of Boundary-Value Problems}
	
\author[1]{Ziad Aldirany}
\author[2]{R\'egis Cottereau}
\author[1]{Marc Laforest}
\author[1]{Serge Prudhomme}
\affil[1]{D\'epartement de math\'ematiques et de g\'enie industriel, \protect\\ Polytechnique Montr\'eal, \protect\\ Montr\'eal, Qu\'ebec, Canada}
\affil[2]{Aix-Marseille Universit\'e, CNRS, Centrale Marseille, \protect\\ LMA UMR 7031, \protect\\ Marseille, France}
	
\date{}
	
\maketitle
	
\begin{abstract}
The solution to partial differential equations using deep learning approaches has shown promising results for several classes of initial and boundary-value problems. However, their ability to surpass, particularly in terms of accuracy, classical discretization methods such as the finite element methods, remains a significant challenge. Deep learning methods usually struggle to reliably decrease the error in their approximate solution. A new methodology to better control the error for deep learning methods is presented here. The main idea consists in computing an initial approximation to the problem using a simple neural network and in estimating, in an iterative manner, a correction by solving the problem for the residual error with a new network of increasing complexity. This sequential reduction of the residual of the partial differential equation allows one to decrease the solution error, which, in some cases, can be reduced to machine precision. The underlying explanation is that the method is able to capture at each level smaller scales of the solution using a new network. Numerical examples in 1D and 2D are presented to demonstrate the effectiveness of the proposed approach. This approach applies not only to physics informed neural networks but to other neural network solvers based on weak or strong formulations of the residual.
\end{abstract}
	
	
%============================================================
\section{Introduction}
\label{sec:introduction}
%============================================================

In recent years, the solution of partial differential equations using deep learning~\cite{sirignano2018dgm,bihlo2022physics,jin2021nsfnets} has gained popularity and is emerging as an alternative to classical discretization methods, such as the finite element or the finite volume methods. Deep learning techniques can be used to either solve a single initial boundary-value problem~\cite{raissi2019physics,weinan2018deep,zang2020weak} or approximate the operator associated with a partial differential equation~\cite{lu2019deeponet,li2020fourier,aldirany2023operator, patel2022variationally}. The primary advantages of deep learning approaches lie in their ability to provide meshless methods, and hence address the curse of dimensionality, and in the universality of their implementation for various initial and boundary-value problems. However, one of the main obstacles remains their inability to consistently reduce the relative error in the computed solution. Although the universal approximation theorem~\cite{cybenko1989approximation, hornik1991approximation} guarantees that a single hidden layer network with a sufficient width should be able to approximate smooth functions to a specified precision, one often observes in practice that the convergence with respect to the number of iterations reaches a plateau, even if the size of the network is increased. This is primarily due to the use of gradient-based optimization methods, e.g.\ Adam~\cite{DBLP:journals/corr/KingmaB14}, for which the solution may get trapped in local minima. These optimization methods applied to classical neural network architectures, e.g.\ feedforward neural networks~\cite{lecun2015deep}, do indeed experience difficulties in controlling the large range of scales inherent to a solution, even with some fine-tuning of the hyper-parameters, such as the learning rate or the size of the network. In contrast, this is one of the main advantages of classical methods over deep learning methods, in the sense that they feature well-defined techniques to consistently reduce the error, using for instance mesh refinement~\cite{bangerth2003adaptive,sevilla2022mesh} or multigrid structures~\cite{hackbusch1985multi}. 

We introduce in this work a novel approach based on the notion of multi-level neural networks, which are designed to consistently reduce the residual associated with a partial differential equation, and hence, the errors in the numerical solution. The approach is versatile and can be applied to various neural network methods that have been developed for the solution of boundary-value problems~\cite{weinan2018deep,zang2020weak}, but we have chosen, for the sake of simplicity, to describe the method on the particular case of physics-informed neural networks (PINNs)~\cite{raissi2019physics}. Once an approximate solution to a linear boundary-value problem has been computed with the classical PINNs, the method then consists in finding a correction, namely, estimating the solution error, by minimizing the residual using a new network of increasing complexity. The process can subsequently be repeated using additional networks to minimize the resulting residuals, hence allowing one to reduce the error to a desired precision. A similar idea has been proposed in~\cite{ainsworth2021galerkin} to control the error in the case of symmetric and positive-definite variational equations. Using Galerkin neural networks, the authors construct basis functions calculated from a sequence of neural networks to generate a finite-dimensional subspace, in which the solution to the variational problem is then approximated. Our approach is more general as the problems do not need to be symmetric. 

The development of the proposed method is based on two key observations. First, each level of the correction process introduces higher frequencies in the solution error, as already discussed in~\cite{ainsworth2021galerkin} and highlighted again in the numerical examples. This is the reason why the sequence of neural networks should be of increasing complexity. Moreover, a key ingredient will be to use the Fourier feature mapping approach~\cite{tancik2020fourier} to accurately approximate the functions featuring high frequencies. Second, the size of the error, equivalently of the residual, becomes at each level increasingly smaller. Unfortunately, feedforward neural networks employing standard parameter initialization, e.g.\ Xavier initialization~\cite{glorot2010understanding} in our case, are tailored to approximate functions whose magnitudes are close to unity. We thus introduce a normalization of the solution error at each level based on the Extreme Learning Method~\cite{huang2011extreme}, which also contributes to the success of the multi-level neural networks.

After finalizing the writing of the manuscript, one has brought to our attention the recent preprint~\cite{Wang-Lai2023-multistage} on multi-stage neural networks. Although the conceptual approach presented in that preprint features many similarities with our method, namely the use of a sequence of networks for the reduction of the numerical errors, the methods developed in our independent work to address the two aforementioned issues are original and sensibly differ from those introduced in~\cite{Wang-Lai2023-multistage}. 

The paper is organized as follows. We briefly describe in Section~\ref{sec:preliminaries} neural networks and their application with PINNs, the deep learning approach that will be used to solve the boundary-value problems at each level of the training. We describe in Section~\ref{sec:error} the two issues that may affect the accuracy of the solutions obtained by PINNs. We motivate in Section~\ref{sec:normalize} the importance of normalization of the problem data and show that it can greatly improve the convergence of the solution. We continue in Section~\ref{sec:high_freq} with the choice of the network architecture and the importance of using the Fourier feature mapping algorithm to approximate high-frequency functions. We then present in Section~\ref{sec:multi-nets} our approach, the multi-level neural network method, and demonstrate numerically with a simple 1D Poisson problem that the method greatly improves the accuracy of the solution, up to machine precision, with respect to the \(L^2\) and the \(H^1\) norms, as in classical discretization methods. We demonstrate further in Section~\ref{sec:numresults} the efficiency of the proposed method on several numerical examples based on the Poisson equation, the convective-diffusion equation, and the Helmholtz equation, in one dimension or two dimensions. We were able to consistently reduce the solution error in these problems using the multi-level neural network method. Finally, we compile concluding remarks about the present work and put forward new directions for research in Section~\ref{sec:conclusion}.

%============================================================
\section{Preliminaries}
\label{sec:preliminaries}
%============================================================

%============================================================
\subsection{Neural networks}

Neural networks have been extensively studied in recent years for solving partial differential equations~\cite{sirignano2018dgm, raissi2019physics}. A neural network can be viewed as a mapping between an input and an output by means of a composition of linear and nonlinear functions with adjustable weights and biases. Training a neural network consists in optimizing the weights and biases by minimizing some measure of the error between the output of the network and corresponding target values obtained from a given training dataset. As a predictive model, the trained network is then expected to provide accurate approximations of the output when considering a wider set of inputs. Several neural network architectures, e.g.\ convolutional neural networks (CNNs)~\cite{krizhevsky2017imagenet} or feedforward neural networks (FNNs)~\cite{lecun2015deep}, are adapted to specific classes of problems.
    
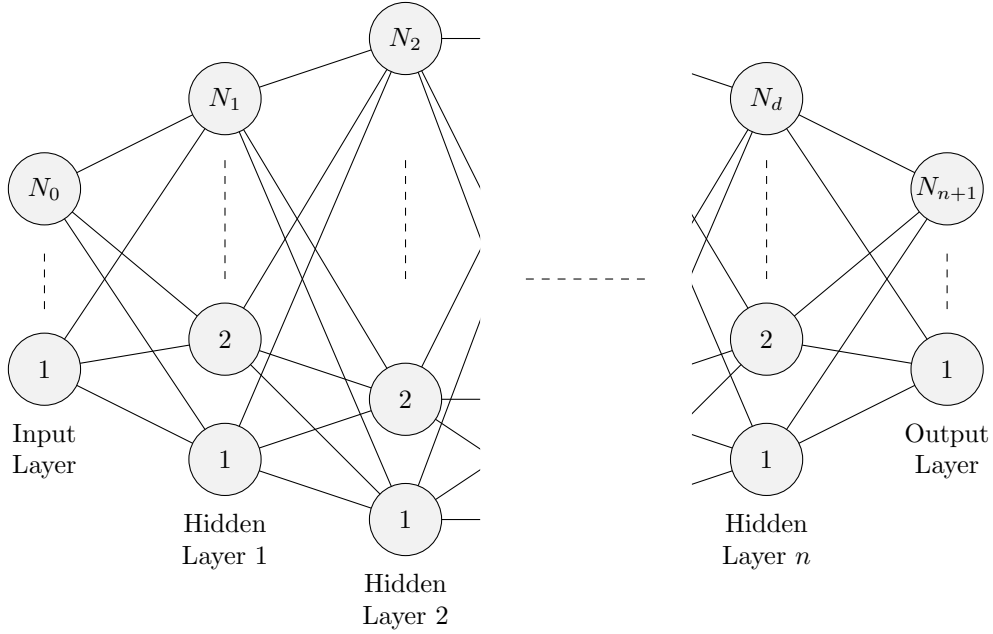
\begin{figure}[tb]
\small
\centering
\begin{tikzpicture}[scale=0.8]
\foreach \i in {1,3,7}{\foreach \j in {2.5,5.5}{\draw (0,\j) -- (3,\i);};}
\foreach \i in {0,2,8}{\foreach \j in {1,3,7}{\draw (3,\j) -- (6,\i);};}
\foreach \i in {0,2,8}{\foreach \j in {0,2,8}{\draw (6,\j) -- (9,\i);};}
\foreach \i in {1,3,7}{\foreach \j in {0,2,8}{\draw (9,\j) -- (12,\i);};}
\foreach \i in {1,3,7}{\foreach \j in {2.5,5.5}{\draw (12,\i) -- (15,\j);};}
\foreach \i in {2.5,5.5}{\fill[gray!10!white,draw=black] (0,\i) circle (0.6);}
\foreach \i in {1,3,7}{\fill[gray!10!white,draw=black] (3,\i) circle (0.6);}
\foreach \i in {0,2,8}{\fill[gray!10!white,draw=black] (6,\i) circle (0.6);}
\foreach \i in {1,3,7}{\fill[gray!10!white,draw=black] (12,\i) circle (0.6);}
\foreach \i in {2.5,5.5}{\fill[gray!10!white,draw=black] (15,\i) circle (0.6);}
\fill[white,draw=white] (7.25,-0.5) rectangle (10.75,8.5);
\foreach \i in {4}{\draw[dashed] (8,\i) -- (10,\i);}
\foreach \i in {3,6,12}{\draw[dashed] (\i,4) -- (\i,6);}
\foreach \i in {0,15}{\draw[dashed] (\i,3.5) -- (\i,4.5);}
\node at (0,1.75) [anchor=north,align=center] {Input \\ Layer};
\node at (0,2.5) {1};
\node at (0,5.5) {$N_0$};
\node at (3,0.25) [anchor=north,align=center] {Hidden \\ Layer 1};
\node at (3,1) {1};
\node at (3,3) {2};
\node at (3,7) {$N_1$};
\node at (6,-0.75) [anchor=north,align=center] {Hidden \\ Layer 2};
\node at (6,0) {1};
\node at (6,2) {2};
\node at (6,8) {$N_2$};
\node at (12,0.25) [anchor=north,align=center] {Hidden \\ Layer $n$};
\node at (12,1) {1};
\node at (12,3) {2};
\node at (12,7) {$N_d$};
\node at (15,1.75) [anchor=north,align=center] {Output \\ Layer};
\node at (15,2.5) {$1$};
\node at (15,5.5) {$N_{n+1}$};
\end{tikzpicture}

\caption{Sketch of a feedforward neural network with $d$ hidden layers of a width \(N_i\), \(i=1,\ldots,n\), an input layer of size \(N_0\), and an output layer of size \(N_{n+1} \).}
\label{fig:FNNs}
\end{figure}

We shall consider here FNNs featuring \(n\) hidden layers, each layer having a width \(N_i\), \(i=1,\ldots,n\), an input layer of width \(N_0\), and an output layer of width \(N_{n+1} \); see Figure~\ref{fig:FNNs}. Denoting the activation function by \(\sigma\), the neural network with input \(\bz_0 \in \mathbb R^{N_0}\) and output \(\bz_{n+1} \in \mathbb R^{N_{n+1}}\) is defined as
\begin{equation}
\label{eqn:FNN}
\begin{aligned}
	&\text{Input layer: }&&  \bz_0, \\
	&\text{Hidden layers: }&&  \bz_{i} = \sigma ( W_i \bz_{i-1} + \bb_i), \quad i=1,\cdots,n, \\ 	
	&\text{Output layer: }&& \bz_{n+1} = W_{n+1} \bz_{n} + \bb_{n+1} ,
\end{aligned}
\end{equation}
where \(W_i\) is the \textit{weights} matrix of size \(N_{i} \times N_{i-1}\) and \(\bb_i \) is the \textit{biases} vector of size \( N_{i}\). To simplify the notation, we combine the weights and biases of the neural network into a single parameter denoted by \(\theta\). The neural network~\eqref{eqn:FNN} generates a finite-dimensional space of dimension $N_\theta=\sum_{i=1}^{n+1} N_{i} (N_{i-1} + 1)$. To keep things simple, throughout this work we shall use the \(\tanh\) activation function and the associated Xavier initialization scheme~\cite{glorot2010understanding} to initialize the weights and biases.

%============================================================
\subsection{Physics-informed neural networks}
\label{subsec:PINNs}

We briefly review the PINNs approach to solving partial differential equations, as described in~\cite{raissi2019physics}. Let \(\Omega \) be an open bounded domain in \(\R^d\), \(d=1,2,\) or \(3\), with boundary \(\partial \Omega\). For two Banach spaces \(U\) and \(V\) of functions over \(\Omega\), we assume a linear differential operator \(A : U \to V\).  Our goal is to find the solution \(u \in U\) that satisfies, for a given \(f\in V\), the partial differential equation cast here in its residual form:
\begin{equation}
\label{eq:residual}
R \big(\bx, u(\bx) \big) := f(\bx) - Au(\bx) =  0 , \quad \forall \bx \in \Omega,
\end{equation} 
and the following boundary conditions:
\begin{equation}
\label{eq:bc}
B \big(\bx, u(\bx)\big) = 0 , \quad \forall \bx \in \partial \Omega.
\end{equation} 

For the sake of simplicity in the presentation, but without loss of generality, we consider here only the case of homogeneous Dirichlet boundary conditions, such that the residual $B$ is given by
\begin{equation}
B\big(\bx,u (\bx)\big) := u(\bx), \quad \forall \bx \in \partial \Omega.
\end{equation} 
%{for some underlying restriction operator.}
The primary objective in PINNs is to use a neural network with parameters \(\theta\) to find an approximation \(\tilde{u}_\theta(\bx)\) of the solution \(u(\bx)\) to problem~\eqref{eq:residual}-\eqref{eq:bc}. For the sake of simplicity in the notation, we shall omit in the rest of the paper the subscript~\(\theta\) when referring to the approximate solutions \(\tilde{u}_\theta\), and thus simply write \(\tilde{u}(\bx)\). The training, i.e.\ the identification of the parameters $\theta$ of the neural network, is performed by minimizing a loss function, defined here as a combination of the residual associated with the partial differential equation and that associated with the boundary condition in terms of the \(L^2\) norm:
\begin{equation}
\label{eq:losspinn}
\mathcal L(\theta) 
:= w_r \int_\Omega R\big(\bx,\tilde{u}(\bx) \big)^2 dx
+ w_{bc} \int_{\partial \Omega} B\big(\bx, \tilde{u}(\bx)\big)^2 dx,
\end{equation}
where \(w_r\) and \(w_{bc}\) are penalty parameters. In other words, by minimizing the loss function~\eqref{eq:losspinn} one obtains a weak solution \(\tilde{u}\) that weakly satisfies the boundary condition. 

Alternatively, the homogeneous Dirichlet boundary condition could be strongly imposed, as done in~\cite{mcfall2009artificial}, by multiplying the output of the neural network by a function \(g(\bx)\) that vanishes on the boundary. For instance, if \(\Omega = (0,\ell) \in \R\), one could  choose \(g(x) = x(\ell-x)\). The trial functions \( \tilde{u} \) would then be constructed, using the feedforward neural network~\eqref{eqn:FNN}, as follows:
\begin{equation}
\label{eqn:FNN_g}
\begin{aligned}
&\text{Input layer:} &&  \bz_0 = \bx, \\
&\text{Hidden layers:} &&  \bz_{i} = \sigma ( W_i \bz_{i-1} + \bb_i), 
\quad i=1,\ldots,n, \\ 	
&\text{Output layer:} && z_{n+1} = W_{n+1} \bz_{n} + \bb_{n+1} ,\\
&\text{Trial function:} && \tilde{u} = g(\bx) z_{n+1}.
\end{aligned}
\end{equation}
where the input and output layers have a width \(N_0 = d\) and \(N_{n+1} = 1\), respectively. The dimension of the finite-dimensional space of functions generated by the neural network~\eqref{eqn:FNN_g} is now given by \(N_\theta=\sum_{i=1}^{n+1} N_{i} (N_{i-1} + 1) = N_1(d+1) + \sum_{i=2}^{n} N_{i} (N_{i-1} +1) + (N_n + 1)\).

For the rest of this work, the boundary conditions will be strongly imposed, so that the loss function will henceforth be
\begin{equation}
\label{eq:losspinn_strong_R}
\mathcal L(\theta) 
= \int_\Omega R\big(\bx,\tilde{u}(\bx) \big)^2 dx.
\end{equation}
The problem that one solves by PINNs can thus be formulated as:
\begin{equation}
\min_{\theta\in \mathbb R^{N_\theta}} \mathcal L(\theta) = \min_{\theta\in \mathbb R^{N_\theta}} \int_\Omega R\big(\bx,\tilde{u}(\bx) \big)^2 dx.
\end{equation}

One advantage of PINNs is that they do not necessarily need the construction of a mesh, which is often a time-consuming process. Instead, the integral in the loss function can be approximated using Monte Carlo integration from randomly generated points in $\Omega$. Another advantage is the ease of implementation of the boundary and initial conditions. On the other hand, one major issue that one faces when using PINNs is that it is very difficult, even impossible, to effectively reduce the $L^2$ or $H^1$ error in the solutions to machine precision. The main reason, from our own experience, is that the solution process may get trapped in some local minima, without being able to converge to the global minimum, when using non-convex optimization algorithms. 
%{Another way to state this is that in a local minimum, the solution becomes very sensitive to small changes in the weights and biases, thereby limiting the effective accuracy that can be reached.} {"Note: But it can also be the opposite where the solution is not sensitive at all. Should we omit this sentence?"} 
We briefly review some commonly used optimizers and study their performance in the next section. 

%============================================================
\subsection{Choice of the optimization algorithm}
\label{sec:optimizer}

The objective functions in PINNs are by nature non-convex, which makes the minimization problems difficult to solve and their solutions highly dependent on the choice of the solver. For these reasons, it is common practice to employ gradient-based methods, such as the Adam optimizer~\cite{DBLP:journals/corr/KingmaB14} or the Broyden–Fletcher–Goldfarb–Shanno (BFGS) algorithm~\cite{fletcher2000practical}. BFGS is a second-order optimizer, but if used alone, has the tendency to converge to a local minimum in the early stages of the training. A widely used strategy to overcome this deficiency is to begin the optimization process using the Adam optimizer and subsequently switch to the BFGS optimizer~\cite{markidis2021old}. In this work, we will actually utilize the so-called L-BFGS optimizer, the limited-memory version of BFGS provided in PyTorch~\cite{NEURIPS2019_9015}. Although L-BFGS is a higher-order method than Adam, the computational cost for each iteration is also much higher than the cost for one iteration of Adam. We actually adopt here the following definition of what we mean by an iteration: in both algorithms, it actually corresponds to a single update of the neural network parameters.

In the following example, we study the performance of the aforementioned strategy, when applied to a simple one-dimensional Poisson problem, and compare the resulting solution with that obtained when using the Adam optimizer only. This numerical example will also serve later as a model problem for further verifications of the underlying principles in our approach.

\begin{example}
\label{example:optimizers}
Given a function \(f(x)\), the problem consists in finding \(u=u(x)\), for all \(x\in[0,1]\), that satisfies 
\begin{equation}
\label{eq:poisson_problem}
\begin{aligned}
-\partial_{xx} u(x) &= f(x), \qquad \forall x\in(0,1), \\
u(0) &=0, \\
u(1) &=0.
\end{aligned}
\end{equation}
For the purpose of the study, the source term \(f\) is chosen such that the exact solution to the problem is given as 
\begin{equation}
\label{eq:solution_exact}
u(x) = e^{\sin(k\pi x)} +x^3 - x -1,   
\end{equation}
where \(k\) is a given integer. We take \(k = 2\) in this example. 

We consider here a network made of only one hidden layer of a width of \(20\), i.e.\ \(n=1\) and \(N_1=20\). Moreover, \(N_0=N_2=1\). The learning rates for the Adam optimizer and L-BFGS are set to \(10^{-2}\) and unity, respectively. In the first experiment, the network is trained for 10,000 iterations using Adam. In the second experiment, it is trained with Adam for 4,000 iterations followed by 100 iterations of L-BFGS. Figure~\ref{fig:prob_optimizer} compares the evolution of the loss function with respect to the number of iterations for these two scenarios. In the first case, we observe that the loss function laboriously reaches a value around \(10^{-2}\) after 10,000 iterations. The loss function further decreases in the second case but still plateaus around \(5 \times 10^{-5}\) after about 30 iterations of L-BFGS. Note that the scale along the \(x\)-axis in the figure on the right has been adjusted in order to account for the large discrepancy in the number of iterations used with Adam and L-BFGS.

\begin{figure}[tb]
\centering
\includegraphics[width=0.32\linewidth]{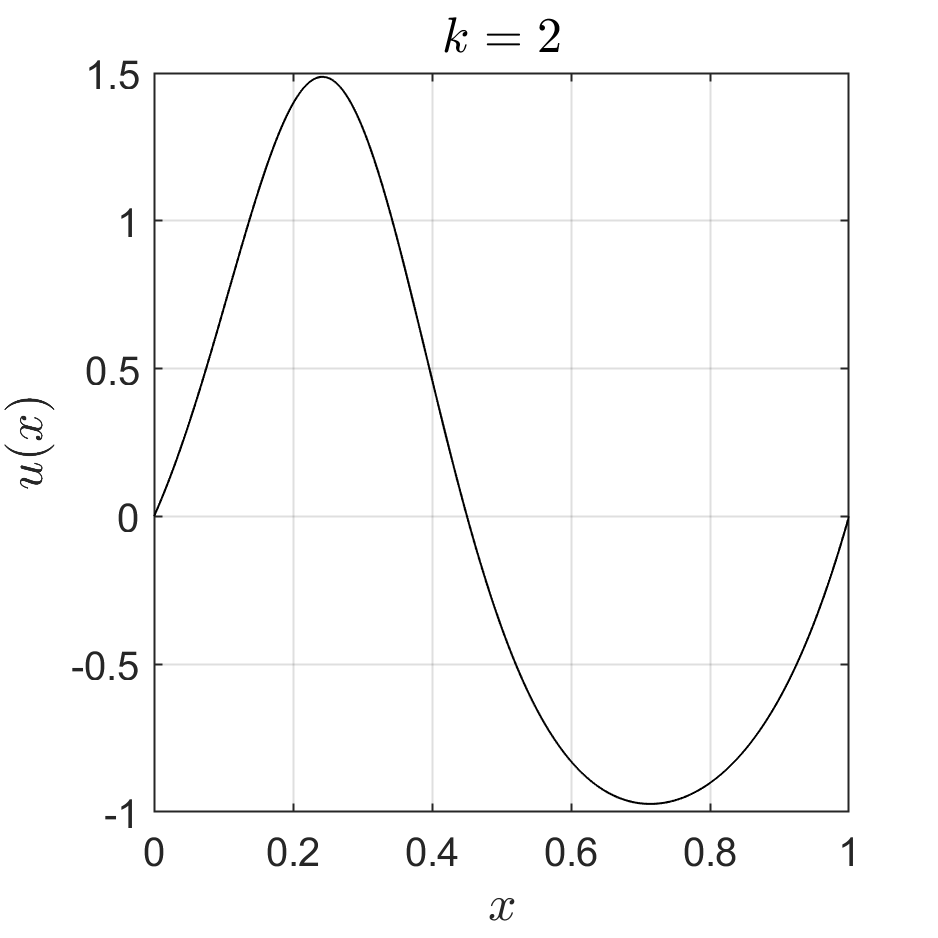}
\includegraphics[width=0.32\linewidth]{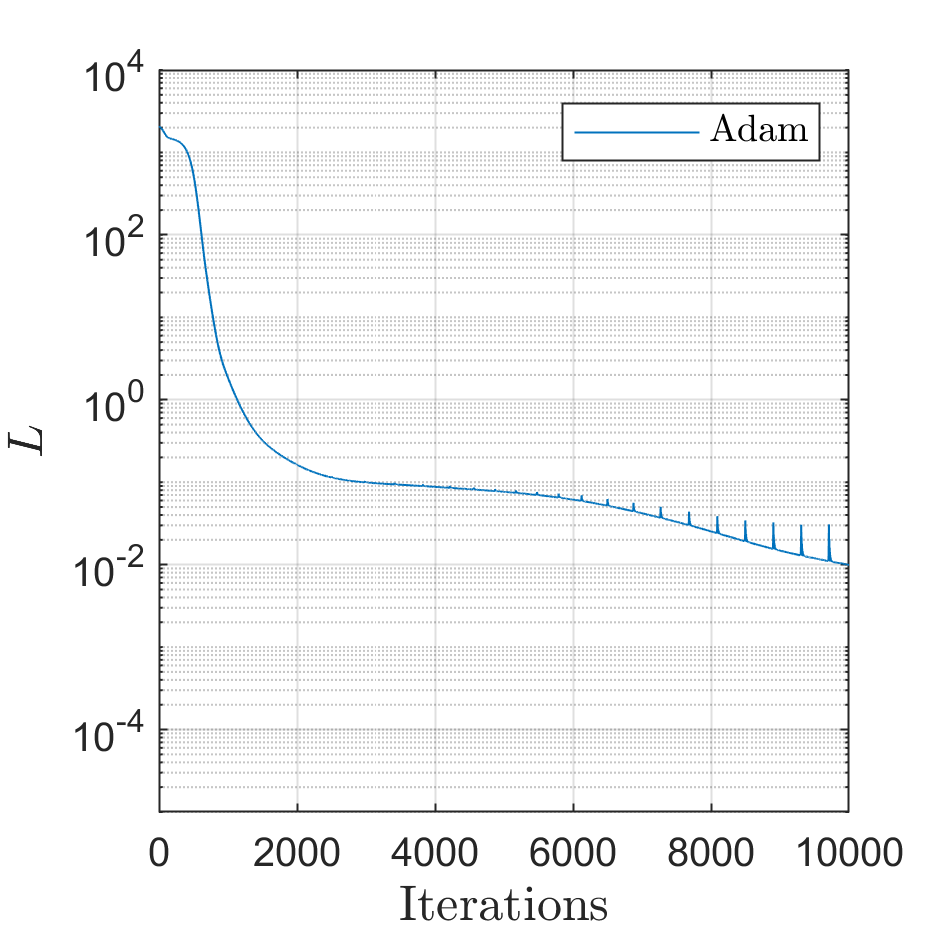}
\includegraphics[width=0.32\linewidth]{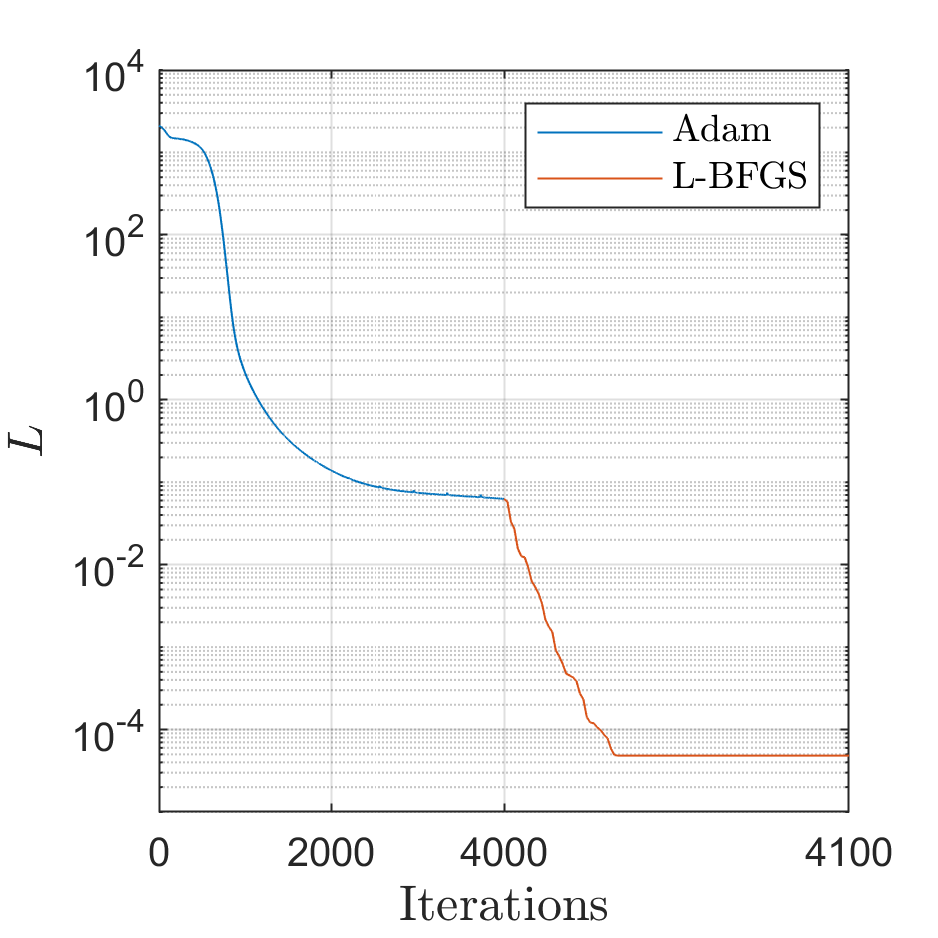}
\caption{Results from Example~\ref{example:optimizers} in Section~\ref{sec:optimizer}. (Left) Exact solution with \(k = 2\). (Middle) Evolution of the loss function using the Adam optimizer only. (Right) Evolution of the loss function using the Adam optimizer and L-BFGS.}
\label{fig:prob_optimizer}
\end{figure}

\end{example}
 
%============================================================
\section{Error analysis in PINNs}
\label{sec:error}

In this section, we further study the numerical errors, and a fortiori, the sources of error, in the solutions obtained with PINNs. We have mentioned in the introduction that two issues may actually affect the quality of the solutions. Indeed, it is well known that the training of the neural networks may perform poorly if the data, in our case the source term in the differential equations, are not properly normalized~\cite{goodfellow2016deep}. Moreover, the accuracy may deteriorate when the solutions to the problem exhibit high frequencies. We briefly review here the state-of-the-art in dealing with those two issues as they will be of paramount importance in the development of the multi-level neural network approach. More specifically, we illustrate on simple numerical examples how these issues can be somewhat mitigated.

%============================================================
\subsection{Data normalization}
\label{sec:normalize}

A major issue when solving a boundary-value problem such as~\eqref{eq:residual}-\eqref{eq:bc} with PINNs is the amplitude of the problem data, in particular, the size of the source term~\(f(x)\). In other words, a small source term naturally implies that the target solution will be also small, making it harder for the training to find an accurate approximation of the solution to the boundary-value problem. This issue will become very relevant and crucial when we design the multi-level neural network approach in  Section~\ref{sec:multi-nets}. Our goal here is to illustrate through a numerical example that the accuracy of the solution clearly depends on the amplitude of the data, and hence of the solution itself, and that it may therefore be necessary to scale the solution before minimizing the cost functional. We hence revisit Example~\ref{example:optimizers} of Section~\ref{sec:optimizer}.

\begin{example}
\label{example:normalization}
We solve again the Poisson problem~\eqref{eq:poisson_problem} in \(\Omega = (0,1)\) with \(k=2\). However, we deliberately divide the source term \(f(x)\) by a factor \(\mu\) such that the exact solution is changed to
\[
u(x) = \frac{1}{\mu} \big(e^{\sin(k\pi x)} +x^3 - x -1 \big).
\]
A large value of \(\mu\) implies a small \(f\), and hence, a small \(u\). We now compare the solutions of the problem for several values of \(\mu\) with different orders of magnitude, namely \(\mu = \{10^{-3},1,10^3\}\). For the training, we consider the Adam optimizer followed by L-BFGS using the same network architecture and hyper-parameters as in Section~\ref{sec:optimizer}.

\begin{figure}[tb]
\centering
\includegraphics[width=0.32\linewidth]{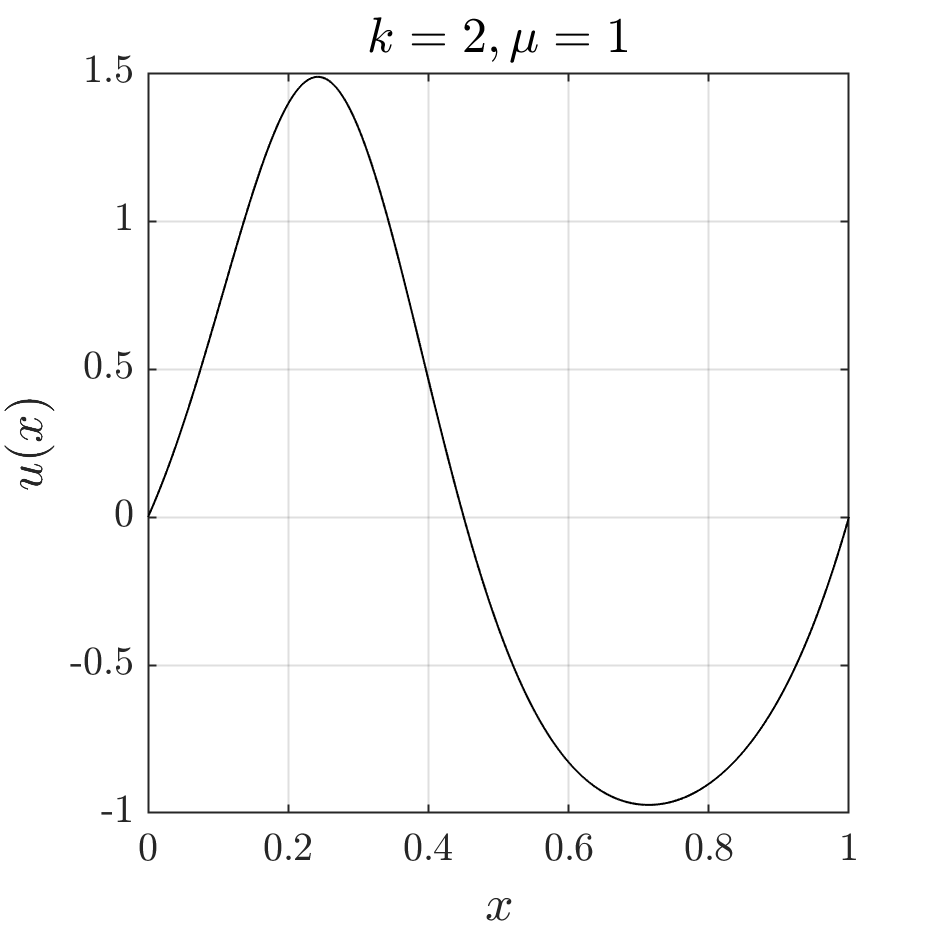}
\includegraphics[width=0.32\linewidth]{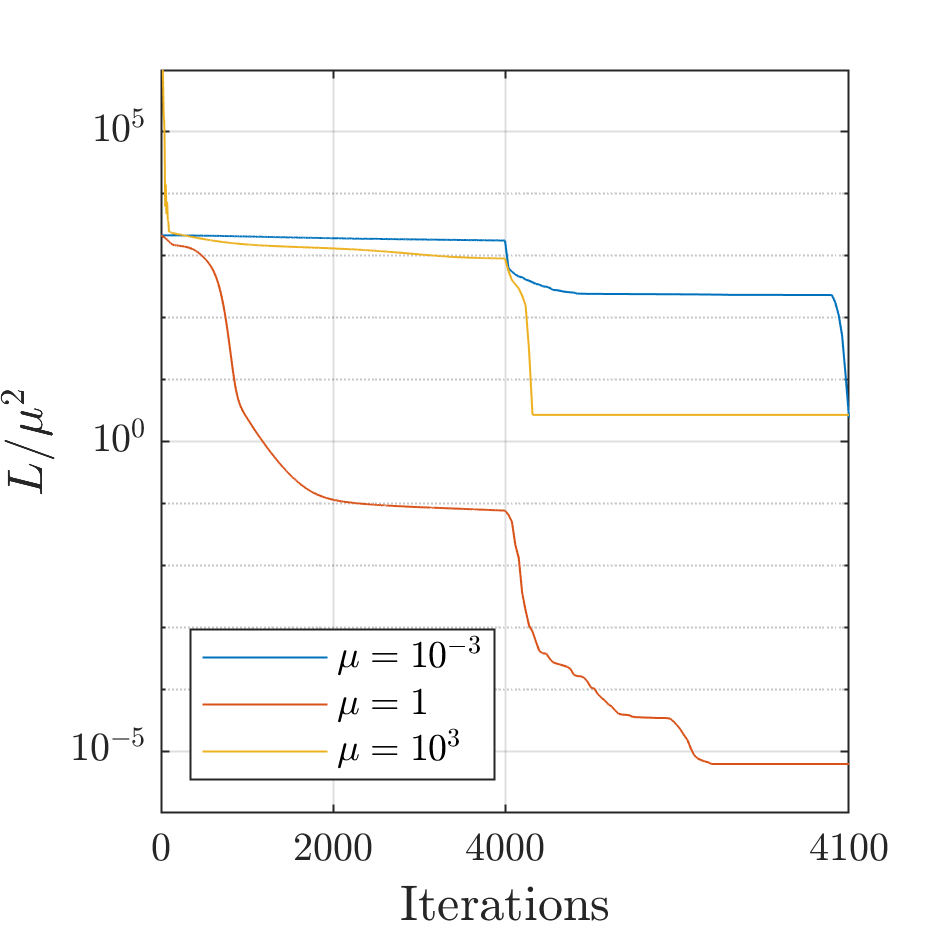}
\includegraphics[width=0.32\linewidth]{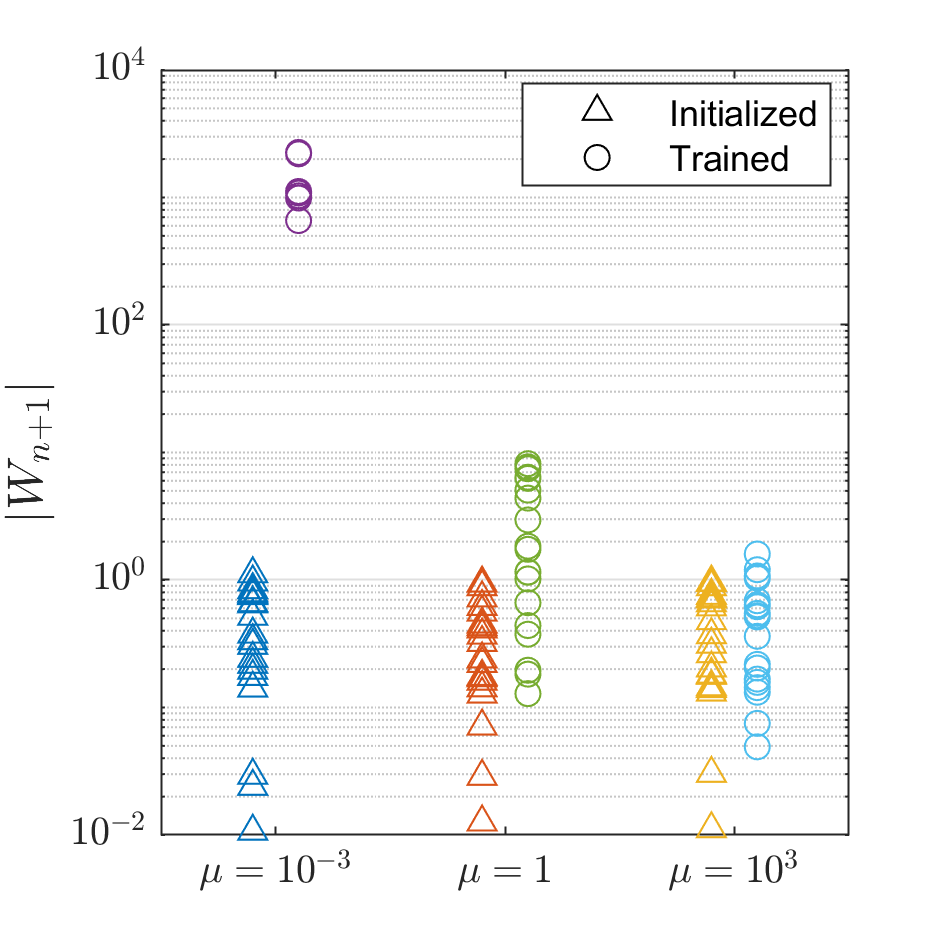}
\caption{Results from Example~\ref{example:normalization} in Section~\ref{sec:normalize}: (Left) Exact solution with \(k = 2\) and \(\mu=1\). (Middle) Evolution of the loss function for \(\mu=10^{-3}\), \(1\), and \(10^3\). (Right) Distribution of the absolute value of the weights in the last layer before and after training. 
%For visibility, the two columns at \(\mu=10^{-3}\), \(1\), and \(10^3\) are slightly offset.
}
\label{fig:normalize}
\end{figure}

\begin{figure}[tb]
\centering
\includegraphics[width=0.32\linewidth]{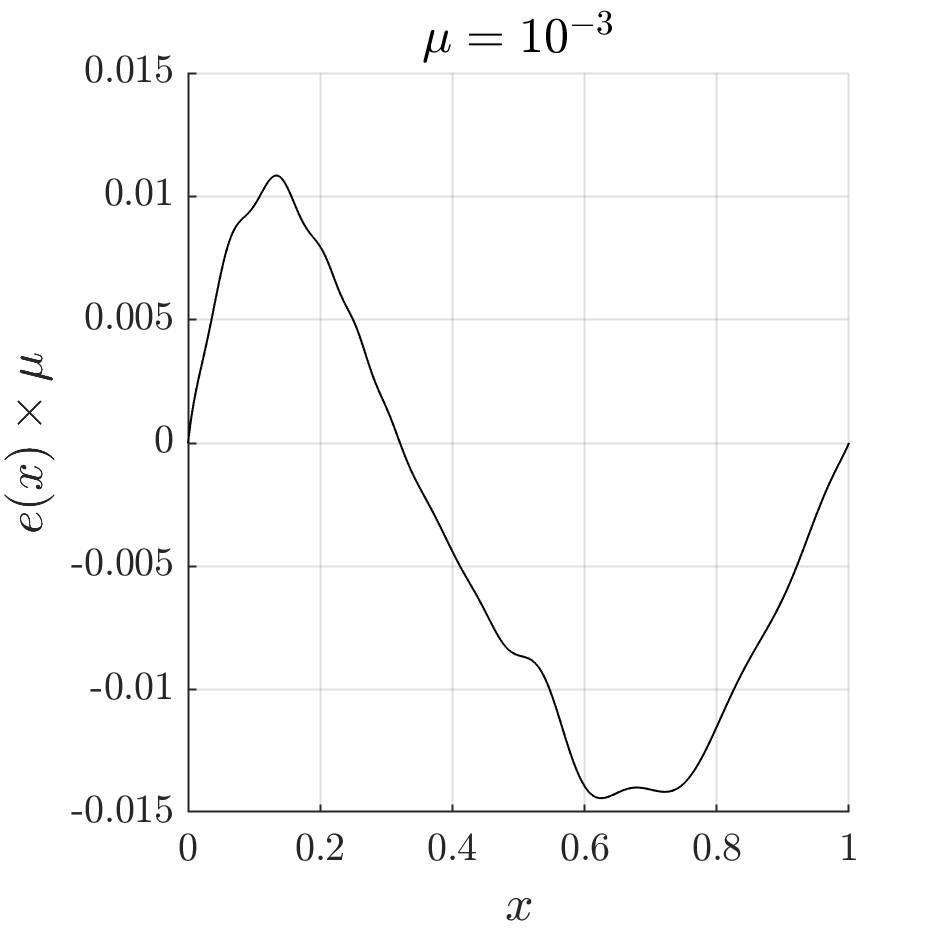}
\includegraphics[width=0.32\linewidth]{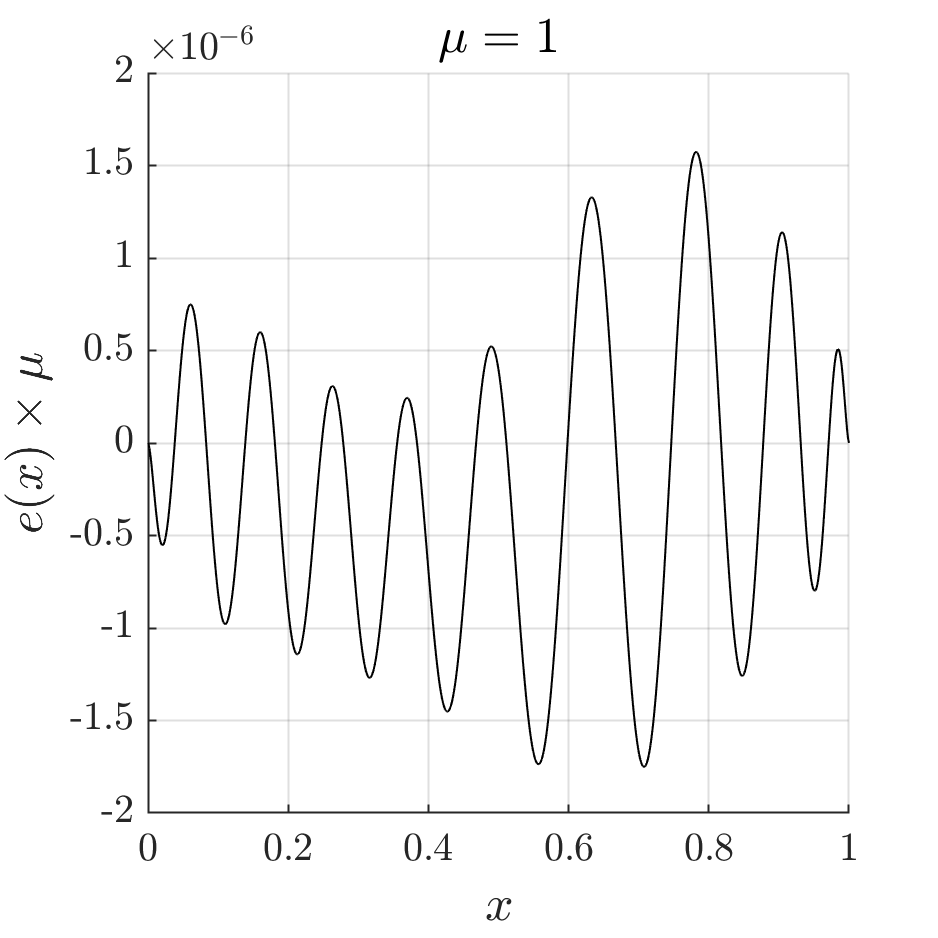}
\includegraphics[width=0.32\linewidth]{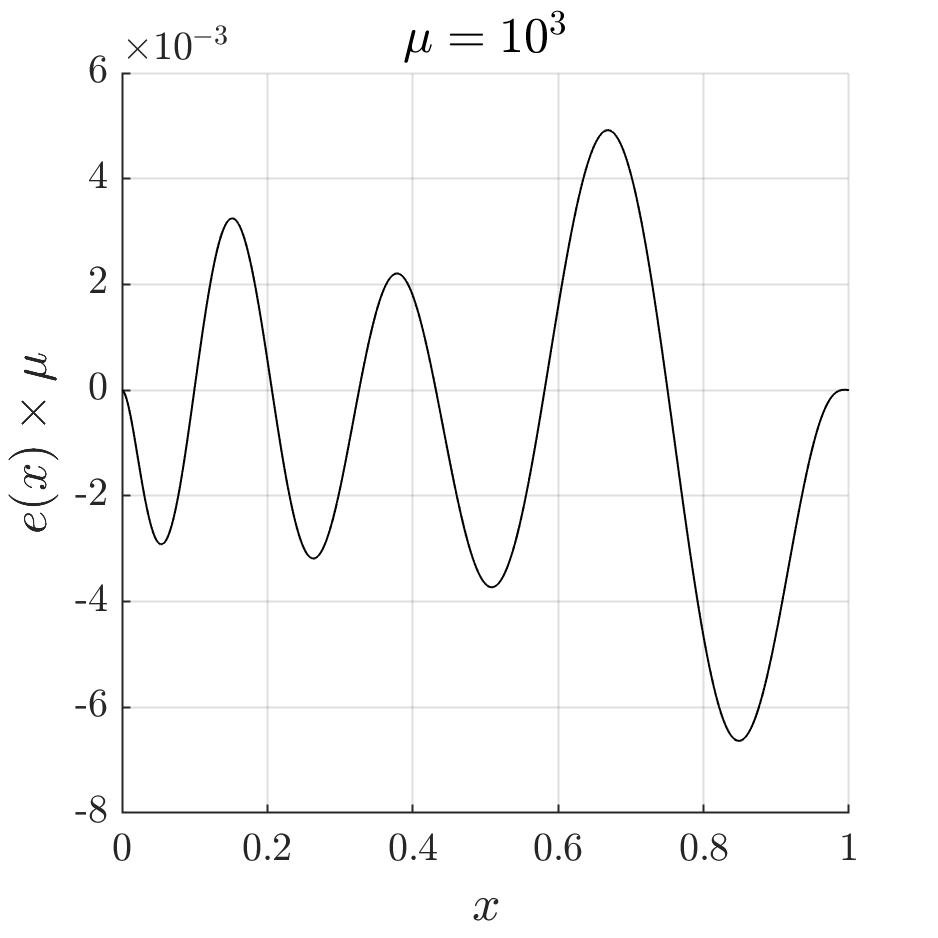}
\caption{Results from Example~\ref{example:normalization} in Section~\ref{sec:normalize}: Pointwise error \(e(x) =  u(x) - \tilde{u}(x)\) for \(\mu=10^{-3}\), \(1\), and \(10^3\).}
\label{fig:normalize-error}
\end{figure}

We show in Figure~\ref{fig:normalize} (left) the solution for \(k=2\) and \(\mu=1\). We want to draw attention to the fact that the maximal amplitude of this solution is roughly unity.

The evolution of the loss function during training is shown in Figure~\ref{fig:normalize} (middle) for the three values of~\(\mu\). We actually plot each loss function as computed but divided by \(\mu^2\) for a clearer comparison. For \(\mu = 1\), we observe that the loss function converges much faster and achieves a much smaller residual at the end of the training, than for \(\mu = 10^3\) and \(\mu=10^{-3}\). We show in Figure~\ref{fig:normalize-error} the errors in the three solutions obtained after training. The error in the solution computed with \(\mu=1\) is indeed several orders of magnitude smaller than the error in the other two solutions. An important observation from these plots is that smaller errors induce higher frequencies. Hence, if one wants to reduce the error even further, it would become necessary to have an algorithm that allows one to capture those higher frequencies. This issue is addressed in the next section.

We remark that the solution obtained with \(\mu=1\) would actually provide a more accurate approximation by simply multiplying it by \(\mu = 10^3\) (resp.\ \(\mu=10^{-3}\)) to the problem with \(\mu = 10^3\) (resp.\ \(\mu=10^{-3}\)). In other words, it illustrates the fact that, when using PINNs, the process of multiplying the source term by \(\mu\), solving, and then re-scaling the solution by \(\mu^{-1}\) is simply not equivalent to the process of simply solving the problem.

The distribution of the weights in the output layer \(|\bW_{n+1} |\) obtained after  initialization and training is shown in Figure~\ref{fig:normalize} (right) for each \(\mu\). First, we observe that the final weights for \(\mu =10^{-3}\) are very different from their initial values and sometimes exceed \(10^{3}\). Second, we would expect the trained parameters for \(\mu =10^{3}\) to be three orders of magnitude smaller than those for \(\mu=1\). However, it seems that the network has difficulty decreasing the values of these weights. As a consequence, the training fails to properly converge in the two cases \(\mu=10^3\) and \(\mu=10^{-3}\). A reasonable explanation is that an accurate solution cannot be obtained if the optimal weights exist far from their initialized values, since in this case the training of the network is more demanding. This implies that, for a very small or very large value of \(\mu\), an efficient  initialization will not suffice to improve the training. One could perhaps adjust the learning rate for the last layer, but finding the proper value of the learning rate is far from a straightforward task. Therefore, a simpler approach would be to normalize the solution being sought so that the output of the neural network is largely of the order of unity. We will propose such an approach in Section~\ref{sec:multi-nets}.
\end{example}

%============================================================
\subsection{Solutions with high frequencies}
\label{sec:high_freq}

A deep neural network usually adheres to the F-principle~\cite{rahaman2019spectral, ronen2019convergence, xu2019frequency}, which states that the neural network tends to approximate the low-frequency components of a function before its high-frequency components. This property explains why networks approximate well functions featuring a low-frequency spectrum while avoiding aliasing, leading to reasonable generalized errors. The F-principle also serves as a filter for noisy data and provides an early stopping criterion to avoid overfitting. When it comes to handling higher frequencies, one is generally exposed to the risk of overfitting and the lack of convexity of the loss function. Unfortunately, there exist few guidelines, to the best of our knowledge, to ensure that the training yields accurate solutions in those cases. As is often the case with PINNs, the quality of the obtained solutions depends on the experience of the user with the initialization of the hyper-parameters.

Several studies, see e.g.~\cite{sitzmann2020implicit, liu2020multi,mildenhall2021nerf}, have put forward some techniques to improve neural networks in approximating high-frequency functions. We start by providing a concise overview of the Fourier feature mapping presented in~\cite{mildenhall2021nerf}, which we shall use in this work, and proceed with an illustration of its performance on a simple one-dimensional example.

In order to  simultaneously approximate the low and high frequencies, the main idea behind the method is to explicitly introduce within the networks high-frequency modes using the so-called Fourier feature mapping. Let \(\boldsymbol{\omega}_M\) denote the vector of \(M\) given wave numbers \(\omega_m\), \(m=1,\ldots,M\), that is $\boldsymbol{\omega}_M = [\omega_1, \ldots, \omega_M]$. The mapping \(\gamma\) for each spatial component \(x_j\) is provided by the row vector of size \(2M\) defined as:
\begin{equation}
\gamma(x_j) = [\cos(\boldsymbol{\omega}_M x_j), \sin(\boldsymbol{\omega}_M x_j)], \qquad j=1,\ldots,d,
\end{equation}
where we have used the shorthand:
\[
\begin{aligned}
&\cos(\boldsymbol{\omega}_M x_j) 
= [\cos(\omega_1 x_j),\cos(\omega_2 x_j),\ldots,\cos(\omega_M x_j)], \\
&\sin(\boldsymbol{\omega}_M x_j) 
= [\sin(\omega_1 x_j),\sin(\omega_2 x_j),\ldots,\sin(\omega_M x_j)].
\end{aligned}
\]
As shown with the Neural Tangent Kernel theory in~\cite{tancik2020fourier}, the Fourier feature mapping helps the network learn the high and low frequencies simultaneously.
The structure of the feedforward neural network~\eqref{eqn:FNN_g} is now modified as follows. Considering a network with an input layer of width \(N_0 = 2M \times d\) and an output layer of width \(N_{n+1} = 1\), the trial functions \( \tilde{u} \) are taken in the form:
\begin{equation}
\label{eqn:FNN_hf1}
\begin{aligned}
&\text{Input layer:} && \bz_0 = [ \gamma(x_1),\ldots, \gamma(x_d) ]^T, \\
&\text{Hidden layers:} &&  \bz_{i} = \sigma ( W_i \bz_{i-1} + \bb_i), 
\quad i=1,\ldots,n, \\ 	
&\text{Output layer:} && z_{n+1} = W_{n+1} \bz_{n} + \bb_{n+1} ,\\
&\text{Trial function:} && \tilde{u} = g(\bx) z_{n+1}.
\end{aligned}
\end{equation}
The dimension of the finite-dimensional space of trial functions is given in this case by \(N_\theta = N_1(2Md + 1) + \sum_{i=2}^{n} N_{i} (N_{i-1} +1) + (N_n + 1)\).

In a similar manner, we will show on a numerical example that using a function \(g(\bx)\) whose spectrum contains both low and high frequencies also improves the convergence of the solutions. In the present work, we only consider one-dimensional problems or two-dimensional problems on rectangular domains so that one can introduce a new mapping \(\gamma_g\) in terms of only the sine functions and thus strongly impose the boundary conditions by:
\begin{equation}
\label{eq:gamma_g} 
\gamma_g(x_j) = [ \sin( \boldsymbol{\omega}_M  x_j)], \qquad j=1,\ldots,d.
\end{equation}
We note here that the wave number vector \(\boldsymbol{\omega}_M \) is the same as in $\gamma$ and should be chosen such that all sine functions vanish on the boundary $\partial\Omega$.
In that case, we consider a feedforward neural network with an input layer of width \(N_0 = 2M\times d\) and an output layer of width \(N_{n+1} = M\), so that the trial functions \( \tilde{u} \) are given by:
\begin{equation}
\label{eqn:FNN_hf2}
\begin{aligned}
&\text{Input layer:} && \bz_0 = [ \gamma(x_1),\ldots, \gamma(x_d) ]^T, \\
&\text{Hidden layers:} &&  \bz_{i} = \sigma ( W_i \bz_{i-1} + \bb_i), 
\quad i=1,\ldots,n, \\ 	
&\text{Output layer:} && \bz_{n+1} = W_{n+1} \bz_{n} + \bb_{n+1} ,\\
&\text{Trial function:} && \tilde{u} = M^{-1}\big( \Pi_{j=1}^d\gamma_g(x_j)\big) \cdot \bz_{n+1},
\end{aligned}
\end{equation}
where the trial function is divided by \(M\) in order to normalize the output. The dimension of the finite-dimensional space of trial functions generated by the neural network~\eqref{eqn:FNN_g} is now given by \(N_\theta = N_1(2Md+1) + \sum_{i=2}^{n} N_{i} (N_{i-1} +1) + M(N_n+1)\).
We reiterate here that the output \(\bz_{n+1}\) needs to be multiplied by sine functions that vanish on the boundary in order to strongly impose the boundary condition. When \(\Omega=(0,\ell)^d \), an appropriate choice for the parameters \(\omega_m\) is given by the geometric series \( \omega_m = 2^{m-1}\pi/\ell\), with \(m=1,\ldots,M\), as suggested in~\cite{mildenhall2021nerf}. 

We now compare the performance of the three approaches using Example~\ref{example:optimizers}, in order to show the importance of introducing high frequencies in the input layer and in the function used to enforce the boundary conditions. The three methods can be summarized as follows:
\begin{itemize}[itemsep=0pt,topsep=5pt,parsep=0pt,leftmargin=20pt]
    \item \textbf{Method 1:} Classical PINNs with input \(x\) and trial functions provided by the neural network~\eqref{eqn:FNN_g} with \(g(x)= x(1-x)\).
    \item \textbf{Method 2:} PINNs using the Fourier feature mapping for the input and trial functions provided by the neural network~\eqref{eqn:FNN_hf1} with \(g(x)= x(1-x)\).
    \item \textbf{Method 3:} PINNs using the Fourier feature mapping for the input and trial functions provided by the neural network~\eqref{eqn:FNN_hf2}.
\end{itemize}

\begin{example}
\label{example:hf}
We solve the Poisson problem~\eqref{eq:poisson_problem} in \(\Omega=(0,1)\) with \(k=10\). The exact solution is given in~\eqref{eq:solution_exact} and shown in Figure~\ref{fig:high_freq} (left). The networks all have a single hidden layer of width \(N_1=10\). As before, the learning rates for Adam and L-BFGS are chosen as \(10^{-2}\) and unity, respectively. The training is performed for 4,000 iterations with Adam and 100 iterations with L-BFGS.
The vector \(\boldsymbol{\omega}_M\) of wave numbers $\omega_m$, $m=1,\ldots,M$, is constructed from a geometric series with \(M = 4\), i.e.\ \(\boldsymbol{\omega}_M = [\pi,2\pi,4\pi,8\pi] \). 

\begin{figure}[tb]
\centering
\includegraphics[width=0.32\linewidth]{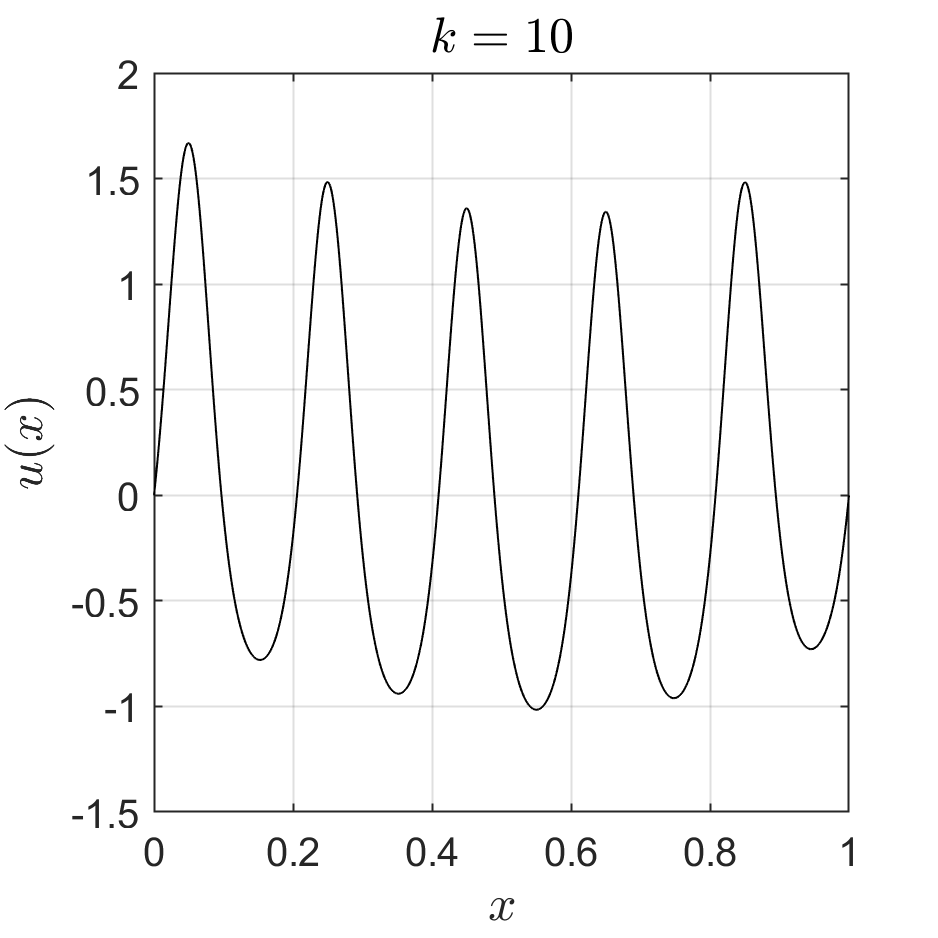}
\includegraphics[width=0.32\linewidth]{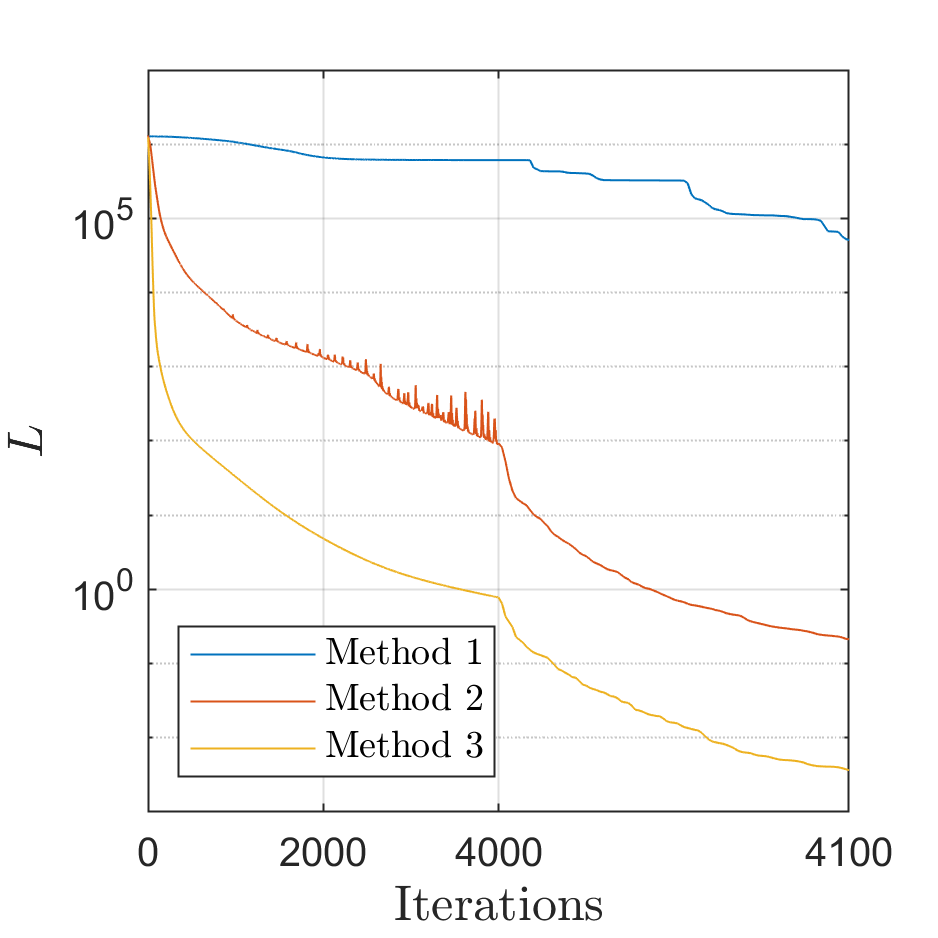}
\includegraphics[width=0.32\linewidth]{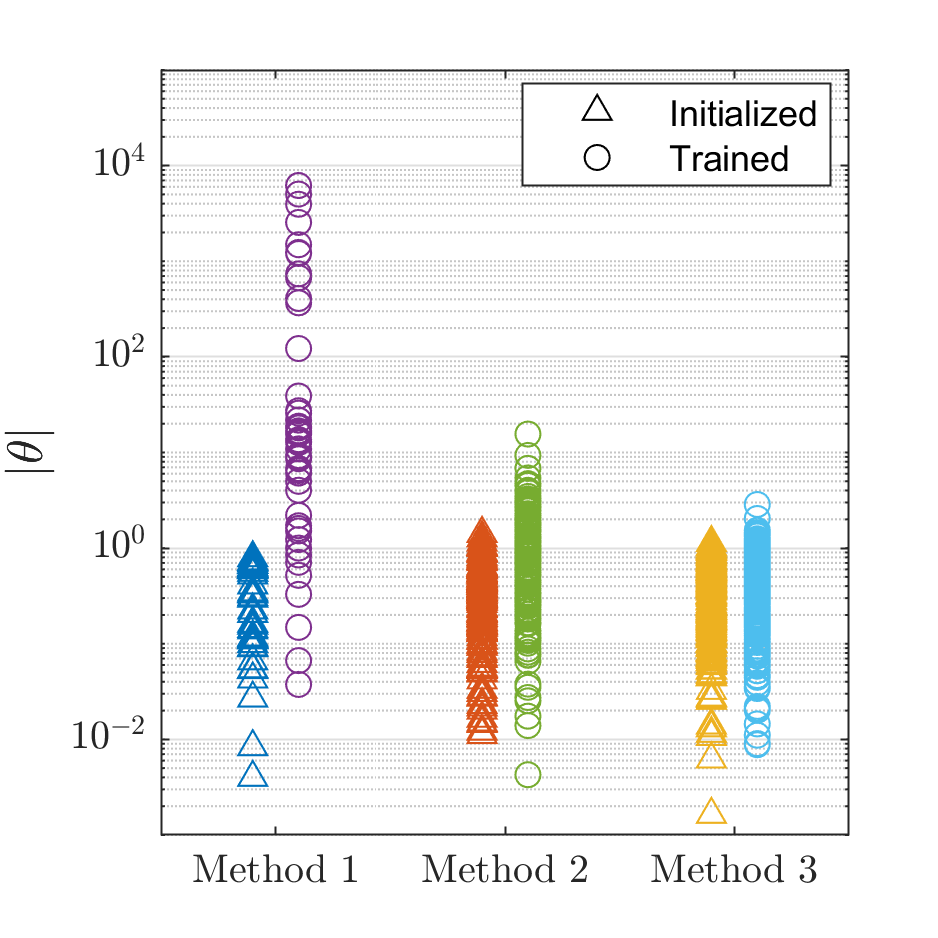}
\caption{Results from Example~\ref{example:hf} in Section~\ref{sec:high_freq}: (Left) Exact solution with \(k = 10\). (Middle) Evolution of the loss function for the three methods. (Right) Distribution of the absolute value of the initialized and trained parameters for the three methods.}
\label{fig:high_freq}
\end{figure}

We first observe in Figure~\ref{fig:high_freq} (middle) that Method~1 fails to converge. On the other hand, Method~2 allows one to decrease the loss function by six orders of magnitude and Method~3 reduces further the loss function by almost two orders of magnitude. This example indicates that it is best to use the architecture given in~\eqref{eqn:FNN_hf2} when dealing with solutions with high frequencies. 

We show in Figure~\ref{fig:high_freq} (right) the absolute value of the initialized and trained parameters for each method. We notice that the trained parameters are very large in the case of the first method, with some of them reaching values as large as \(10^4\). In contrast, the values remain much smaller in the case of Methods~2 and~3, with the parameters of Method~3 staying closer to the initialized parameters when compared to those obtained by Method~2. For Method 1, the weights in the hidden layers need to be large so as to be able to capture the high frequencies, as seen in Figure~\ref{fig:high_freq} (right). If one uses  Method~2 to obtain an approximation of the exact solution~\eqref{eq:solution_exact}, the function computed by the output layer in~\eqref{eqn:FNN_hf1} should converge to the function:
\[
\dfrac{u(x)}{x(1-x)} = \dfrac{e^{\sin(k\pi x)} +x^3 - x -1}{x(1-x)}.
\]
However, this function takes on large values near the boundary. Indeed, when \(x\) tends to \( 0\), the limit is equal to \(k \pi - 1\), which becomes large for large values of \(k\). Hence, the parameters of the network after training will tend to take large values in order to approximate well the solution, as explained in Section~\ref{sec:normalize}. In order to avoid these issues, we have thus introduced the architecture~\eqref{eqn:FNN_hf2}, such that the  functions used to enforce the boundary conditions contain a mix of low and high frequencies. Method~3 thus allows one to get a solution whose trained parameters remain of the same order as the initial ones, as observed in Figure~\ref{fig:high_freq} (right).

\begin{figure}[tb]
\centering
\includegraphics[width=0.32\linewidth]{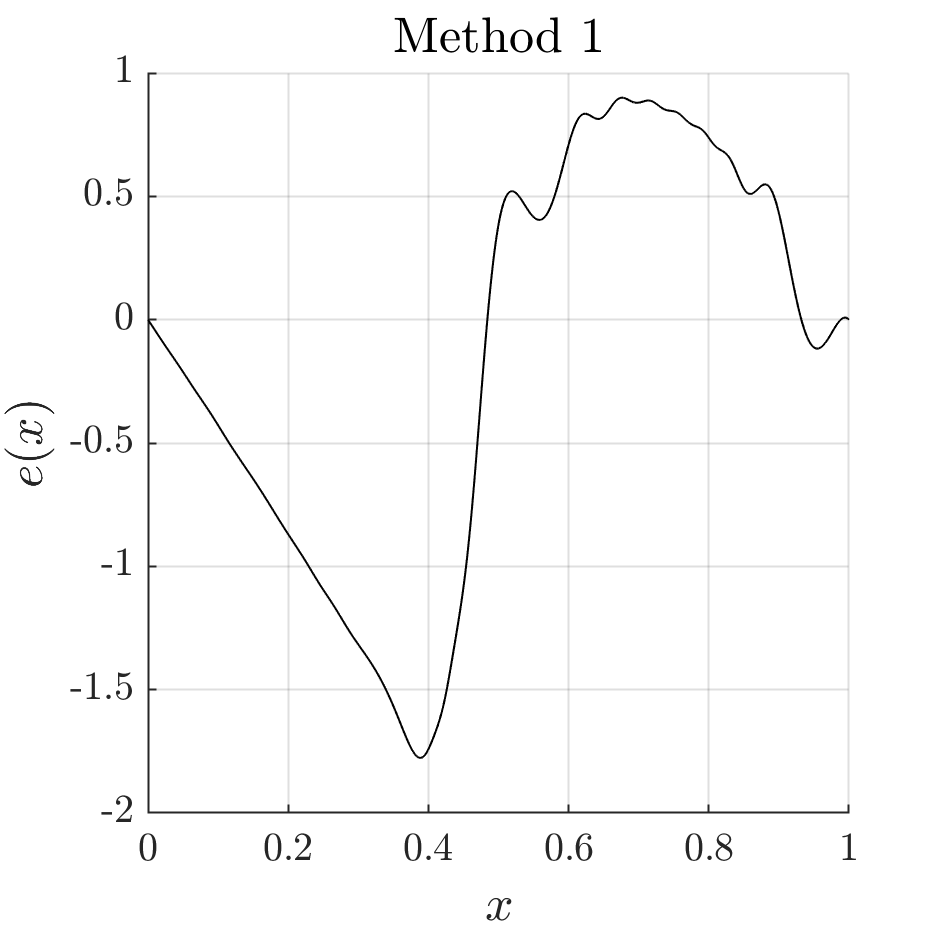}
\includegraphics[width=0.32\linewidth]{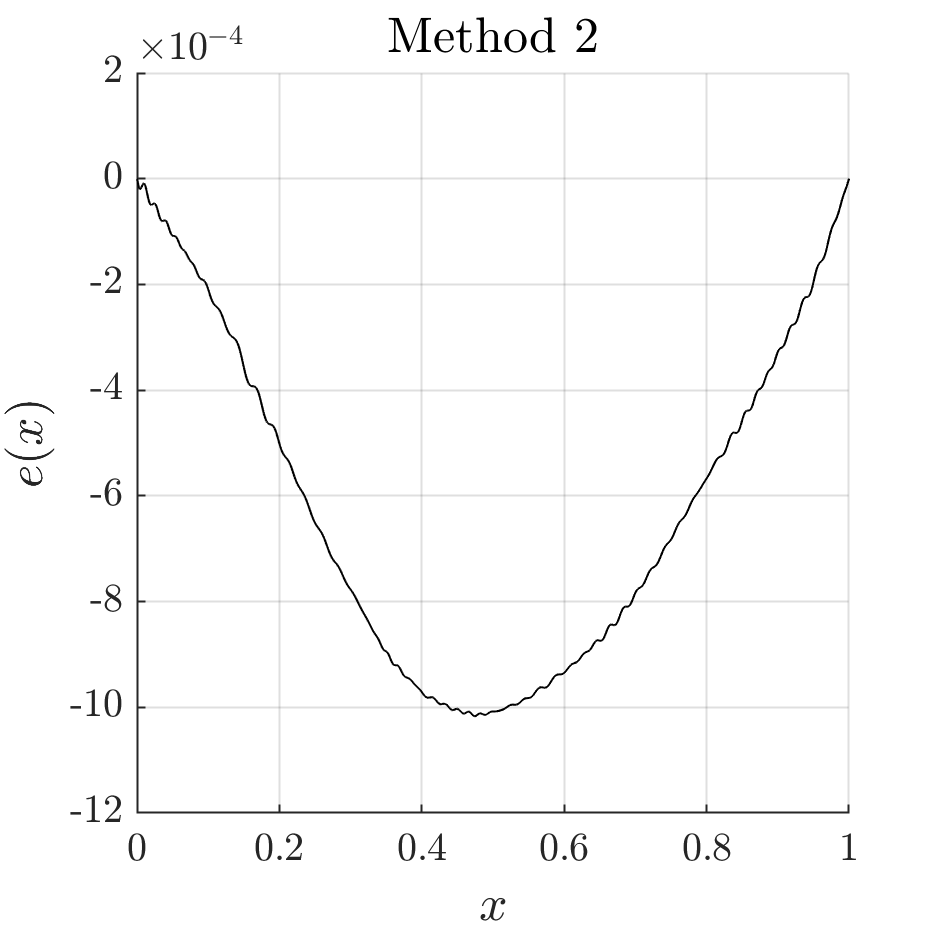}
\includegraphics[width=0.32\linewidth]{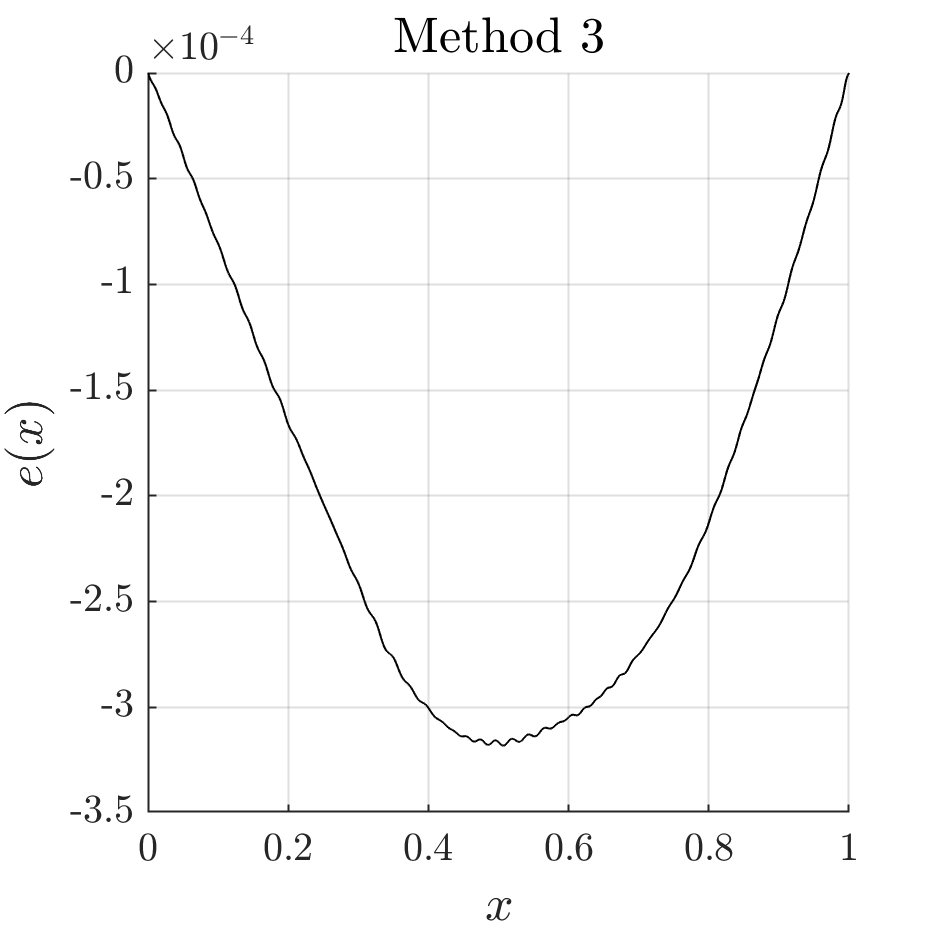}
\caption{Results from Example~\ref{example:hf} in Section~\ref{sec:high_freq}: Pointwise error \(e(x) = u(x) - \tilde{u}(x) \) for Methods 1, 2, and 3. }
\label{fig:high_freq_error}
\end{figure}

\begin{figure}[tb]
\centering
\includegraphics[width=0.32\linewidth]{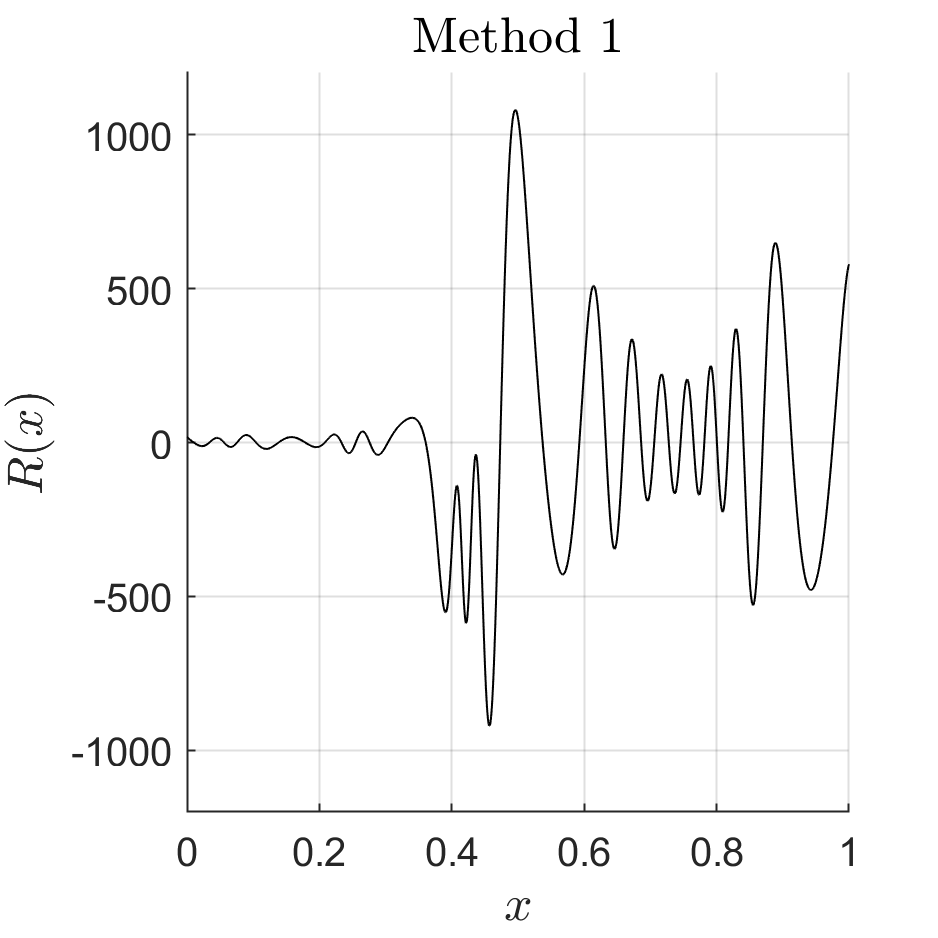}
\includegraphics[width=0.32\linewidth]{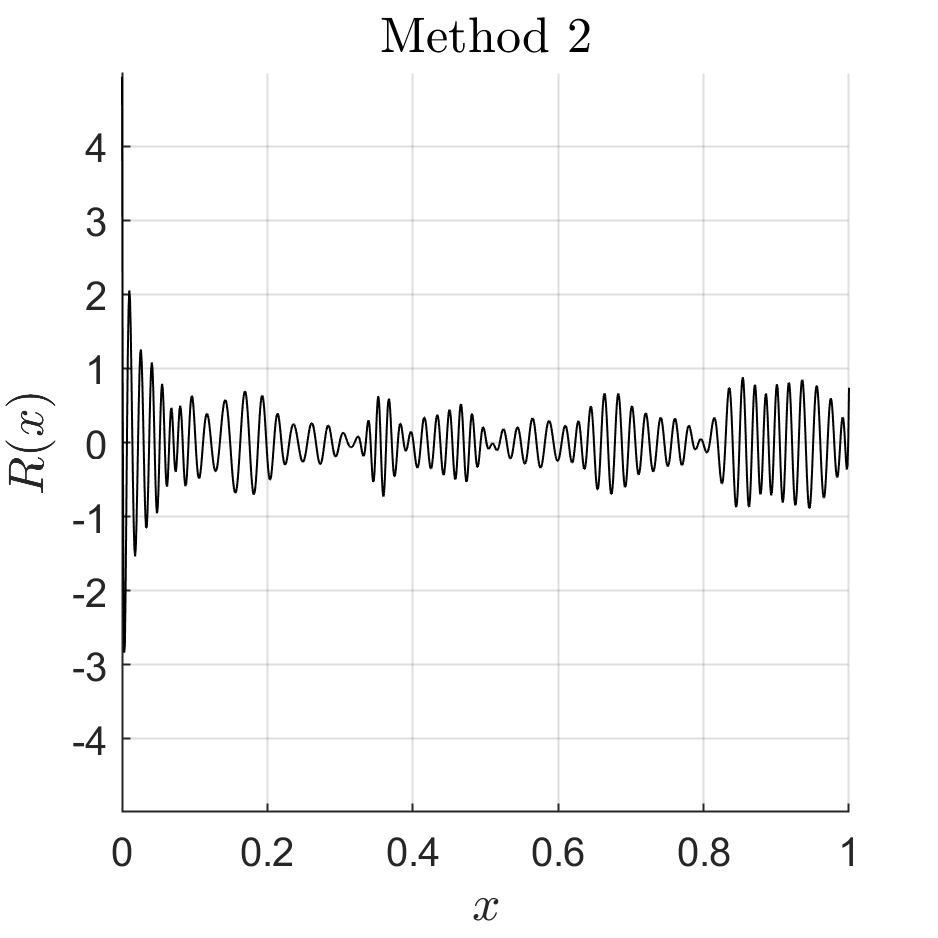}
\includegraphics[width=0.32\linewidth]{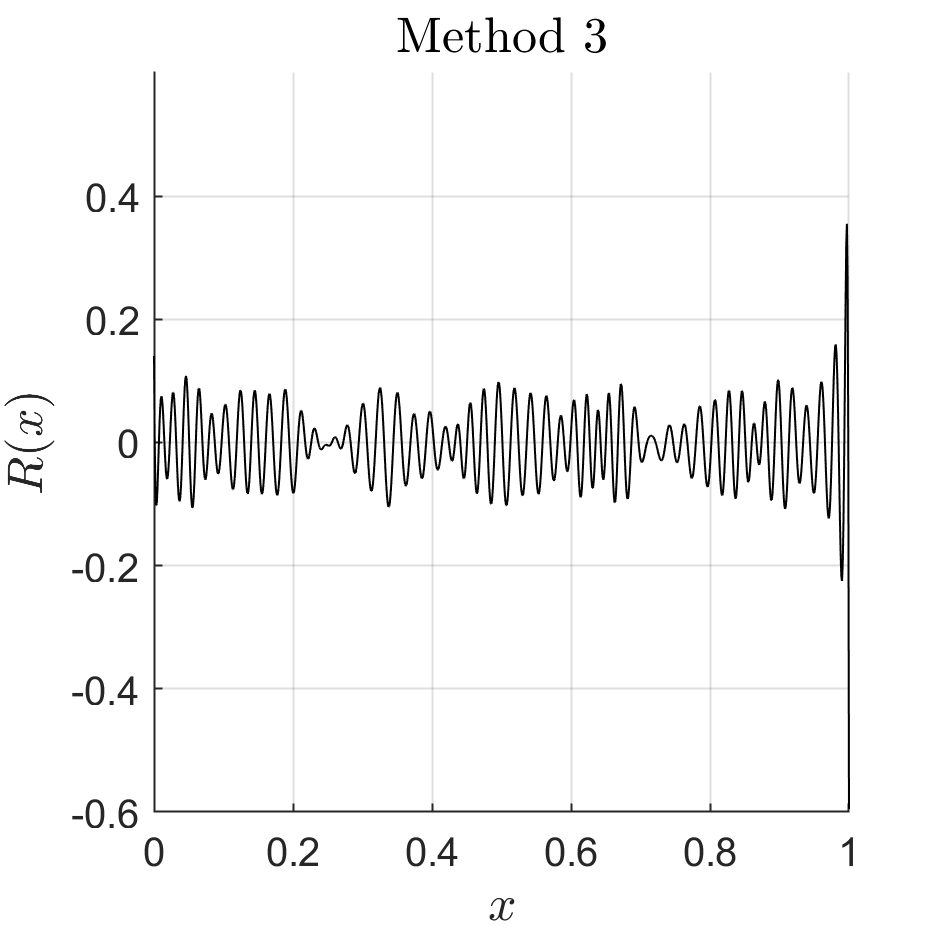}
\caption{Results from Example~\ref{example:hf} in Section~\ref{sec:high_freq}: Residual \(R(x)\) associated with the partial differential equation at the end of the training using Methods 1, 2, and 3. Note that the scale along the \(y\)-axis is different from one plot to the other.}
\label{fig:high_freq_residual}
\end{figure}

Finally, we show in Figure~\ref{fig:high_freq_error} the pointwise error \(e(x)= u(x)-\tilde{u}(x)\) obtained at the end of the training for Methods~1, 2, and~3. Note that the scale along the \(y\)-axis is different on the graphs. As expected, the pointwise error obtained by Method~1 is of the same order as the solution itself. Moreover, we observe that the maximum value of \(|e(x)|\) using Method~3 is smaller than that obtained with Method~2. Hence, the architecture presented in Method~3 yields a better solution when compared to the other two methods. We observed in Example~\ref{example:normalization} that smaller approximation errors contained higher frequencies. The picture is slightly different here. If we closely examine the pointwise error obtained by Methods~2 or~3, we observe that the error contains both a low-frequency component of large amplitude and a high-frequency component of small amplitude. In order to explain this phenomenon, we plot in Figure~\ref{fig:high_freq_residual} the residual \(R(x)\) associated with the partial differential equation for the three methods. For Method 1, we observe that the residual is still very large by the end of the training since the method did not converge. For Methods 2 and 3, the residual is a high-frequency function as the second-order derivatives of the solution tend to amplify its high-frequency components, as confirmed by the trivial calculation:
\[
\dfrac{d^2}{dx^2}\sin(\omega x) = -\omega^2 \sin(\omega x).
\]
It follows that the high-frequency components of the solution will be reduced first since the training is based on minimizing the residual of the partial differential equation.
On the other hand, the error, see e.g.\ Figure~\ref{fig:high_freq_error} (middle) or (right), includes some low-frequency contributions, which are imperceptible in the plot of the residual. To further reduce the pointwise error, the objective should then be to reduce the low-frequency modes alone, without the need to reduce the high frequencies whose amplitudes are smaller.
\end{example}

In summary, we have seen through numerical experiments that the accuracy of the solutions may be affected by the scale of the problem data and the range of frequencies inherent to the solutions. The methodology that we describe below allows one to address these issues, namely to control the error within machine precision in neural network solutions using the PINNs approach.
%{i.e.\ residual reduction}.

%============================================================
\section{Multi-level neural networks}
\label{sec:multi-nets}

In this section, we describe the multi-level neural networks, whose main objective is to improve the accuracy of the solutions obtained by PINNs. Supposing that an approximation~\(\tilde{u}\) of the solution~\(u\) to  Problem~\eqref{eq:residual}-\eqref{eq:bc} has been computed, the error in \(\tilde{u}\) is defined as \(e(\bx)=u(\bx)-\tilde{u}(\bx)\) and satisfies:
\[
\begin{aligned}
R(\bx,u(\bx)) = f(\bx) - Au(\bx) 
= f(\bx) - A\tilde{u}(\bx) - Ae(\bx) 
= R(\bx,\tilde{u}(\bx)) - Ae(\bx) &= 0, \quad \forall \bx \in \Omega, \\
B(\bx, u(\bx)) = B(\bx,\tilde{u}(\bx)) +  B(\bx,e(\bx)) = B(\bx,e(\bx)) &= 0, \quad \forall \bx \in \partial \Omega,  
\end{aligned}
\]
where we have used the fact that \(A\) and \(B\) are linear operators and $\tilde{u}$ strongly verifies the boundary condition. In other words, the error function $e(x)$ satisfies the new problem in the residual form:
\begin{align}
\tilde{R}(\bx,e(\bx)) = R(\bx,\tilde{u}(\bx)) - Ae(\bx) &= 0, \quad \forall \bx \in \Omega, \\
B(\bx,e(\bx)) &= 0, \quad \forall \bx \in \partial \Omega.  
\end{align}
We notice that the above problem for the error has exactly the same structure as the original problem, with maybe two exceptions: 1) the source term \(R(\bx,\tilde{u}(\bx))\) in the error equation may be small, 2) the error $e(\bx)$ may be prone to higher frequency components than in \(\tilde{u}\). Our earlier observations suggest we find an approximation $\tilde{e}$ of the error using the PINNs approach after normalizing the source term by a scaling parameter \(\mu\), in a way that scales the error to a unit maximum amplitude. The new problem becomes:
\begin{align}
\tilde{R}(\bx,e(\bx)) = \mu R(\bx,\tilde{u}(\bx)) - Ae(\bx) &= 0, \quad \forall \bx \in \Omega, \\
B(\bx,e(\bx)) &= 0, \quad \forall \bx \in \partial \Omega.  
\end{align}
The dimension of the new neural network to approximate \(e\) should be larger than that used to find \(\tilde{u}\), due to the presence of higher frequency modes in \(e\). In particular, the number of wave numbers \(M\) in the Fourier feature mapping should be increased. The idea is to some extent akin to a posteriori error estimation techniques developed for Finite Element methods, see e.g.~\cite{bank-weiser-1993,prudhomme-oden-1999,ainsworth2000,prudhomme-oden-2001,bangerth2003adaptive}. Finally, one should expect that the optimization algorithm should once again reach a plateau after a certain number of iterations and that the process should be repeated to estimate a new correction to the error \(e\).

We thus propose an iterative procedure, referred here to as the ``multi-level neural network method", in order to improve the accuracy of the solutions when using PINNs (or any other neural network procedure based on residual reduction). We start by modifying the notation due to the iterative nature of the process. As mentioned in the previous section, the source term \(f\) may need to be normalized by a scaling parameter $\mu_0$, so that we reconsider the initial solution \(u_0\) satisfying a problem in the form:
\begin{align}
R_0(\bx,u_0(\bx)) = \mu_0 f(\bx) - Au_0(\bx) &= 0, \quad \forall \bx \in \Omega, \\
B(\bx,u_0(\bx)) &= 0, \quad \forall \bx \in \partial \Omega.  
\end{align}
The above problem can then be approximated using a neural network to obtain an approximation \(\tilde{u}_0\) of \(u_0\). 
Hence, the first approximation \(\tilde{u}\) to \(u\) reads after scaling \(\tilde{u}_0\) with \(\mu_0\):
\begin{equation}
\tilde{u}(\bx) = \dfrac{1}{\mu_0}\tilde{u}_0(\bx).
\end{equation}
We would like now to estimate the error in \(\tilde{u}\). However, we find it easier to work in terms of \(\tilde{u}_0\). Therefore, we look for a new correction \(u_1\) that solves the problem:
\begin{align}
R_1(\bx,u_1(\bx)) = \mu_1 R_0(\bx,\tilde{u}_0(\bx)) - Au_1(\bx) &= 0, \quad \forall \bx \in \Omega, \\
B(\bx,u_1(\bx)) &= 0, \quad \forall \bx \in \partial \Omega.  
\end{align}
Once again, one can compute an approximation \(\tilde{u}_1\) of \(u_1\) using PINNs. Since \(\tilde{u}_1\) can be viewed as the normalized correction to the error in \(\tilde{u}_0(\bx)\), the new approximation to \(u\) is now given by:
\begin{equation}
\tilde u(\bx) = \dfrac{1}{\mu_0}\tilde u_0(\bx) + \dfrac{1}{\mu_0\mu_1}\tilde u_1(\bx).
\end{equation}
The process can be repeated \(L\) times to find corrections \(u_i\) at each level \(i=1,\ldots,L\) given the prior approximations \(\tilde{u}_0\), \(\tilde{u}_1\), \ldots, \(\tilde{u}_{i-1}\).
Each new correction \(u_i\) then satisfies the boundary-value problem:
\begin{align}
R_i(\bx,u_i(\bx)) = \mu_i R_{i-1}(\bx,\tilde{u}_{i-1}(\bx)) - A u_i(\bx) &= 0, \quad \forall \bx \in \Omega, \\
B(\bx,u_i(\bx)) &= 0, \quad \forall \bx \in \partial \Omega.  
\end{align}
After finding an approximation \(\tilde{u}_i\) to each of those problems up to level~\(i\), one can obtain a new approximation \(\tilde{u}\) of \(u\) such that:
\begin{equation}
\tilde{u}(\bx) = 
\dfrac{1}{\mu_0}\tilde u_0(\bx) 
+ \dfrac{1}{\mu_0\mu_1}\tilde u_1(\bx) 
+ \ldots
+ \dfrac{1}{\mu_0\mu_1\ldots\mu_i}\tilde u_i(\bx).
\end{equation}
Once the approximations \(\tilde{u}_i\) have been found at all levels \(i=0,\ldots,L\), the final approximation at the end of the process would then be given by:
\begin{equation}
\tilde{u}(\bx) = \sum_{i=0}^{L} \dfrac{1 }{\Pi_{j=0}^i \mu_j} \tilde{u}_i(\bx).
\end{equation}

Using PINNs, the neural network approximation \(\tilde{u}_i\) (which implicitly depends on the network parameters \(\theta\)) for each error correction will be obtained by solving the following minimization problem:
\begin{equation}
\label{eq:losspinn_strong_Ri}
\min_{\theta \in \mathbb R^{N_{\theta,i}}} 
\mathcal L_i(\theta) 
= 
\min_{\theta \in \mathbb R^{N_{\theta,i}}} 
\int_\Omega R_i\big(\bx,\tilde{u}_i(\bx) \big)^2 dx,
\end{equation}
where \(N_{\theta,i}\) denotes the dimension of the function space generated by the neural network used at level~\(i\). We recall that the boundary conditions are strongly imposed and, hence, do not appear in the loss functions \(\mathcal L_i(\theta)\). Since each correction $\tilde{u}_i$ is expected to have higher frequency contents, the size \(N_{\theta,i}\) of the networks should be increased at each level. Moreover, the number of iterations used in the optimization algorithms Adam and L-BFGS will be increased as well, since more iterations are usually needed to approximate higher frequency functions. For illustration purposes, we consider a simple one-dimensional numerical example and use once again the setting of Example~\ref{example:optimizers} in Section~\ref{sec:normalize}.

\begin{example}
\label{example:multi_level1}
We solve Problem~\eqref{eq:poisson_problem} with \(k=2\) whose exact solution is given by~\eqref{eq:solution_exact}. We consider three levels of the multi-level neural networks, i.e.\ \(L=3\), in addition to the initial solve, so that the approximation \(\tilde{u}\) will be obtained using four sequential neural networks. We choose networks with a single hidden layer of width \(N_1\) given by \(\{10,20,40,20\}\). The networks are first trained with 4,000 iterations of Adam followed by \(\{200,400,600,0\}\) iterations of L-BFGS. The mappings of the input and boundary conditions are chosen with \(M = \{1,3,5,1\}\) wave numbers. In this example, the scaling parameter \(\mu_i\), \(i=0,\ldots,3\), for \(\tilde{u}_0\) and the three corrections \(\tilde{u}_i\), are chosen here as \(\mu_i = \{1,10^3,10^3,10^2\}\). In the next section, we will present a simple approach to evaluate these normalization factors. We note that the last network has been designed to approximate functions with a low-frequency content. This choice will be motivated below. 

\begin{figure}[tb]
\centering
\includegraphics[width=0.32\linewidth]{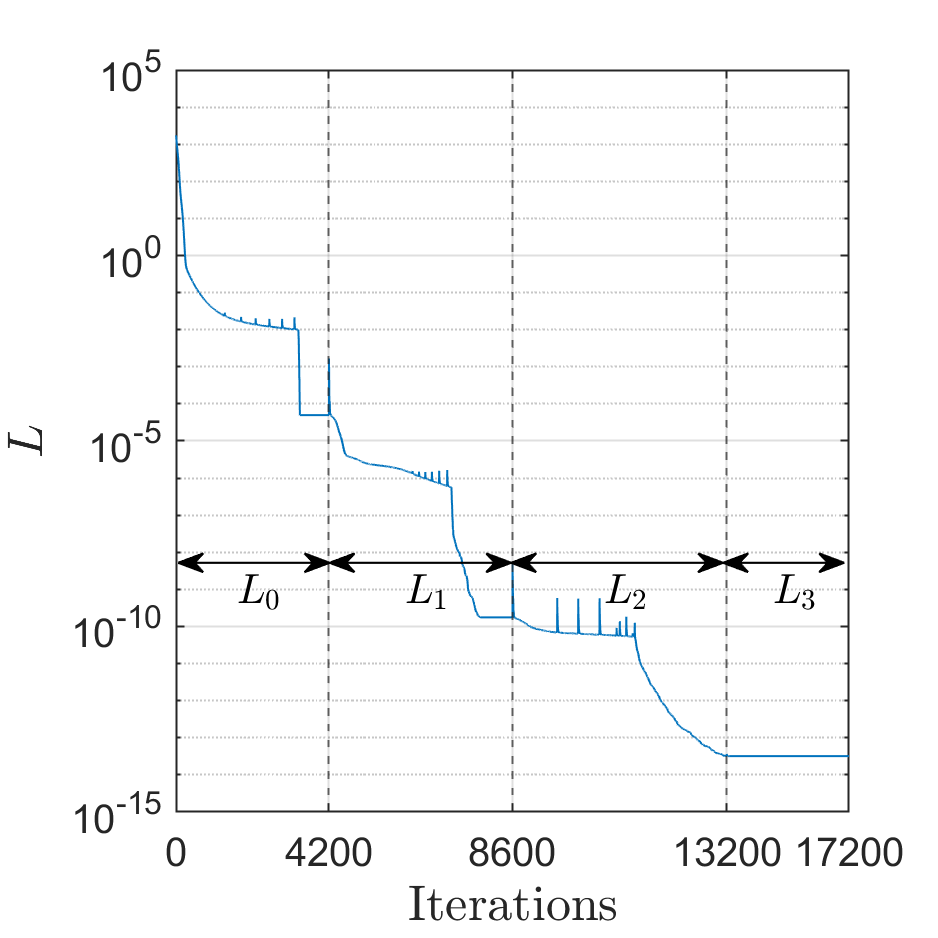}
\includegraphics[width=0.32\linewidth]{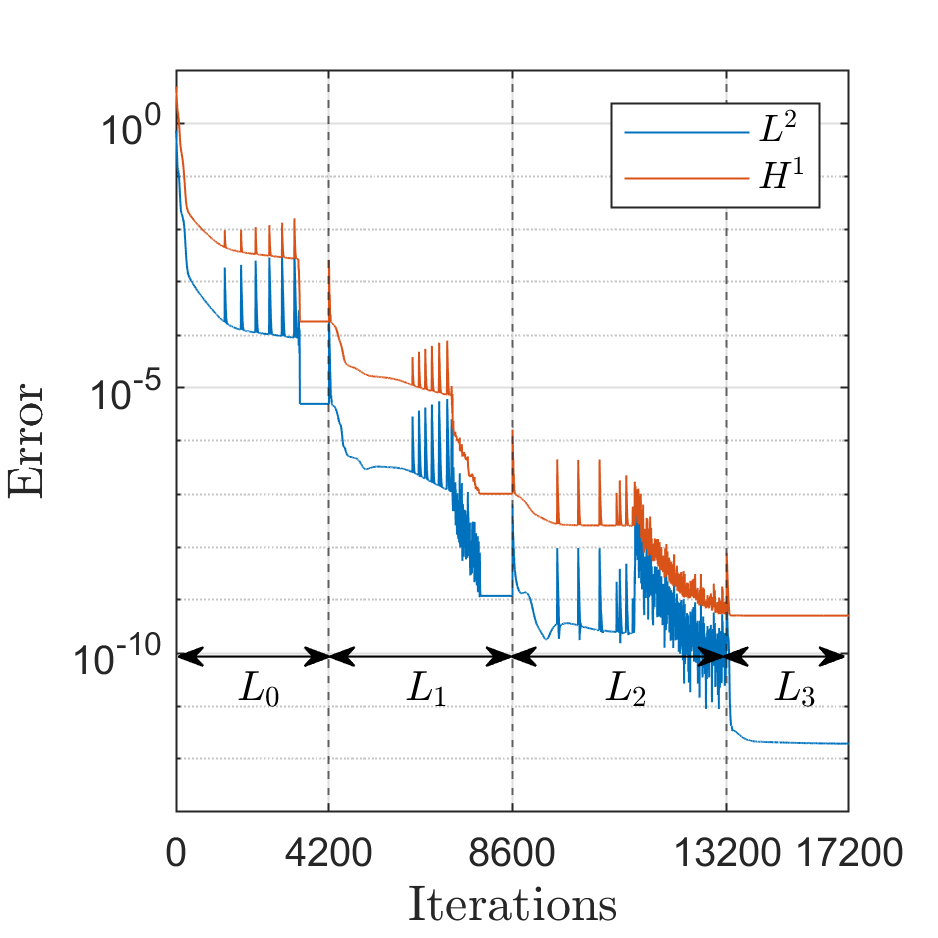}
\includegraphics[width=0.32\linewidth]{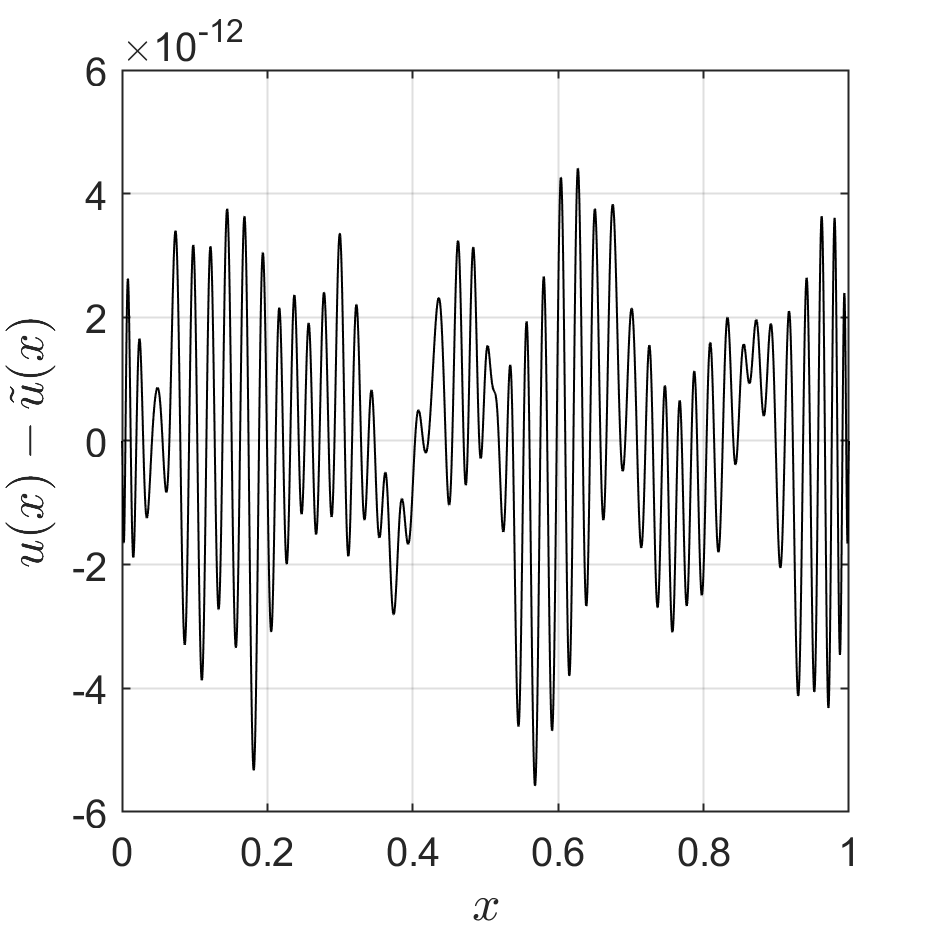}
\caption{Results from Example~\ref{example:multi_level1} in Section~\ref{sec:multi-nets}: (Left) Evolution of the loss function. (Middle) Evolution of the error \(e(x)=u(x)-\tilde{u}(x)\) measured in the \(L^2\) and \(H^1\) norms. (Right) Pointwise error after three error corrections. The regions shown by \(L_i\), \(i=0,\ldots,3\), in the first two graphs indicate the region in which the neural network at level~\(i\) is trained.}
\label{fig:example1_1}
\end{figure}

We present in Figure~\ref{fig:example1_1} the evolution of the loss function and of the errors in the \(L^2\) and \(H^1\) norms with respect to the number of optimization iterations, along with the pointwise error at the end of the training. We first observe that each error correction allows one to converge closer to the exact solution. More precisely, we gain almost seven orders of magnitude in the \(L^2\) error thanks to the introduced corrections. Indeed, after three corrections, the maximum pointwise error is around \(6\times10^{-12}\), which is much smaller than the error we obtained with \(\tilde{u}_0\) alone. To better explain our choice of the number \(M\) of wave numbers at each level, we show in Figure~\ref{fig:example1_2} the computed corrections \(\tilde{u}_i\). We observe that each correction approximates higher frequency functions than the previous one, except \(\tilde{u}_3\). In fact, once we start approximating the high-frequency errors, it becomes harder to capture the low-frequencies with larger amplitudes. This phenomenon was actually observed and described in Example~\ref{example:hf}. Here, we see that the loss function eventually decreases during the training of \(\tilde{u}_2\), but that the \(L^2\) error has the tendency to oscillate while slightly decreasing. It turns out that this behavior can be attributed to the choice of the loss function \(\mathcal L_2\), in which the higher frequencies are penalized more than the lower ones. In other words, we have specifically designed the last network to approximate only low-frequency functions and be trained using Adam only. Thanks to this architecture, the \(L^2\) error significantly decreases during the training of \(\tilde{u}_3\), without a noticeable decrease in the loss function. As a remark, longer training for \(\tilde{u}_2\) would correct the lower frequencies while correcting high frequencies with smaller amplitudes, but it is in our opinion more efficient, with respect to the number of iterations, to simply introduce a new network targeting the low frequencies.

\begin{figure}[tb]
\centering
\includegraphics[width=0.35\linewidth]{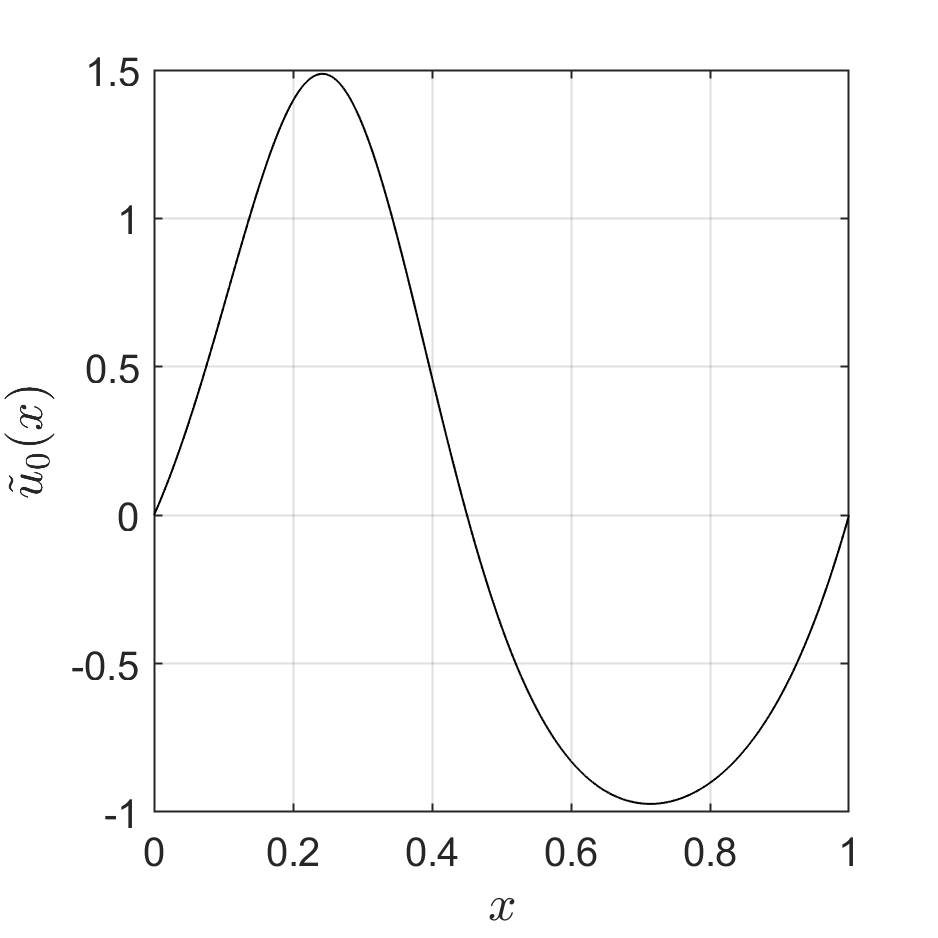}
\includegraphics[width=0.35\linewidth]{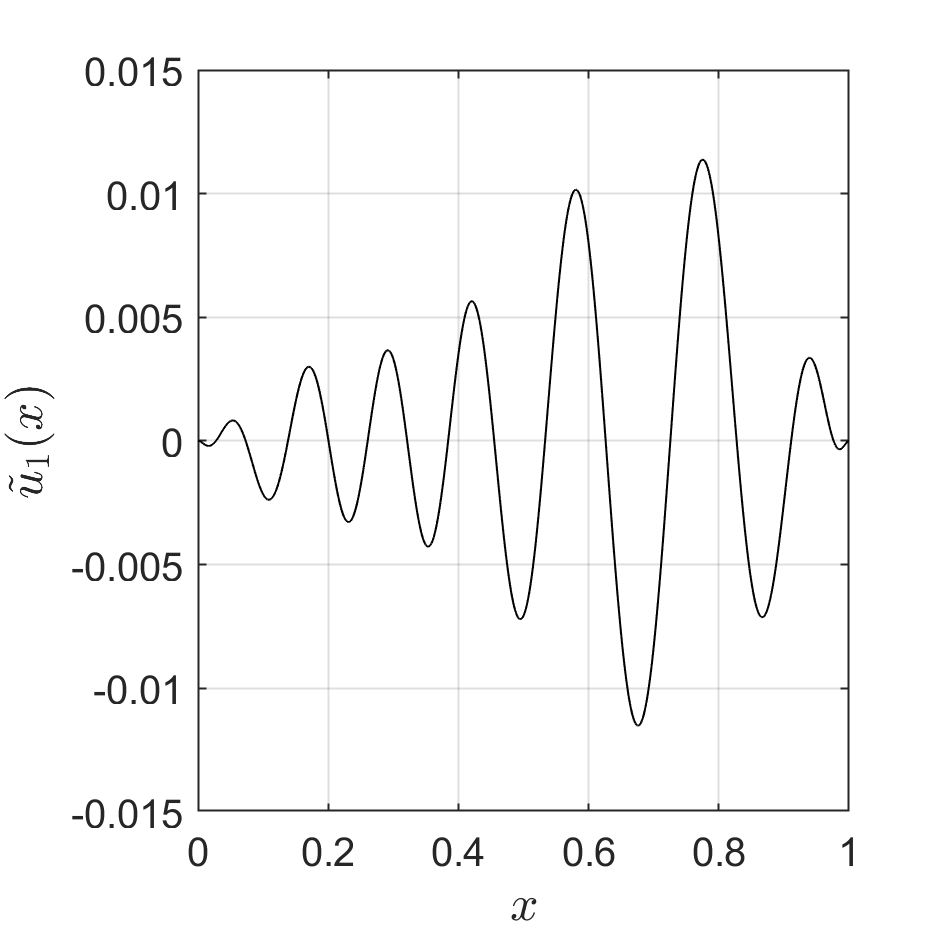}
\\
\includegraphics[width=0.35\linewidth]{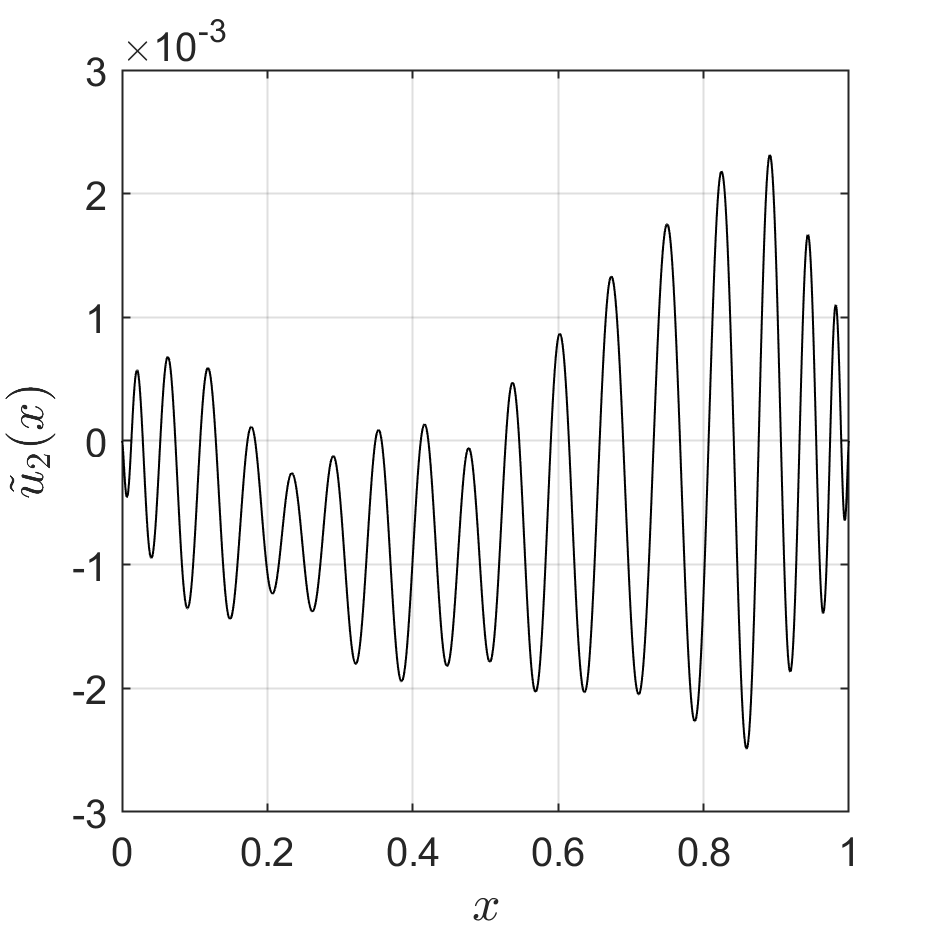}
\includegraphics[width=0.35\linewidth]{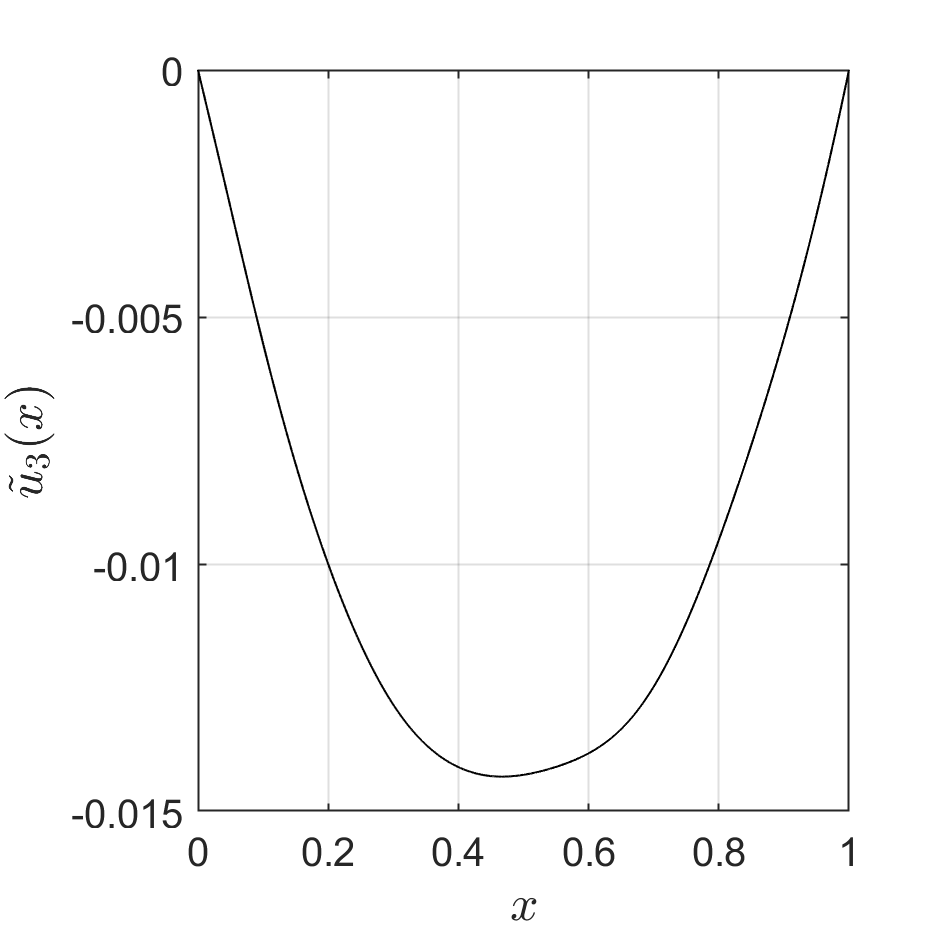}
\caption{Results from Example~\ref{example:multi_level1} in Section~\ref{sec:multi-nets}: Approximation \(\tilde{u}_0(x)\) and corrections \(\tilde{u}_i(x)\), \(i=1,2,3\).}
\label{fig:example1_2}
\end{figure}

\end{example}

In the last example, we observe that the maximal values in the three corrections \(\tilde{u}_i\), \(i=1,2,3\), range in absolute value from \(0.002\) to \(0.015\), see Figure~\ref{fig:example1_2}, whereas their values should ideally be of order one if proper normalization were used.
The reason is that we provided a priori values for the scaling parameters \(\mu_i\). It appears that those values were not optimal, i.e.\ too small, yielding solutions whose amplitudes were two to three orders of magnitude smaller than the ones we should expect. The multi-level neural network approach was still able to improve the accuracy of the solution despite sub-optimal values of the scaling parameters.

Ideally, one would like to have a method to uncover appropriate values of the scaling parameters. Unfortunately, it is not a straightforward task, that of predicting the amplitude of the remaining error in order to correctly normalize the residual term in the partial differential equation. We propose here a simple approach based on the Extreme Learning Method~\cite{huang2011extreme}. The main idea of the Extreme Learning Method is to use a neural network with a single hidden layer, to fix the weight and bias parameters of the hidden layer, and to minimize the loss function with respect to the parameters of only the output layer by a least squares method. We propose here to utilize the Extreme Learning Method to obtain a coarse prediction for each correction \(\tilde{u}_i\). The solution might not be very accurate, but it should provide a reasonable estimate of the amplitude of the correction function, which can be employed to adjust the normalization parameter \(\mu_i\). Moreover, the approach has the merit of being very fast and scale-independent. We assess its performance in the next numerical example and show that it allows one to further improve the accuracy of the multi-level neural network solution.

\begin{example}
\label{example:multi_level2}
We use the exact same setting as described in Example~\ref{example:multi_level1} but we now employ the Extreme Learning Method to normalize the residual terms, as explained above. We observe in Figure~\ref{fig:example2_1} that the proposed normalization technique leads at the end of the training to errors \(e(x)=u(x)-\tilde{u}(x)\) within machine precision. We actually gain about two orders of magnitude in the error with respect to both norms over the results obtained in Example~\ref{example:multi_level1}. Even more striking, the approximations of the corrections \(\tilde{u}_i\) all have amplitudes very close to unity, see Figure~\ref{fig:example2_2}, which confirms the efficiency of the proposed approach.

\begin{figure}[tb]
\centering
\includegraphics[width=0.32\linewidth]{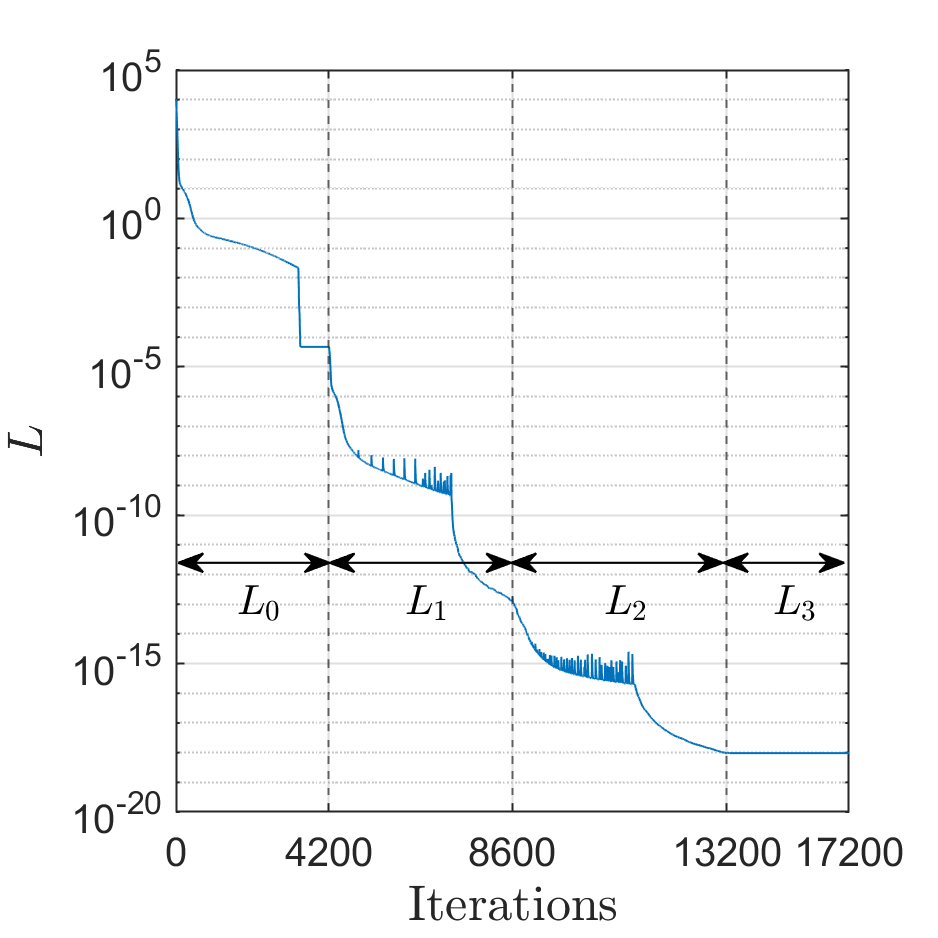}
\includegraphics[width=0.32\linewidth]{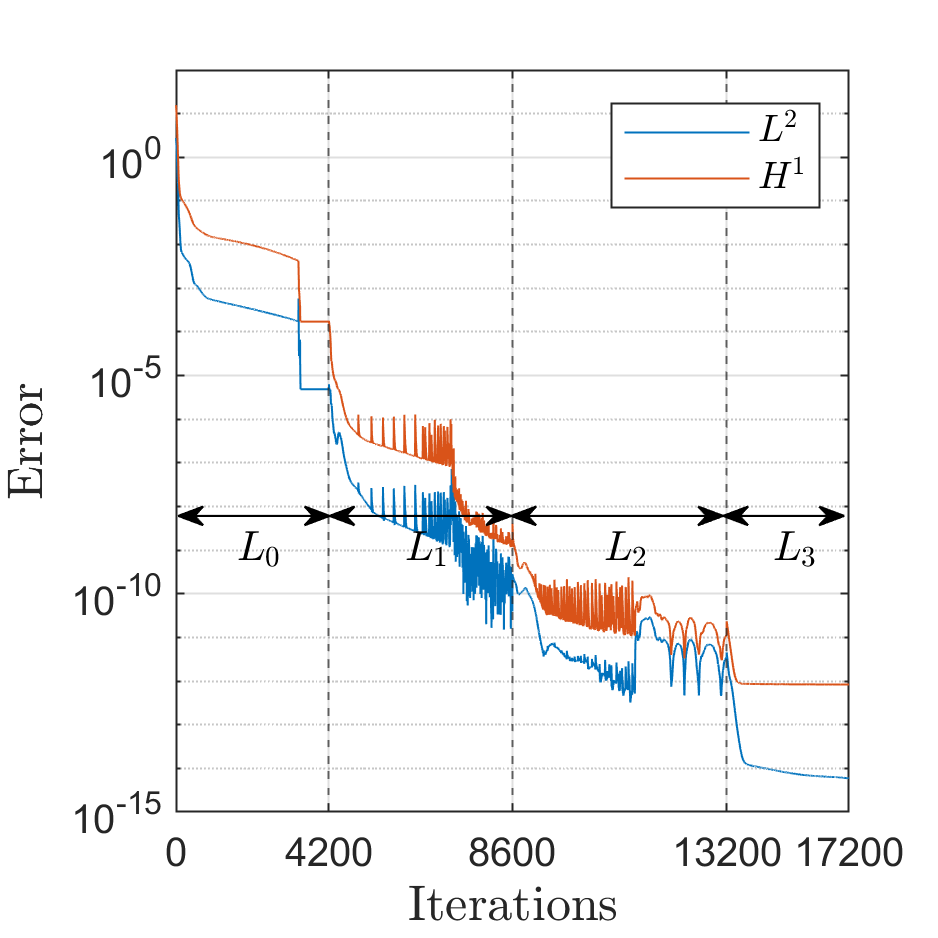}
\includegraphics[width=0.32\linewidth]{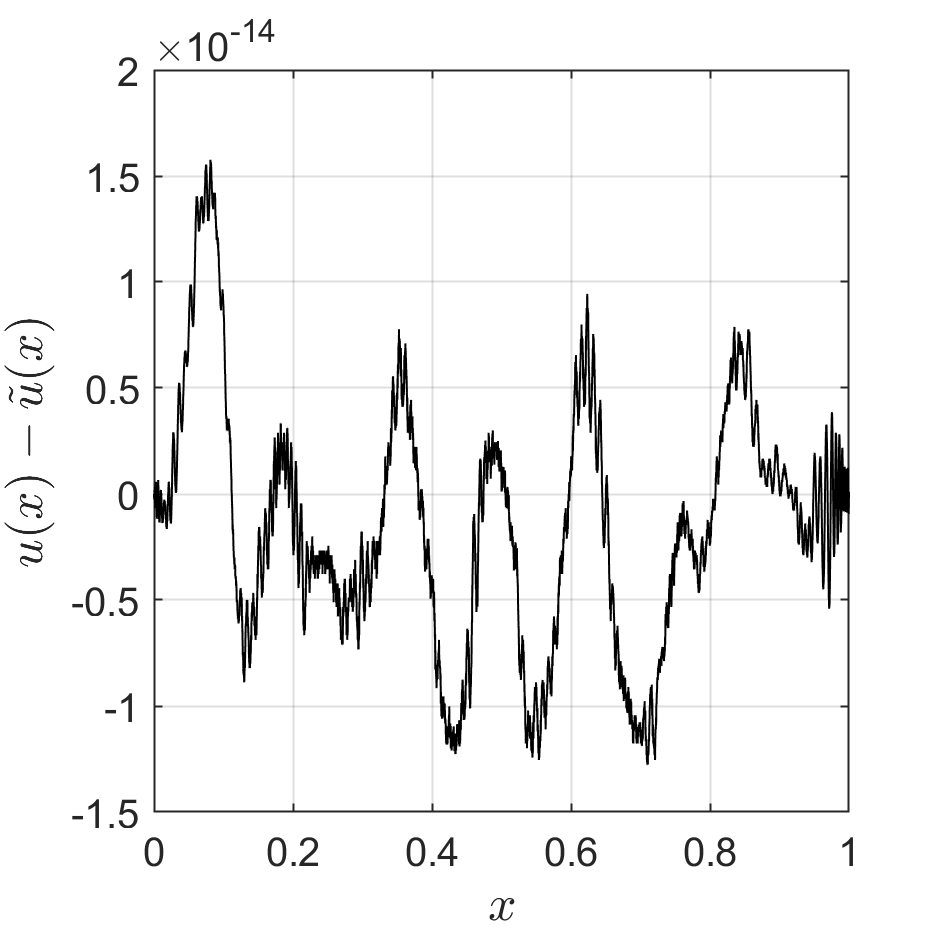}
\caption{Results from Example~\ref{example:multi_level2} in Section~\ref{sec:multi-nets}: (Left) Evolution of the loss function. (Middle) Evolution of the error \(e(x)=u(x)-\tilde{u}(x)\) measured in the \(L^2\) and \(H^1\) norms. (Right) Pointwise error after three error corrections.}
\label{fig:example2_1}
\end{figure}

\begin{figure}[tb]
\centering
\includegraphics[width=0.35\linewidth]{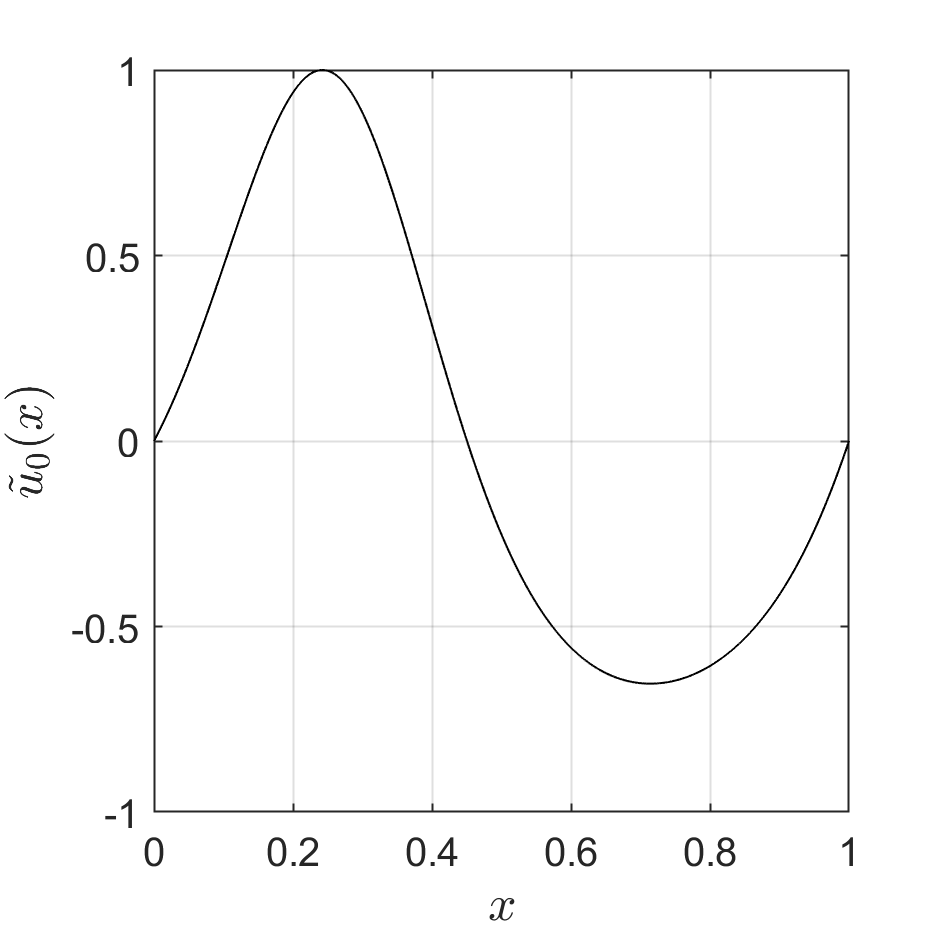}
\includegraphics[width=0.35\linewidth]{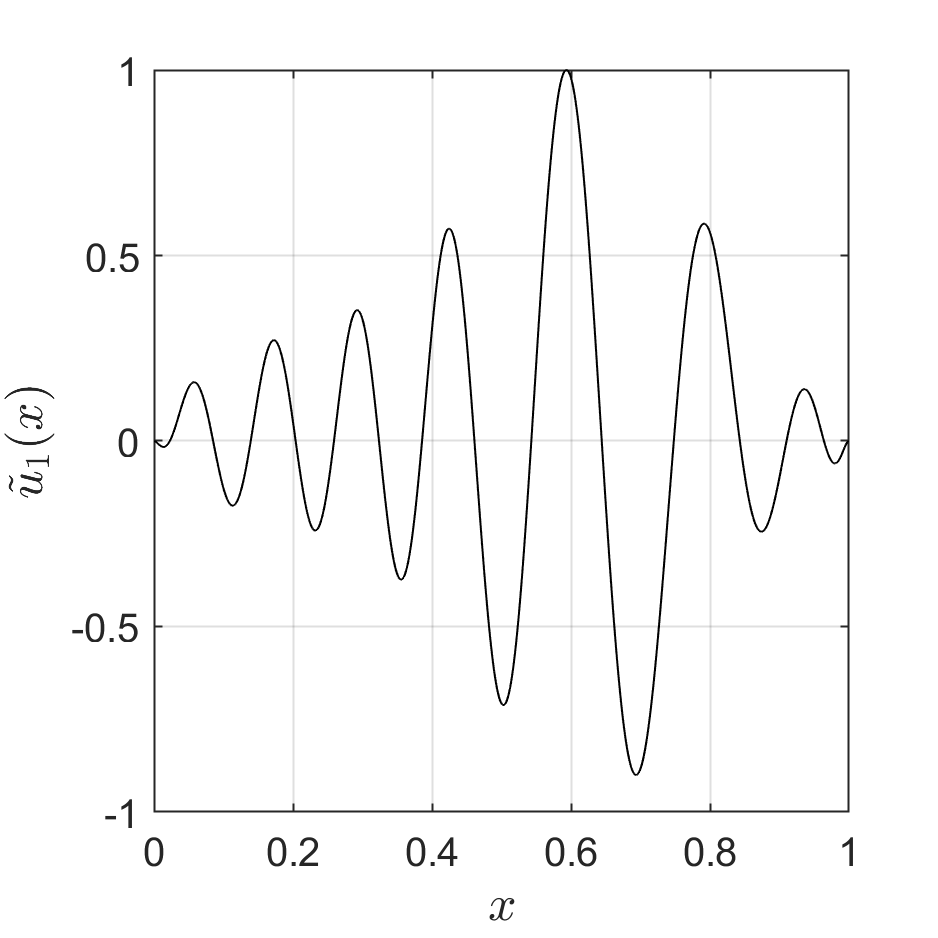}
\\
\includegraphics[width=0.35\linewidth]{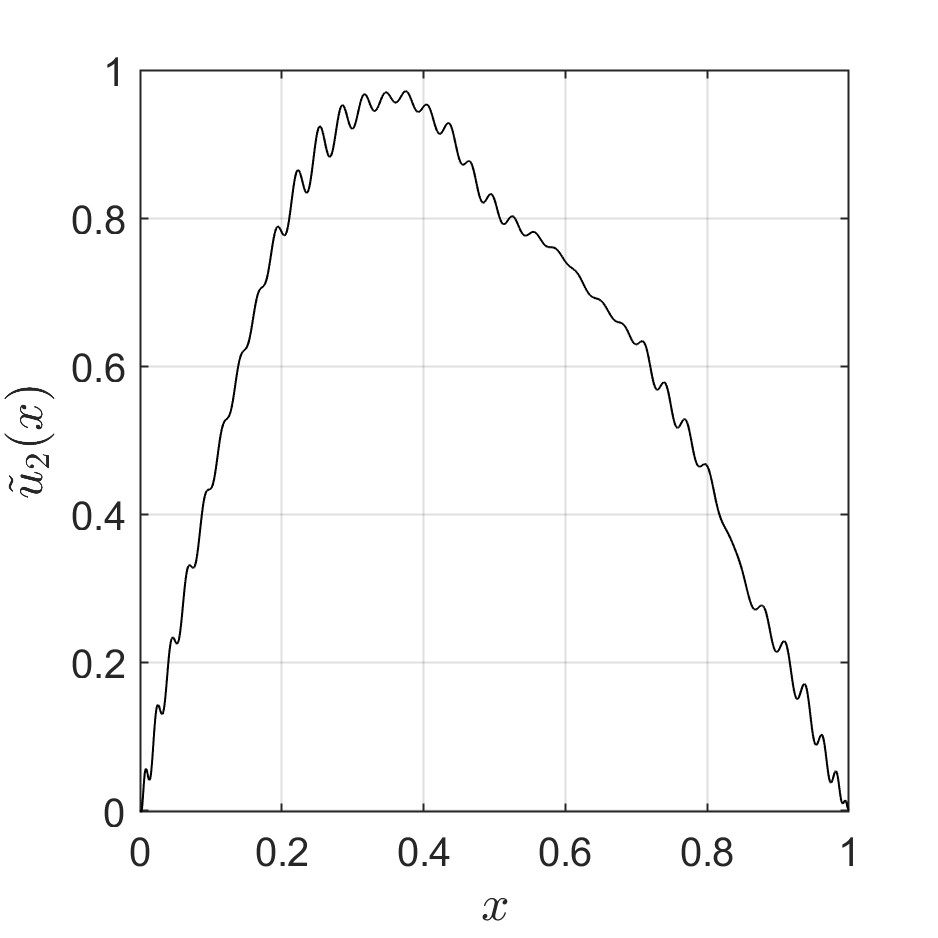}
\includegraphics[width=0.35\linewidth]{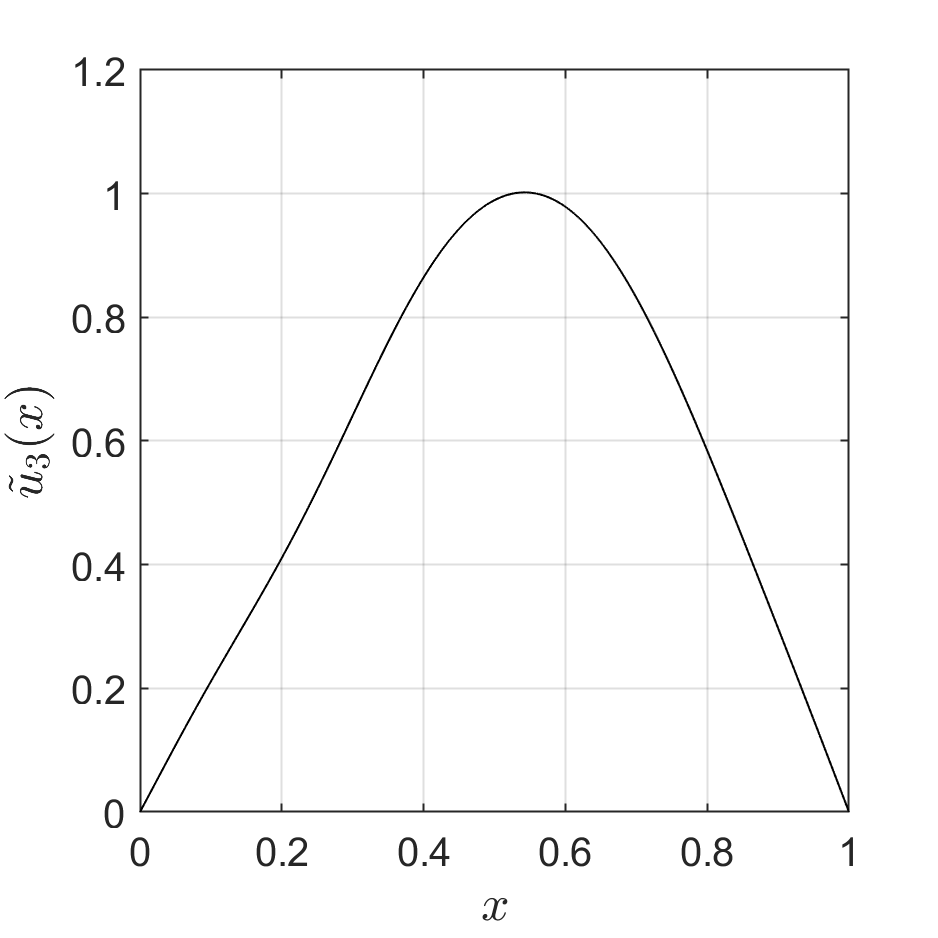}
\caption{Results from Example~\ref{example:multi_level2} in Section~\ref{sec:multi-nets}: Approximation \(\tilde{u}_0(x)\) and corrections \(\tilde{u}_i(x)\), \(i=1,2,3\).}
\label{fig:example2_2}
\end{figure}

\end{example}

%============================================================
\section{Numerical results}
\label{sec:numresults}
%============================================================

In this section, we present a series of numerical examples to illustrate the whole potential of the multi-level neural networks to reduce errors in neural network approximations of boundary-value problems in one and two dimensions. The computational domain in each of the examples is defined as \(\Omega = [0,1]^d\), with \(d=1\) or~\(2\). The solutions to the problems will all be submitted to homogeneous Dirichlet boundary conditions, unless explicitly stated otherwise. The solutions and error corrections computed with the multi-level neural network approach shall be consistently approximated by the neural network architecture provided in~\eqref{eqn:FNN_hf2}, for which the vector of wave numbers \(\boldsymbol{\omega}_M\) is constructed from a geometric series, as described in Section~\ref{sec:high_freq}. The normalization of the source terms  is implemented through the use of the Extreme Learning Method, as described in Section~\ref{sec:multi-nets}. We again emphasize that the scaling along the horizontal axis in the convergence plots is actually different for the Adam iterations and L-BFGS iterations. The reason is simply to provide a clearer visualization of the L-BFGS training phase, given the notable difference in the number of iterations used in the Adam and L-BFGS algorithms. Finally, for each example, the number of levels is set to \(L=3\) and the values of the hyper-parameters for each network (number of hidden layers \(n\) and widths \(N_i\), number of Adam and L-BFGS iterations, number of wave numbers \(M\)) will be collected in a table. 

%============================================================
\subsection{Poisson problem in 1D}
\label{subsec:poissonk10}

We revisit once again Problem~\eqref{eq:poisson_problem} described in Example~\ref{example:optimizers}, this time with \(k=10\). Our objective here is to demonstrate the performance of the multi-level neural networks even in the case of solutions with high-frequency components. The solution is approximated using four sequential networks whose hyper-parameters are reported in Table~\ref{tbl:example_p10}. Similarly to the Example~\ref{example:multi_level1}, the last network is chosen in such a way that only the low-frequency modes are approximated, using the Adam optimizer only.

\begin{table}[tb]
\label{tbl:example_p10}
\centering
\renewcommand{\arraystretch}{1.4}
\begin{tabular}{c|c|c|c|c}
Hyper-parameters
& \(\tilde{u}_0\) & \(\tilde{u}_1\) 
& \(\tilde{u}_2\) & \(\tilde{u}_3\) \\ \hline
{\# Hidden layers \(n\)}    
& 1     & 1     & 1     & 1         \\ 
{Width \(N_1\)}            
& 10    & 20    & 40    & 40        \\ 
{\# Adam iterations}  
& 4,000 & 4,000 & 4,000 & 10,000    \\ 
{\# L-BFGS iterations} 
& 500   & 1,000 & 1,500 & 0         \\ 
{\# wave numbers \(M\)}
& 4     & 6     & 8     & 2         \\ \hline
\end{tabular}
\caption{Hyper-parameters used in the example of Section~\ref{subsec:poissonk10}.}
\end{table}

\begin{figure}[tb]
\centering
\includegraphics[width=0.32\linewidth]{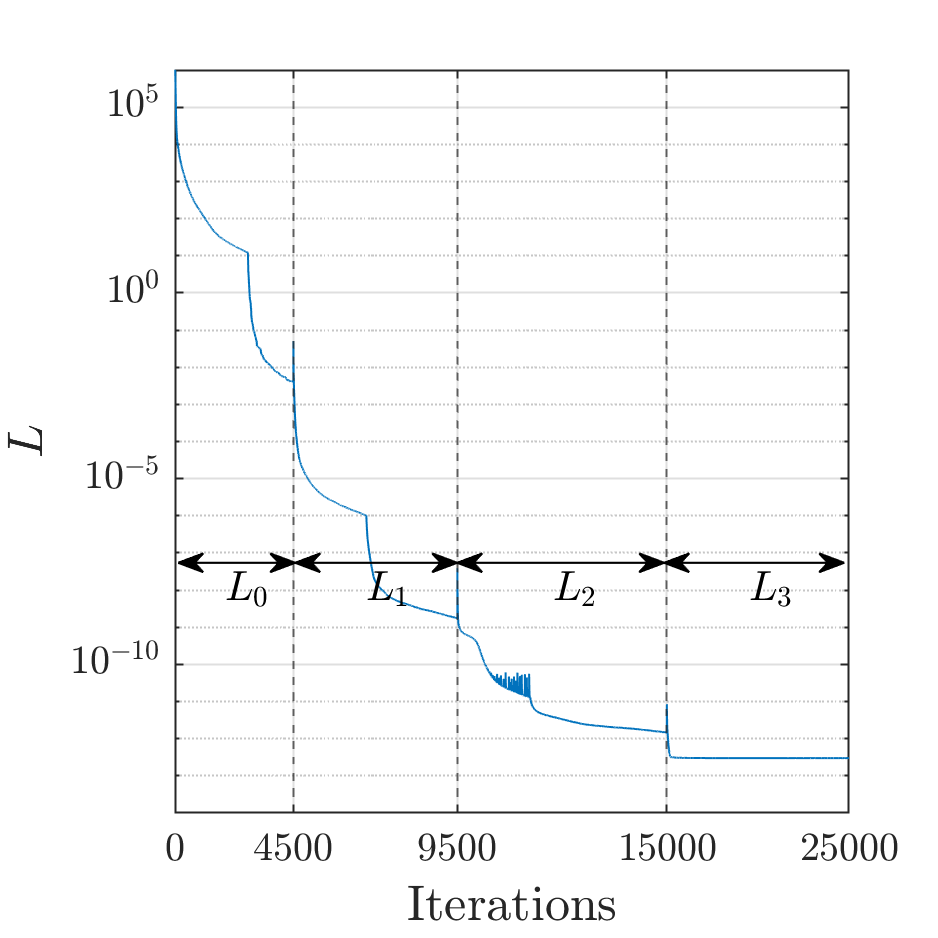}
\includegraphics[width=0.32\linewidth]{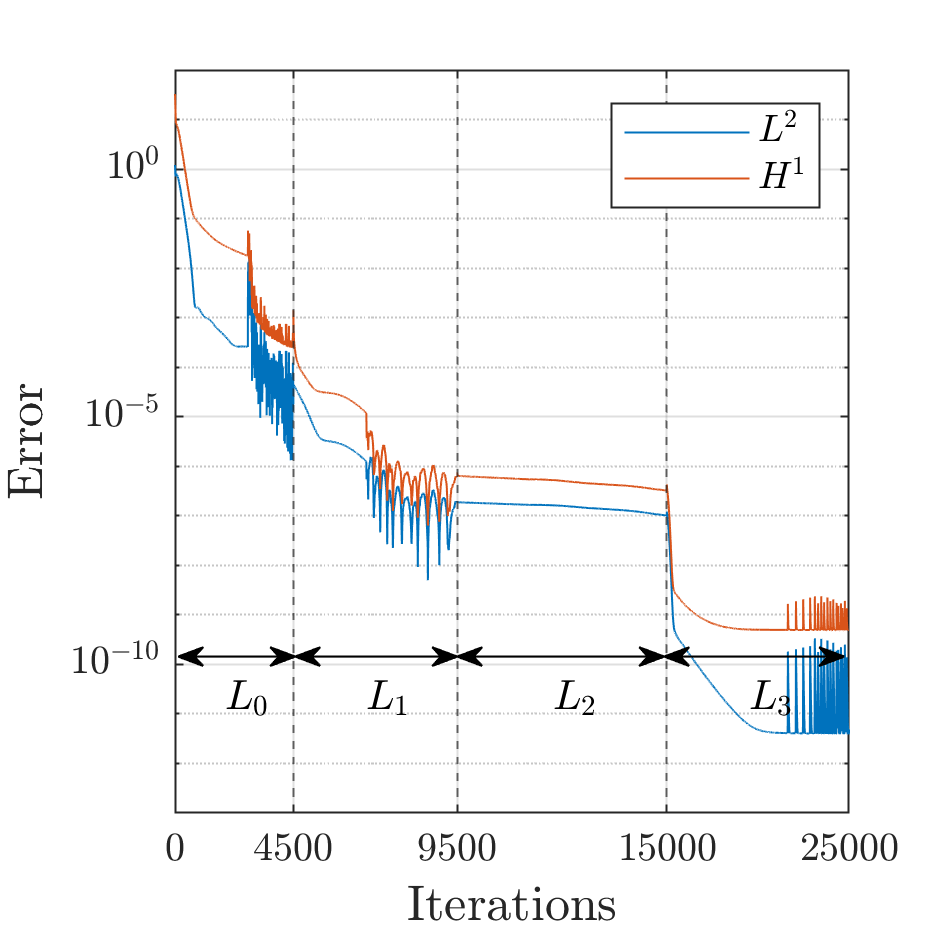}
\includegraphics[width=0.32\linewidth]{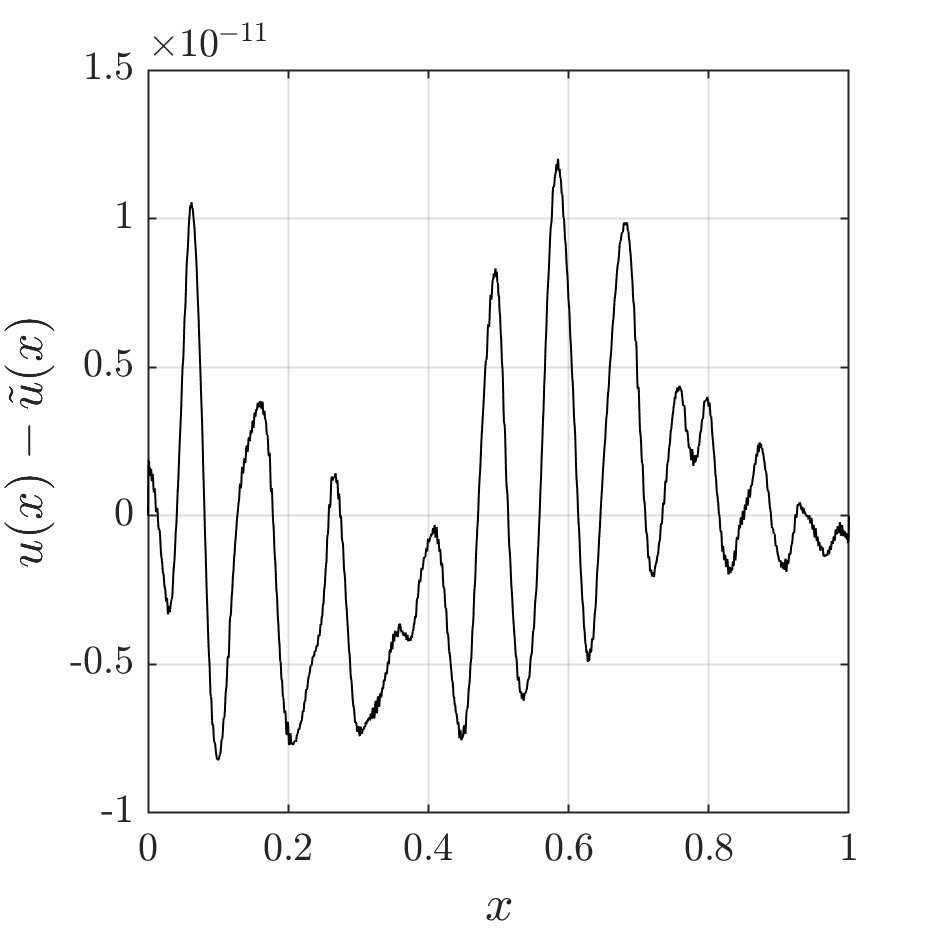}
\caption{
Example of Section~\ref{subsec:poissonk10}:  (Left) Evolution of the loss function. (Middle) Evolution of the error \(e(x)=u(x)-\tilde{u}(x)\) measured in the \(L^2\) and \(H^1\) norms. (Right) Pointwise error after three error corrections.}
\label{fig:example_p10}
\end{figure}

We plot in Figure~\ref{fig:example_p10} the evolution, during training, of the loss function (left) and the errors in the \(L^2\) and \(H^1\) norms (middle), along with the pointwise error at the end of the training (right). We observe that the loss function and the errors in both norms are reduced when using two corrections. During the second error correction, we notice that the reduction in the loss function did not yield a significant decrease in the \(L^2\) and \(H^1\) errors. As described in Example~\ref{example:multi_level1}, this is a consequence of our choice of the loss function, where higher frequencies are penalized more than the lower ones, which yields large low-frequency errors. This issue is addressed in the third correction that helps decrease the \(L^2\) and the \(H^1\) errors without significantly decreasing the loss function, since the role of the last network is mainly to capture the low frequencies. As described in Example~\ref{example:multi_level1}, this is a consequence of the specific choice of the hyper-parameters for the last network. In this example, we are able to attain a maximum pointwise error of around \(10^{-11}\) with four successive networks.

%============================================================
\subsection{Boundary-layer problem}
\label{subsec:BL}

\begin{figure}[tb]
\centering
\includegraphics[width=0.4\linewidth]{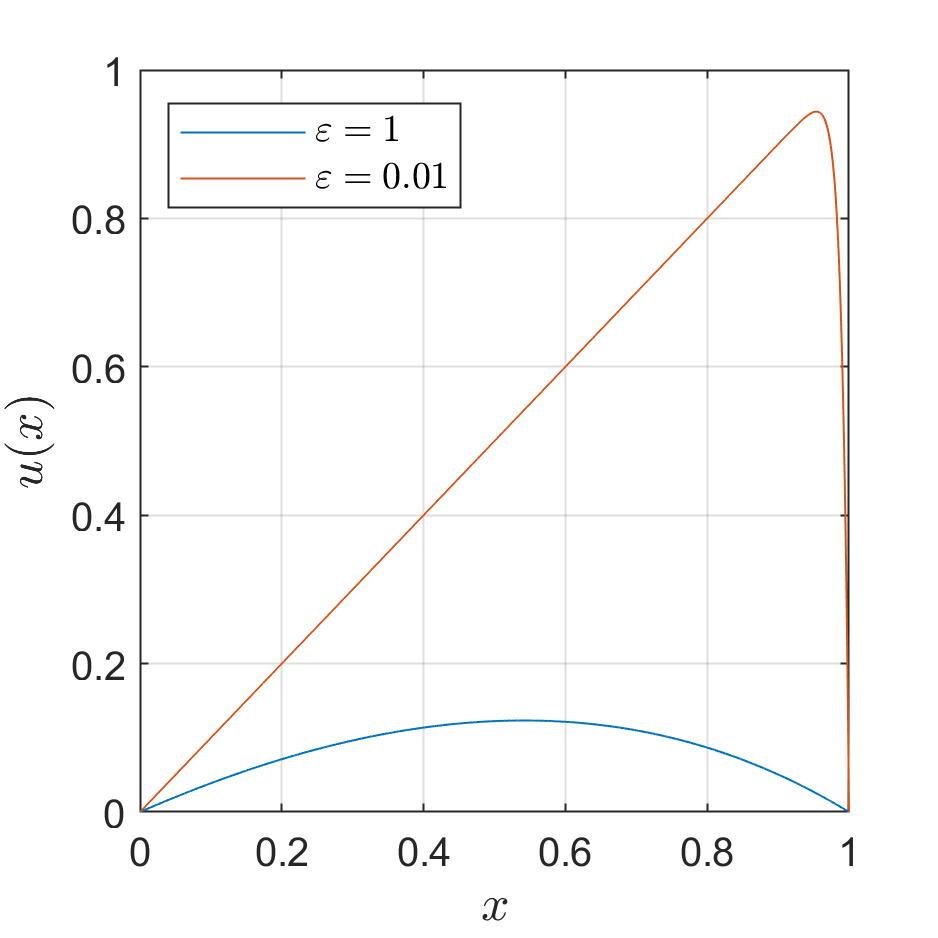}
\caption{Exact solutions when \(\varepsilon = 1\) and \(\varepsilon = 0.01\) to the boundary-layer problem of Section~\ref{subsec:BL}.}
\label{fig:example_p13_exact}
\end{figure}

In this section, we consider the convection-diffusion problem given by
\begin{equation}
\label{eq:BL_problem}
\begin{aligned}
-\varepsilon \partial_{xx} u(x) +\partial_{x} u(x)  &= 1, \qquad \forall x\in(0,1), \\
u(0) &=0, \\
u(1) &=0,
\end{aligned}
\end{equation}
where \(\varepsilon\) denotes a viscosity coefficient. We show in Figure~\ref{fig:example_p13_exact} the exact solutions to the problem when \(\varepsilon = 1\) and \(\varepsilon = 0.01\). As \(\varepsilon\) gets smaller, a sharp boundary layer is formed in the vicinity of \(x=1\), which makes the problem more challenging to approximate. Finite element approximations of the same problem without using any stabilization technique actually exhibit large oscillations whenever the mesh size is not fine enough to capture the boundary layer.
We apply the multi-level neural network method to both cases using the hyper-parameters given in Table~\ref{tbl:example_p13} for \(\varepsilon = 1\) and Table~\ref{tbl:example_p14} for \(\varepsilon = 0.01\).

\begin{table}[tb]
\centering
\renewcommand{\arraystretch}{1.4}
\begin{tabular}{c|c|c|c|c}
Hyper-parameters
& \(\tilde{u}_0\) & \(\tilde{u}_1\) 
& \(\tilde{u}_2\) & \(\tilde{u}_3\) \\ \hline
{\# Hidden layers \(n\)}    
& 1     & 1     & 1     & 1         \\ 
{Width \(N_1\)}            
& 5     & 10    & 20    & 40        \\ 
{\# Adam iterations}  
& 4,000 & 4,000 & 4,000 & 10,000    \\ 
{\# L-BFGS iterations} 
& 200   & 500   & 800   & 0         \\ 
{\# wave numbers \(M\)}
& 1     & 3     & 5     & 2         \\ \hline
\end{tabular}
\caption{Hyper-parameters used in the example of Section~\ref{subsec:BL} for \( \varepsilon = 1\).}
\label{tbl:example_p13}
\end{table}

\begin{figure}[tb]
\centering
\includegraphics[width=0.32\linewidth]{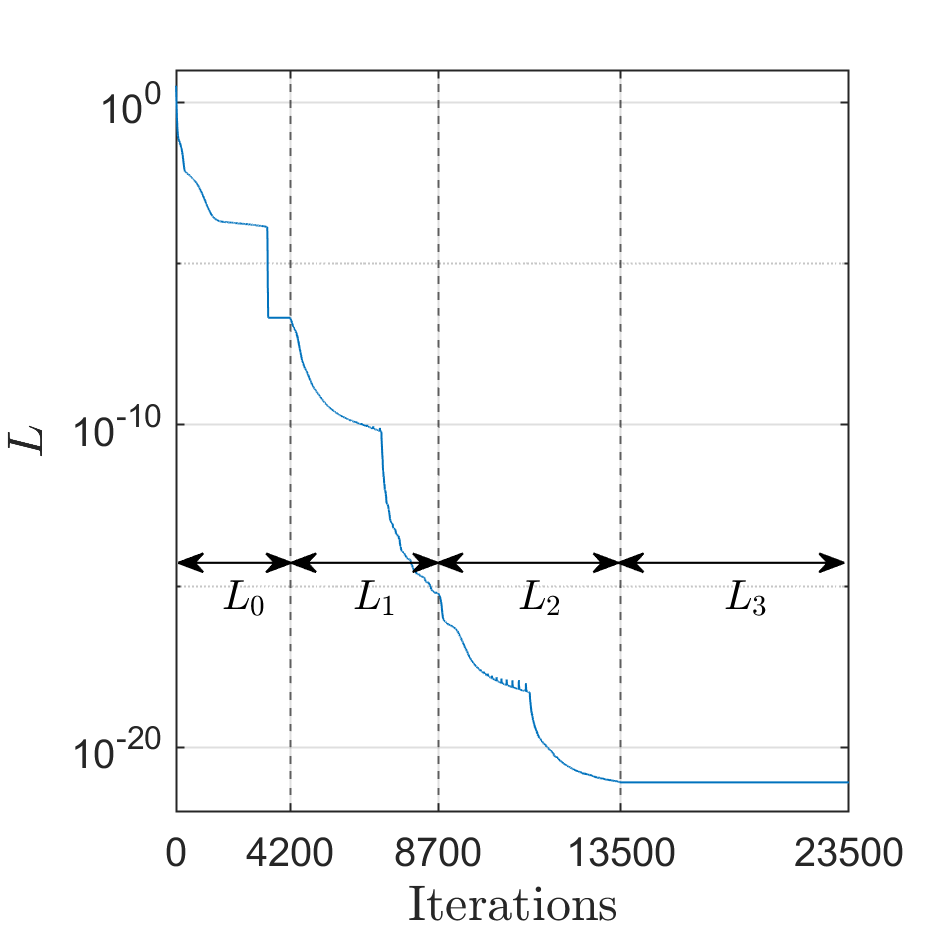}
\includegraphics[width=0.32\linewidth]{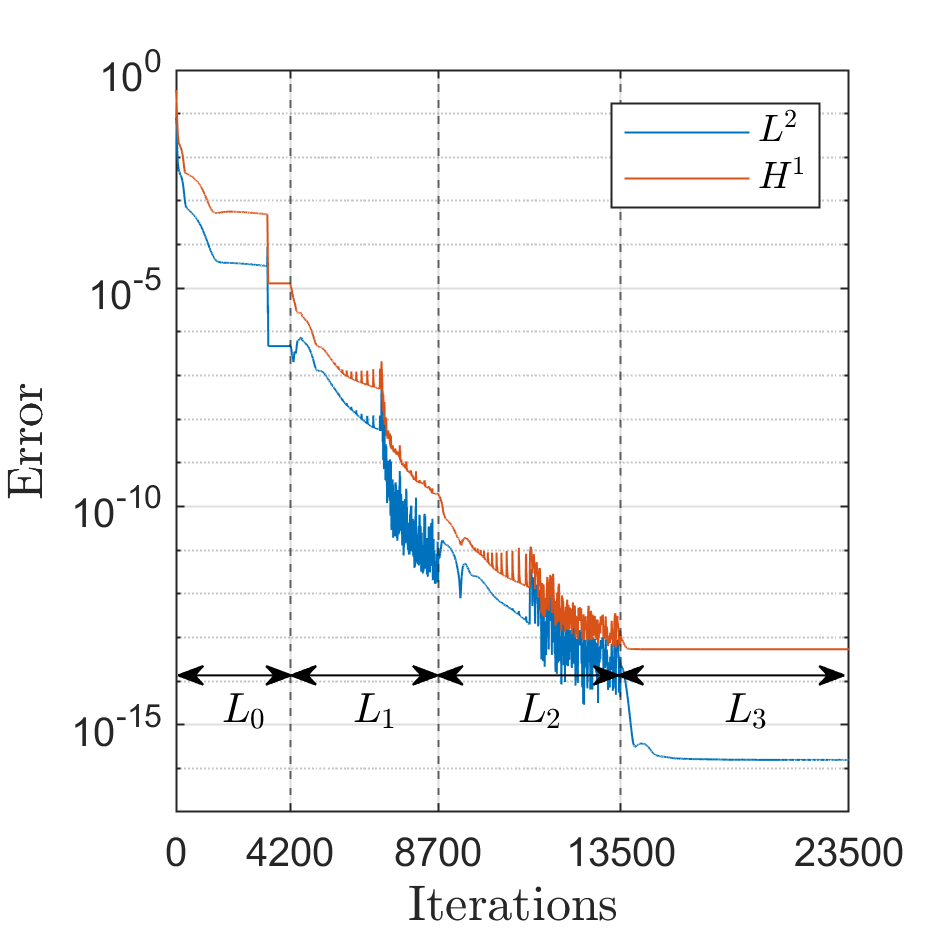}
\includegraphics[width=0.32\linewidth]{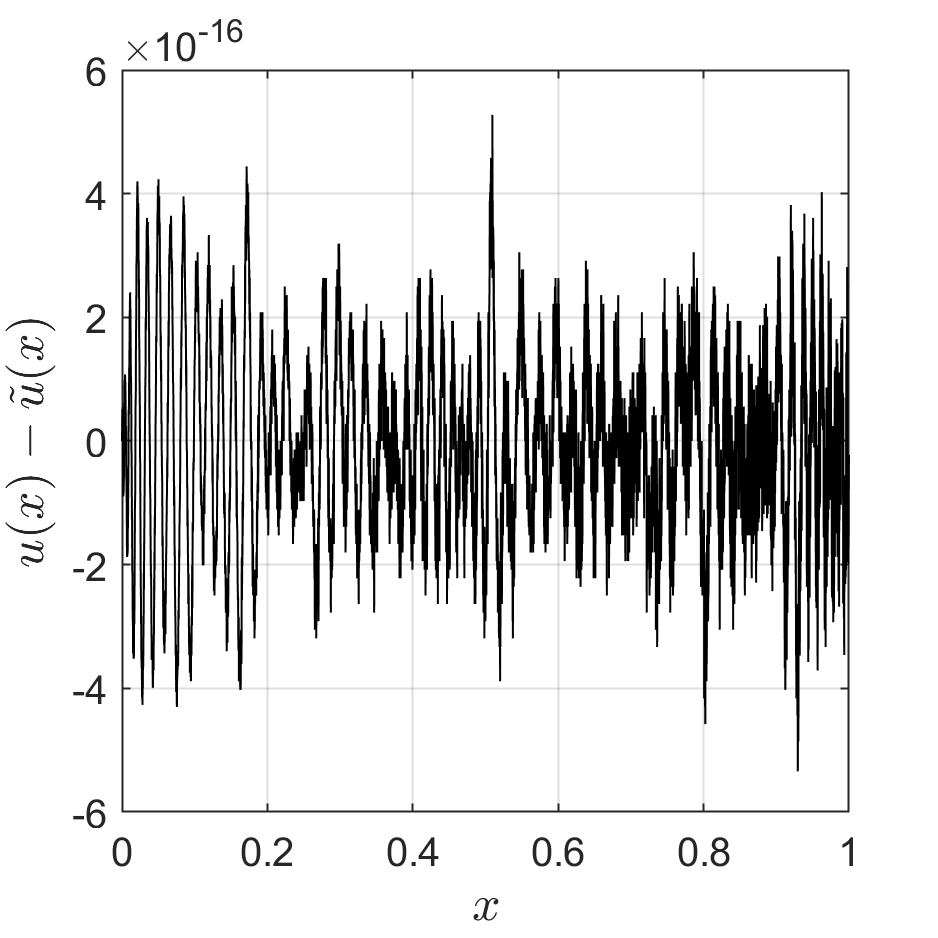}
\caption{Example of Section~\ref{subsec:BL} with \(\varepsilon = 1\): (Left) Evolution of the loss function. (Middle) Evolution of the error \(e(x)=u(x)-\tilde{u}(x)\) measured in the \(L^2\) and \(H^1\) norms. (Right) Pointwise error after three error corrections. 
}
\label{fig:example_p13}
\end{figure}

We first plot in Figure~\ref{fig:example_p13} the convergence results and the pointwise error for \(\varepsilon = 1\). We observe that in this case, we were able to gain at least eight orders of magnitude in the \(L^2\) and the \(H^1\) errors with three error corrections. At the end of the training, the pointwise error is of the order of the machine precision.

\begin{figure}[tb]
\centering
\includegraphics[width=0.32\linewidth]{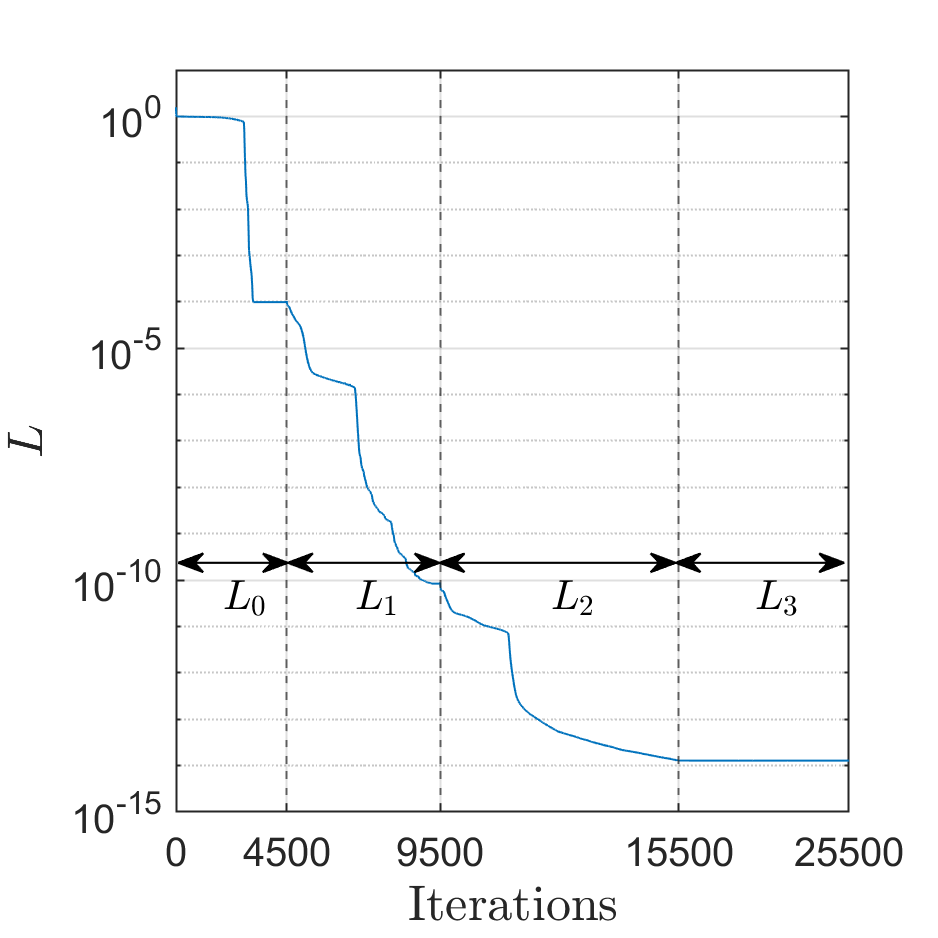}
\includegraphics[width=0.32\linewidth]{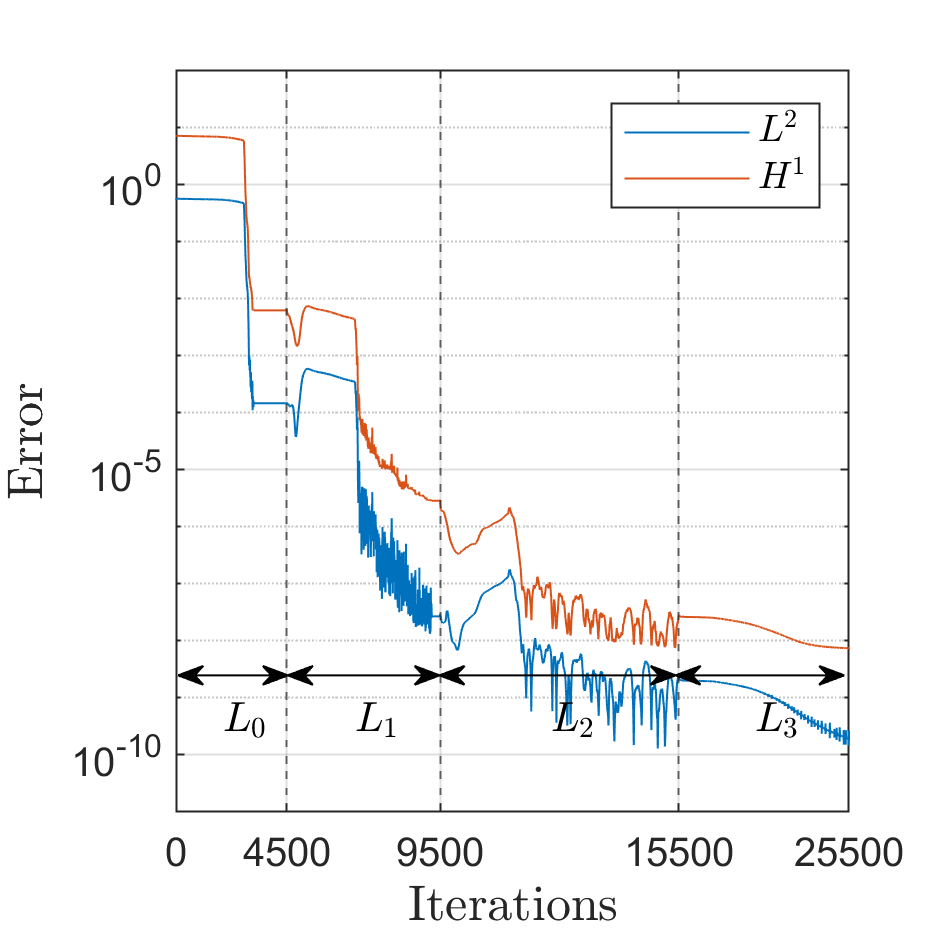}
\includegraphics[width=0.32\linewidth]{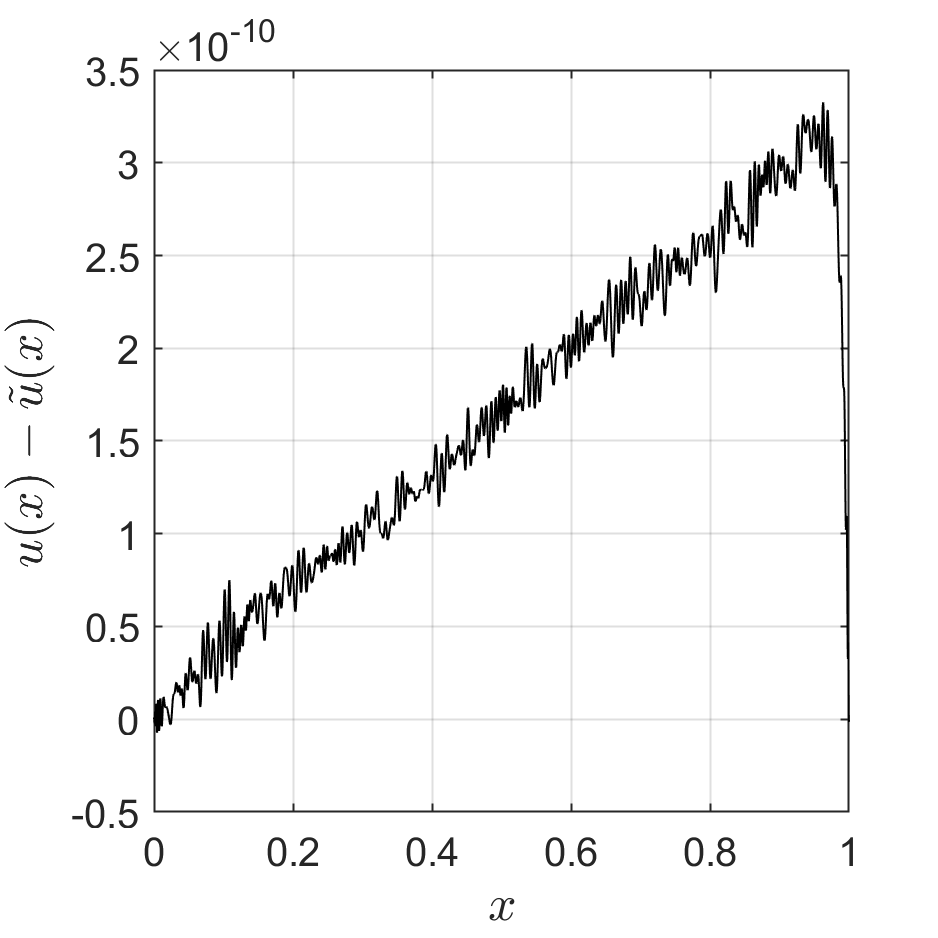}
\caption{Example of Section~\ref{subsec:BL} with \(\varepsilon = 0.01\): (Left) Evolution of the loss function. (Middle) Evolution of the error \(e(x)=u(x)-\tilde{u}(x)\) measured in the \(L^2\) and \(H^1\) norms. (Right) Pointwise error after three error corrections.}
\label{fig:example_p14}
\end{figure}

In Figure~\ref{fig:example_p14}, we show the convergence results and the pointwise error for \(\varepsilon = 0.01\). As expected, the convergence with four networks is slower than in the case with \(\varepsilon = 1\) and plateaus for larger values of the loss function and the errors. As a matter of fact, the loss function after the training of \(\tilde{u}_0\) stagnates around a value of \(10^{-4}\). But using the multi-level neural network method, we are able to decrease the loss function down to \(10^{-14}\), which is also accompanied by a reduction of the \(L^2\) and \(H^1\) errors.

\begin{table}[tb]
\centering
\renewcommand{\arraystretch}{1.4}
\begin{tabular}{c|c|c|c|c}
Hyper-parameters
& \(\tilde{u}_0\) & \(\tilde{u}_1\) 
& \(\tilde{u}_2\) & \(\tilde{u}_3\) \\ \hline
{\# Hidden layers \(n\)}    
& 1     & 1     & 1     & 1         \\ 
{Width \(N_1\)}            
& 10    & 10    & 20    & 20        \\ 
{\# Adam iterations}  
& 4,000 & 4,000 & 4,000 & 10,000    \\ 
{\# L-BFGS iterations} 
& 500   & 1,000 & 2,000 & 0         \\ 
{\# wave numbers \(M\)}
& 3     & 5     & 7     & 3         \\ \hline
\end{tabular}
\caption{Hyper-parameters used in the example of Section~\ref{subsec:BL} for \( \varepsilon = 0.01\).}
\label{tbl:example_p14}
\end{table}

%============================================================
\subsection{Helmholtz Equation}
\label{subsec:helmholtz}

We are now looking for the field \(u=u(x)\) governed by the one-dimensional Helmholtz equation:
\begin{equation}
\label{eq:helmholtz}
- \partial_{xx} u(x) - \kappa^2 u(x)  = 0, \qquad \forall x\in(0,1),
\end{equation}
where the value of the wave number is chosen here equal to \(\kappa=\sqrt{9200}\approx 95.91\),
and subject to the Dirichlet boundary conditions:
\begin{equation}
\begin{aligned}
u(0) &=0, \\
u(1) &=1.
\end{aligned}
\end{equation}
The Dirichlet boundary condition is non-homogeneous at \(x=1\). We thus introduce the lift function \(\bar{u}(x) = x\) to account for the boundary condition, so that we consider trial functions for the initial solution \(\tilde{u}_0\) in the form:
\[
\tilde{u}_0(x) = x + \dfrac{\big(\Pi_{j=1}^d\gamma_g(x_j)\big) \cdot \bz_{n+1}}{M}
\]
where \(\gamma_g\) is defined~\eqref{eq:gamma_g}. Since \(\tilde{u}_0\) strongly verifies the two boundary conditions, the corrections \(\tilde{u}_i\), for \(i \geq 1\), will therefore be subjected to homogeneous Dirichlet boundary conditions, i.e.\ \(\tilde{u}_i(0) = \tilde{u}_i(1) = 0 \). The main objective of this example is to show that the multi-level neural network method can actually recover a high-frequency solution resulting from the large value of the wave number \(\kappa\). 
The hyper-parameters of the multi-level neural networks are provided in Table~\ref{tbl:example_p16}. 

\begin{table}[tb]
\centering
\renewcommand{\arraystretch}{1.4}
\begin{tabular}{c|c|c|c|c}
Hyper-parameters
& \(\tilde{u}_0\) & \(\tilde{u}_1\) 
& \(\tilde{u}_2\) & \(\tilde{u}_3\) \\ \hline
{\# Hidden layers \(n\)}    
& 1     & 1     & 1     & 1         \\ 
{Width \(N_1\)}            
& 10    & 20    & 40    & 10        \\ 
{\# Adam iterations}  
& 10,000 & 10,000 & 10,000 & 30,000    \\ 
{\# L-BFGS iterations} 
& 400   & 800   & 1,600 & 0         \\ 
{\# wave numbers \(M\)}
& 5     & 7     & 9     & 5         \\ \hline
\end{tabular}
\caption{Hyper-parameters used in the example of 
Section~\ref{subsec:helmholtz}.}
\label{tbl:example_p16}
\end{table}

\begin{figure}[tb]
\centering
\includegraphics[width=0.32\linewidth]{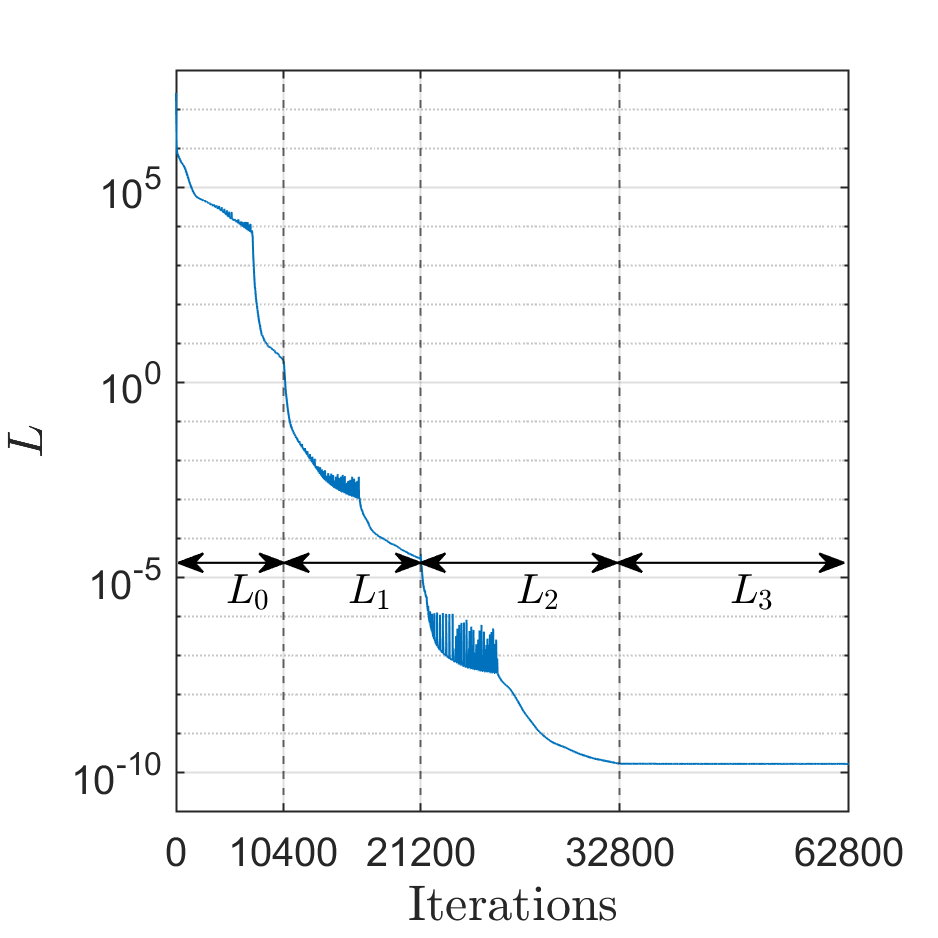}
\includegraphics[width=0.32\linewidth]{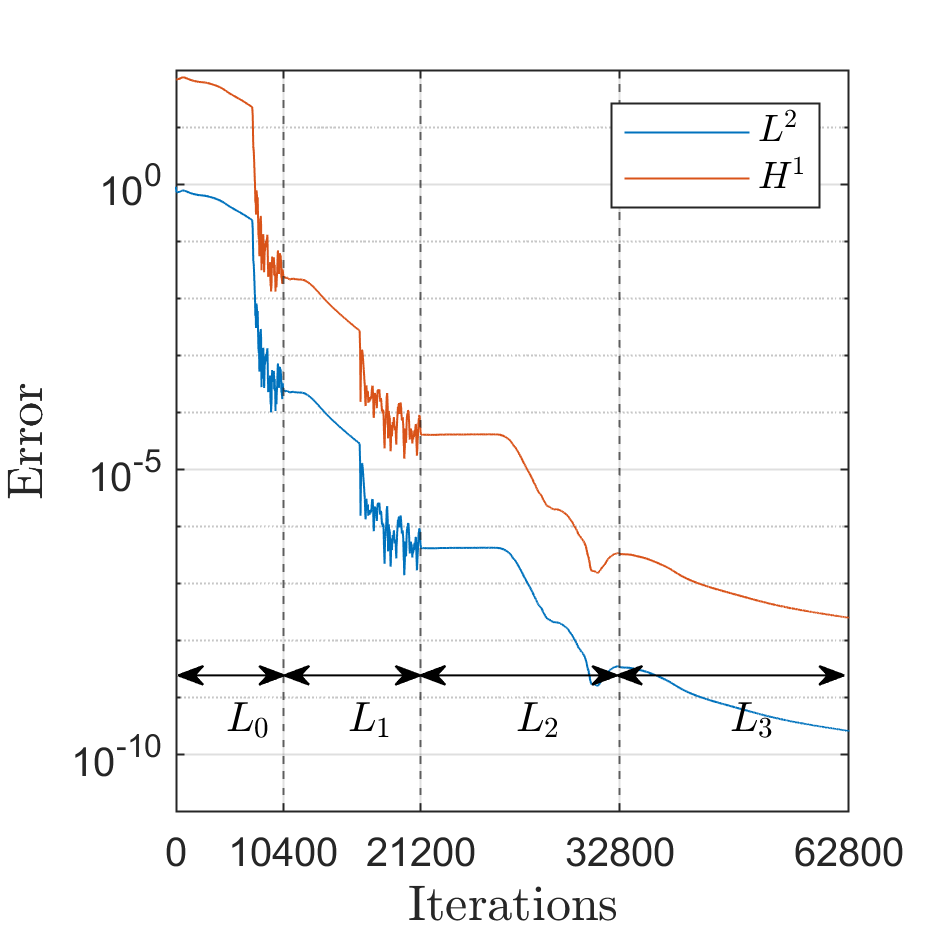}
\includegraphics[width=0.32\linewidth]{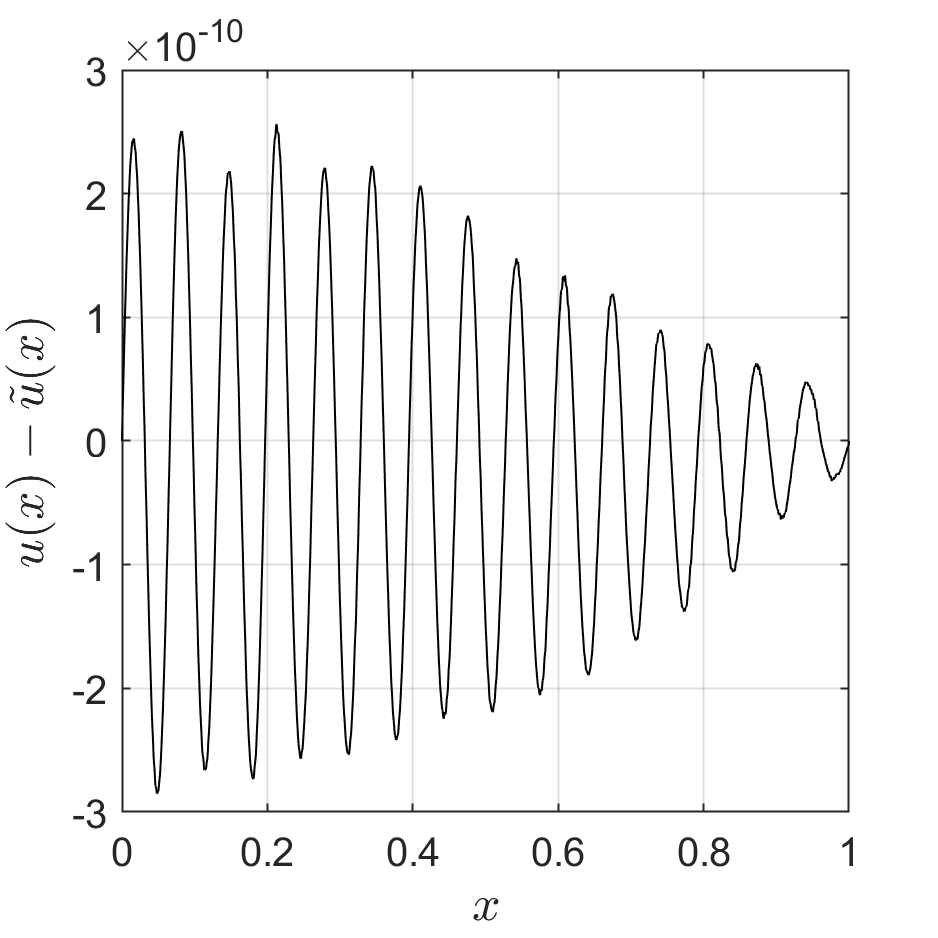}
\caption{Example of Section~\ref{subsec:helmholtz}: (Left) Evolution of the loss function. (Middle) Evolution of the error \(e(x)=u(x)-\tilde{u}(x)\) measured in the \(L^2\) and \(H^1\) norms. (Right) Pointwise error after three error corrections. }
\label{fig:example_p16_1}
\end{figure}

As before, we plot in Figure~\ref{fig:example_p16_1} the convergence of the loss function and the errors along with the pointwise error. We observe that the use of the multi-level neural networks leads to a significant reduction of the error as the absolute pointwise error in the final approximation \(\tilde{u}\) never exceeds \(3\times 10^{-10}\).

\begin{figure}[tb]
\centering
\includegraphics[width=0.32\linewidth]{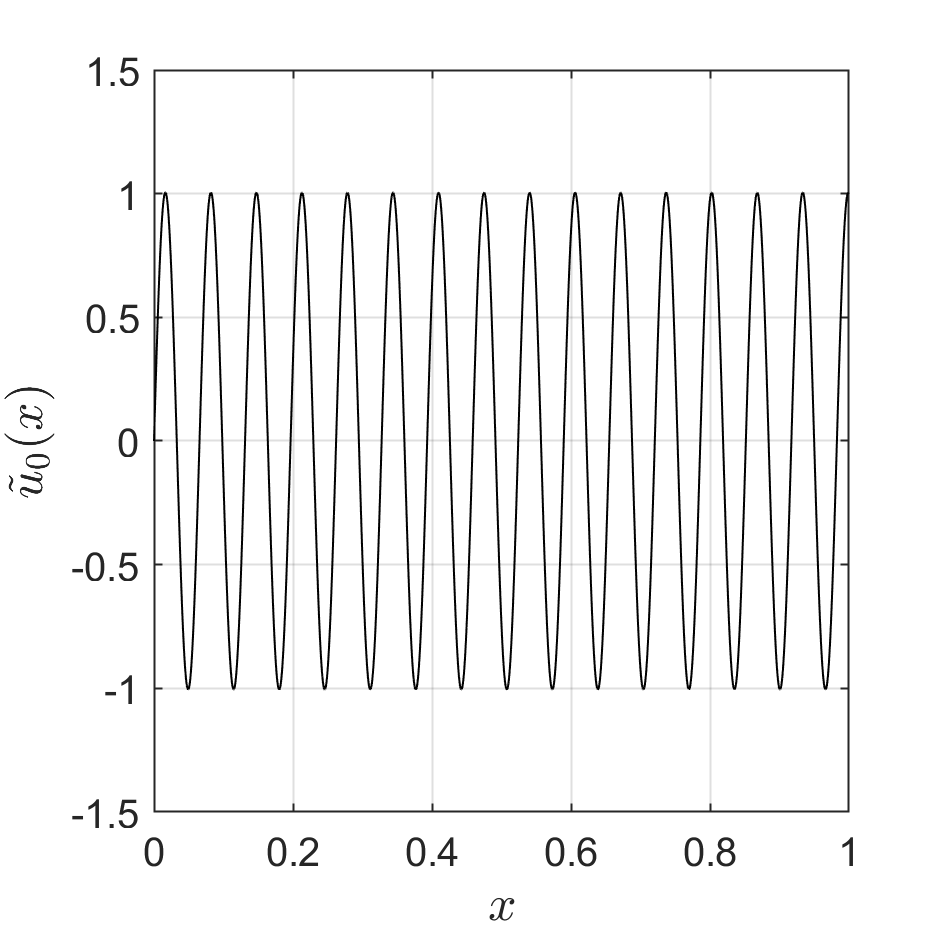}
\includegraphics[width=0.32\linewidth]{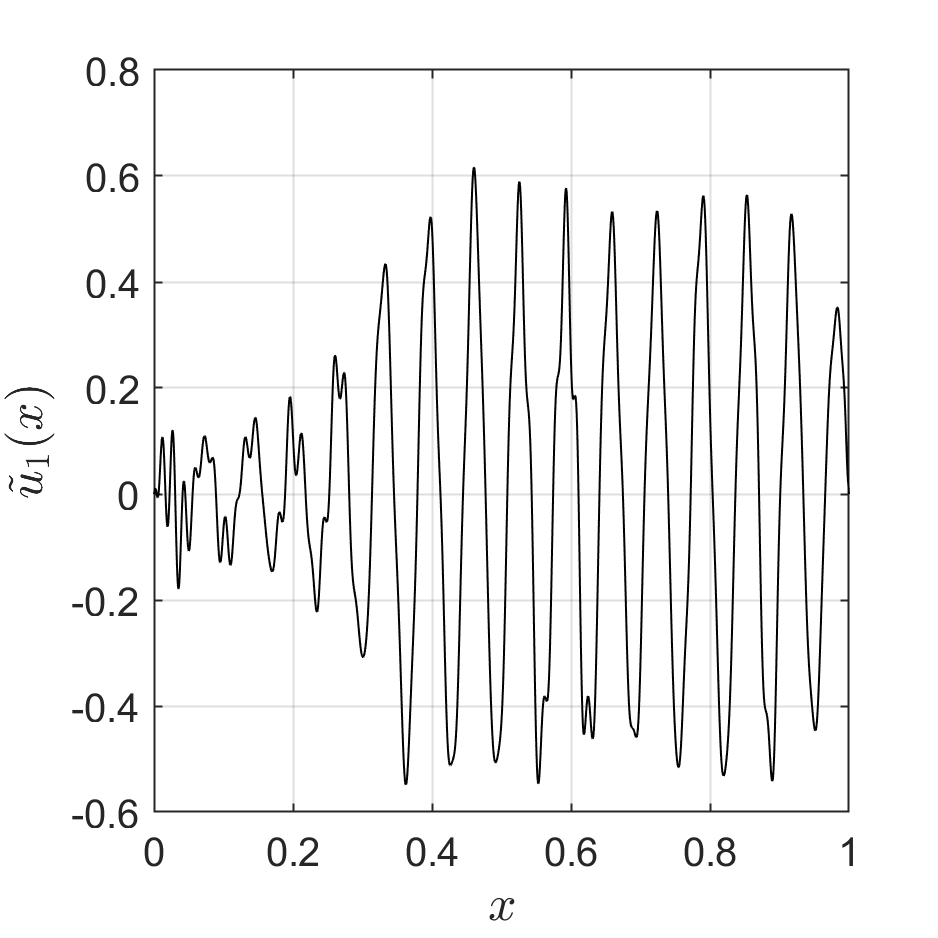}
\includegraphics[width=0.32\linewidth]{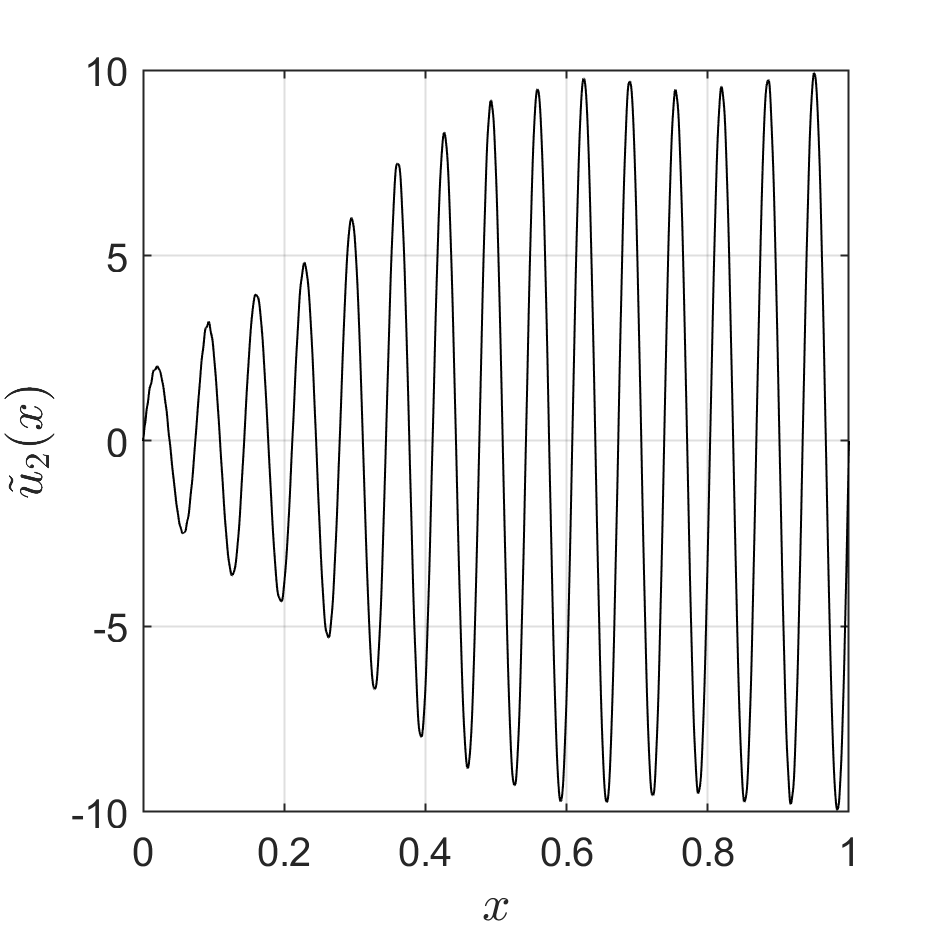}
\caption{Example of Section~\ref{subsec:helmholtz}: Approximation \(\tilde{u}_0(x)\) and corrections \(\tilde{u}_i(x)\), \(i=1,2\).}
\label{fig:example_p16_2}
\end{figure}

\begin{figure}[tb]
\centering
\includegraphics[width=1\linewidth]{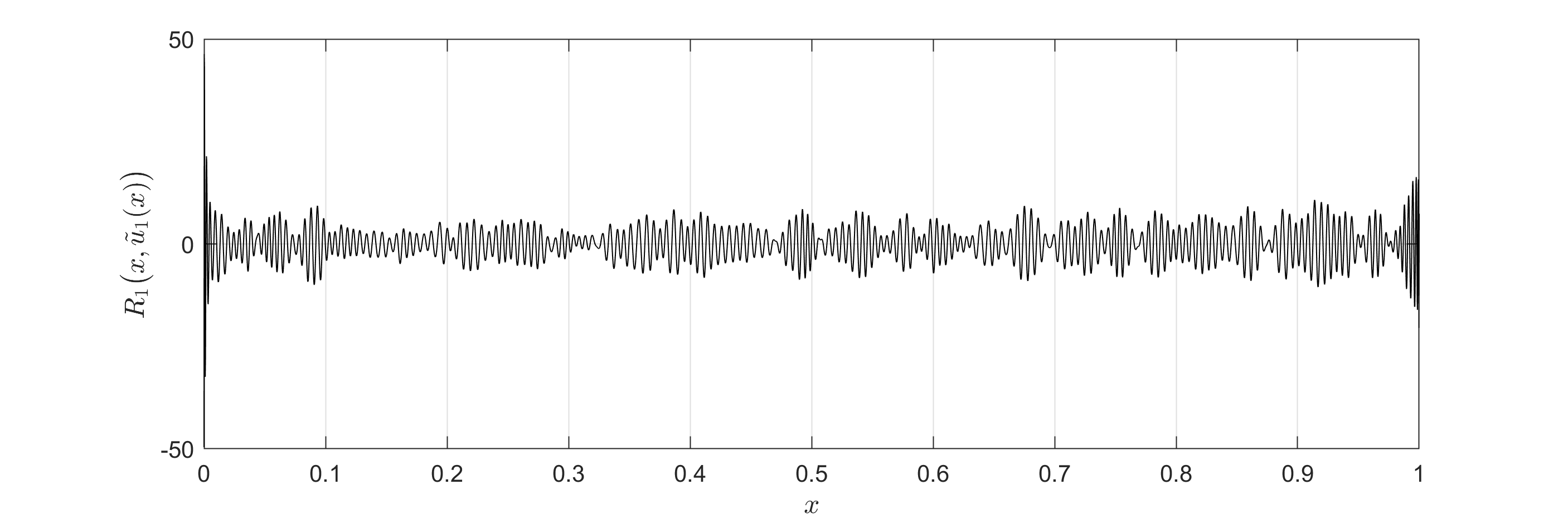}
\caption{Example of Section~\ref{subsec:helmholtz}: Residual \(R_1(x,\tilde{u}_1(x))\) associated with the Hemlholtz equation obtained after the training of \(\tilde{u}_1(x)\).}
\label{fig:example_p16_3}
\end{figure}

In this example, the last correction is constructed using the Fourier feature mapping with \(M=4\) wave numbers. This is in contrast to the previous examples where we chose lower frequencies for the last correction. The reason is that, as observed in Figure~\ref{fig:example_p16_2}, the dominant frequency in \(\tilde{u}_1\) and \(\tilde{u}_2\) is comparable to that of the solution. Therefore, in order to reduce this frequency without reducing the larger frequencies whose amplitudes are smaller, we actually select an architecture for \(\tilde{u}_3\) similar to that of \(\tilde{u}_0\). Using this architecture, we see that the errors in the solution significantly decrease even if the loss function remains virtually unchanged. We show in Figure~\ref{fig:example_p16_3} the residual \(R_1(x,\tilde{u}_1(x))\) associated with the approximation \(\tilde{u}_1\). We have already mentioned that the error corrections \(\tilde{u}_1\) and \(\tilde{u}_2\) have a dominant frequency similar to that of \(\tilde{u}_0\). However, we observe that the residual clearly features higher-frequency modes, whose amplitudes, although small in the approximation \(\tilde{u}_1\), are in fact amplified due to the second-order derivatives. For that reason, it is desirable to consider larger networks with a larger number of wave numbers in the Fourier feature mapping for the approximations \(\tilde{u}_1\) and \(\tilde{u}_2\). We have thus used in this experiment \(N_1=20\) and \(M=7\) for \(\tilde{u}_1\) and \(N_1=40\) and \(M=9\) for \(\tilde{u}_2\).

%============================================================
\subsection{Poisson problem in 2D}
\label{subsec:Poisson2D}

In this final example, we consider the two-dimensional Poisson equation in \(\Omega=(0,1)^2\), with homogeneous Dirichlet boundary conditions prescribed on the boundary \(\partial \Omega\) of the domain. The boundary-value problem consists then in solving for \(u=u(x,y)\) satisfying:
\begin{equation}
\label{eq:poisson_2D}
\begin{aligned}
- \laplacian u(x,y) = f(x,y),& \qquad \forall x\in\Omega, \\
u(x,y) = 0,& \qquad \forall x \in \partial \Omega, 
\end{aligned}
\end{equation}
where the source term \(f(x,y)\) is chosen such that the exact solution is given by:
\[
u(x,y) = \sin(\pi x)\sin(\pi y).
\]

\begin{table}[tb]
\centering
\renewcommand{\arraystretch}{1.4}
\begin{tabular}{c|c|c|c|c}
Hyper-parameters
& \(\tilde{u}_0\) & \(\tilde{u}_1\) 
& \(\tilde{u}_2\) & \(\tilde{u}_3\) \\ \hline
{\# Hidden layers \(n\)}    
& 2     & 2     & 2     & 1         \\ 
{Widths \(N_1\) and \(N_2\)}            
& 10    & 20    & 40    & 40        \\ 
{\# Adam iterations}  
& 2,500 & 5,000 & 10,000 & 4,000    \\ 
{\# L-BFGS iterations} 
& 200   & 400   & 600   & 0         \\ 
{\# wave numbers \(M\)}
& 1     & 3     & 5     & 1         \\ \hline
\end{tabular}
\caption{Hyper-parameters used in the example of Section~\ref{subsec:Poisson2D}.}
\label{tbl:example_p15}
\end{table}

The problem is solved using four networks whose hyper-parameters are given in Table~\ref{tbl:example_p15}. We note that, for this two-dimensional problem, we increase the depth of the networks at levels 0, 1, and 2, to two hidden layers, both having the same width \(N_1=N_2\) at each level.

\begin{figure}[tb]
\centering
\includegraphics[width=0.3\linewidth]{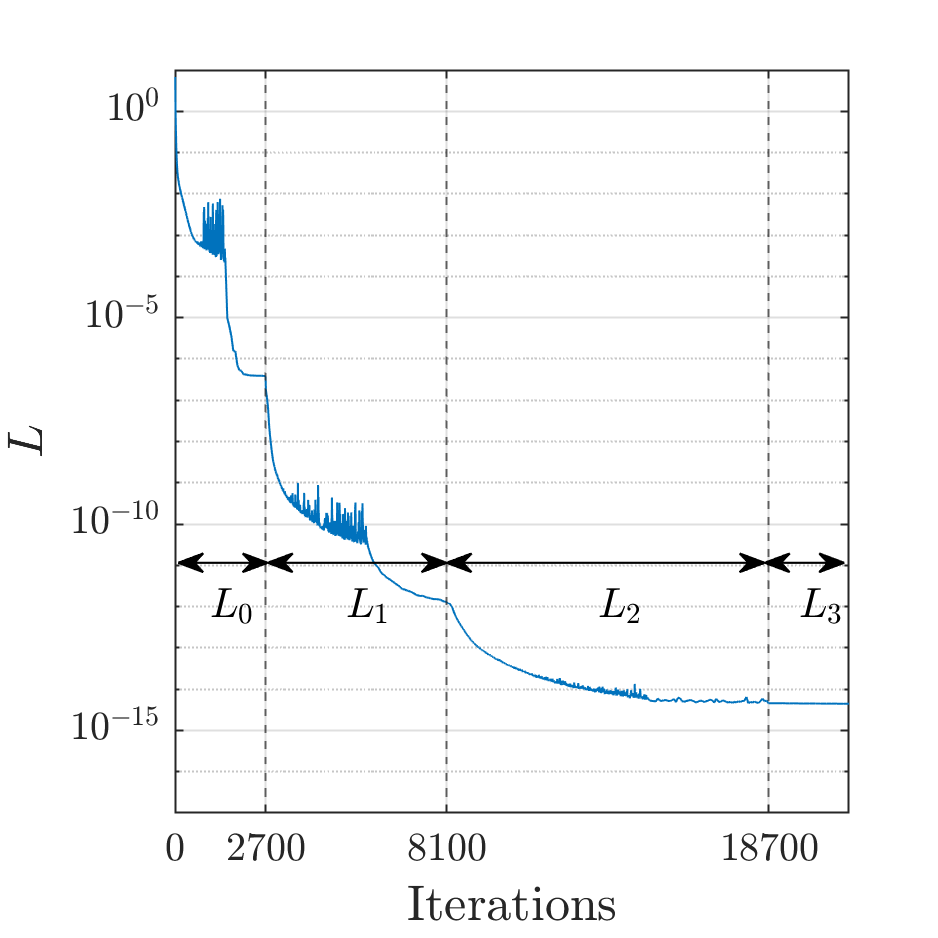}
\includegraphics[width=0.3\linewidth]{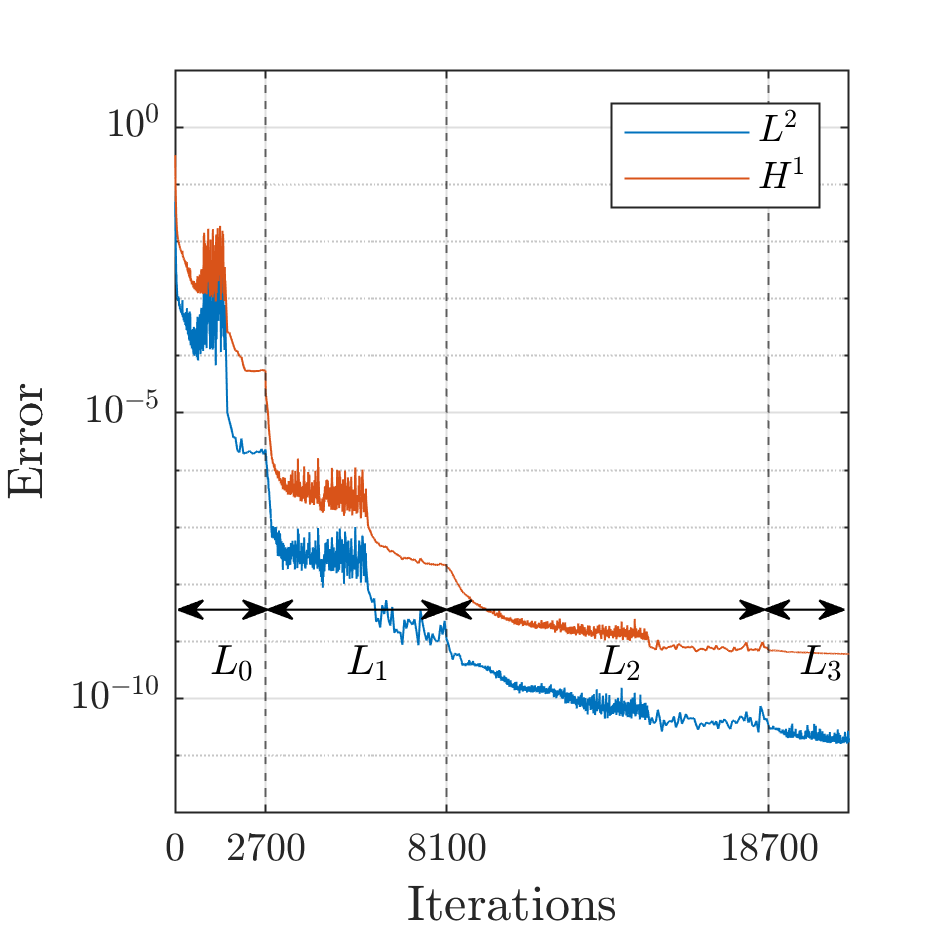}
\includegraphics[width=0.38\linewidth]{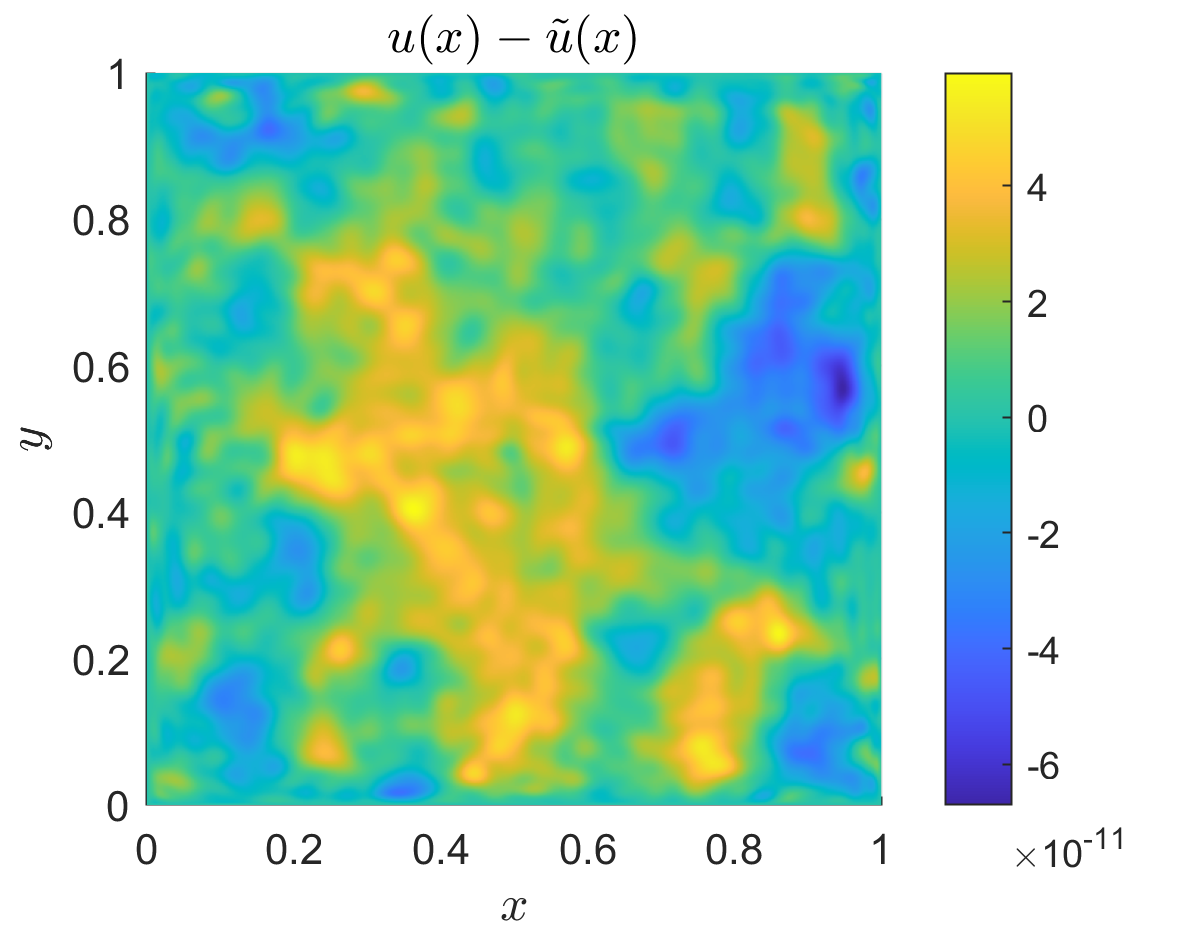}
\caption{Example of Section~\ref{subsec:Poisson2D}: (Left) Evolution of the loss function. (Middle) Evolution of the error \(e(x)=u(x)-\tilde{u}(x)\) measured in the \(L^2\) and \(H^1\) norms. (Right) Pointwise error after three error corrections.}
\label{fig:example_p15_1}
\end{figure}

\begin{figure}[tb]
\centering
\includegraphics[width=0.4\linewidth]{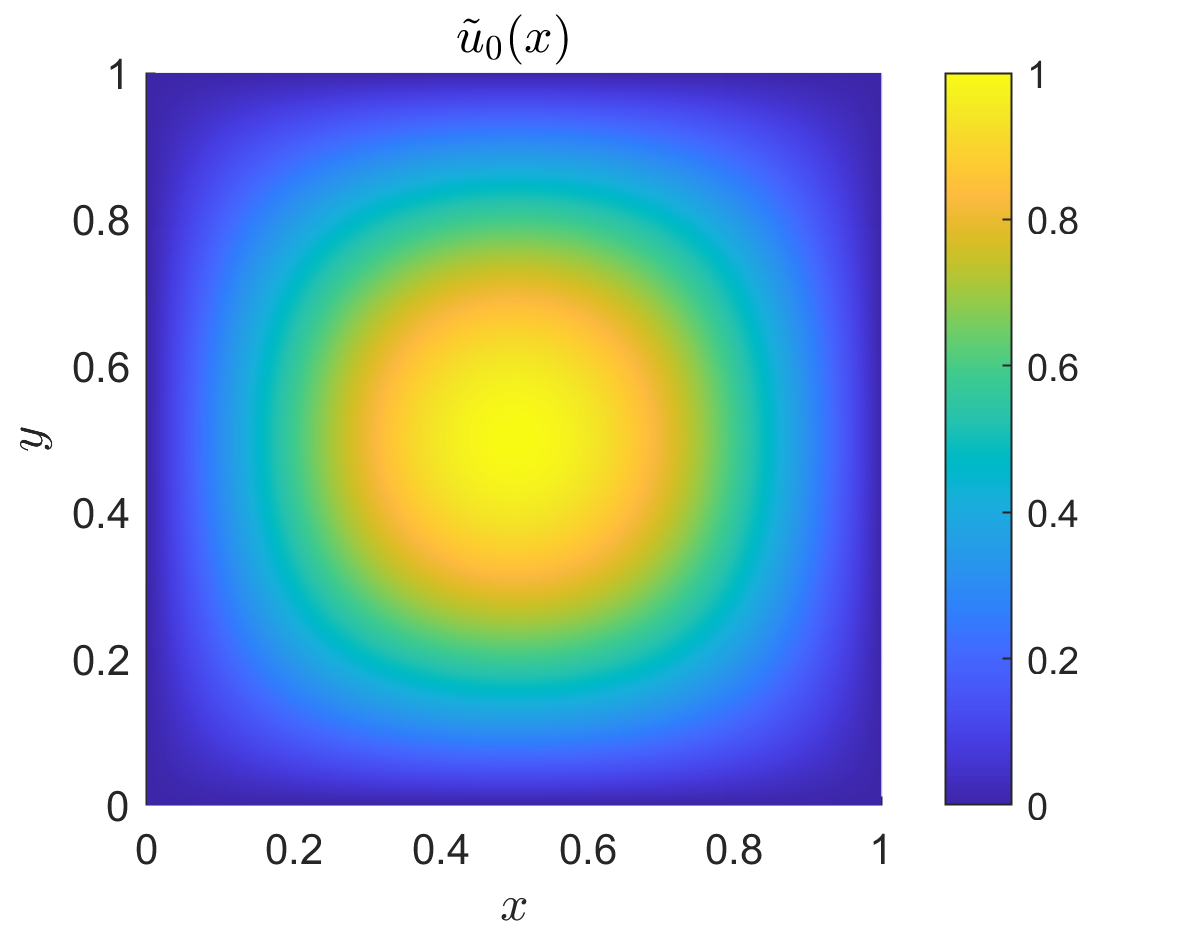}
\includegraphics[width=0.4\linewidth]{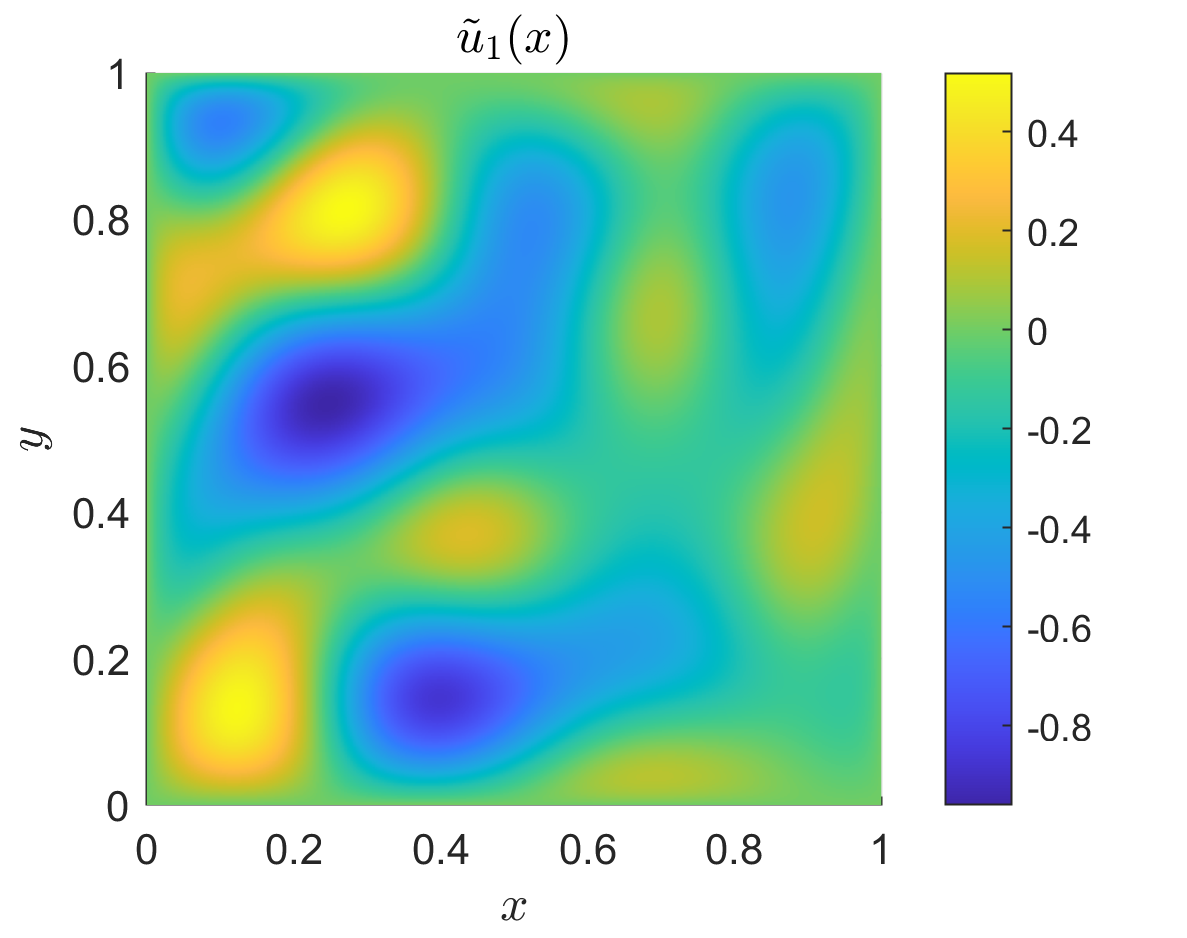}
\\
\includegraphics[width=0.4\linewidth]{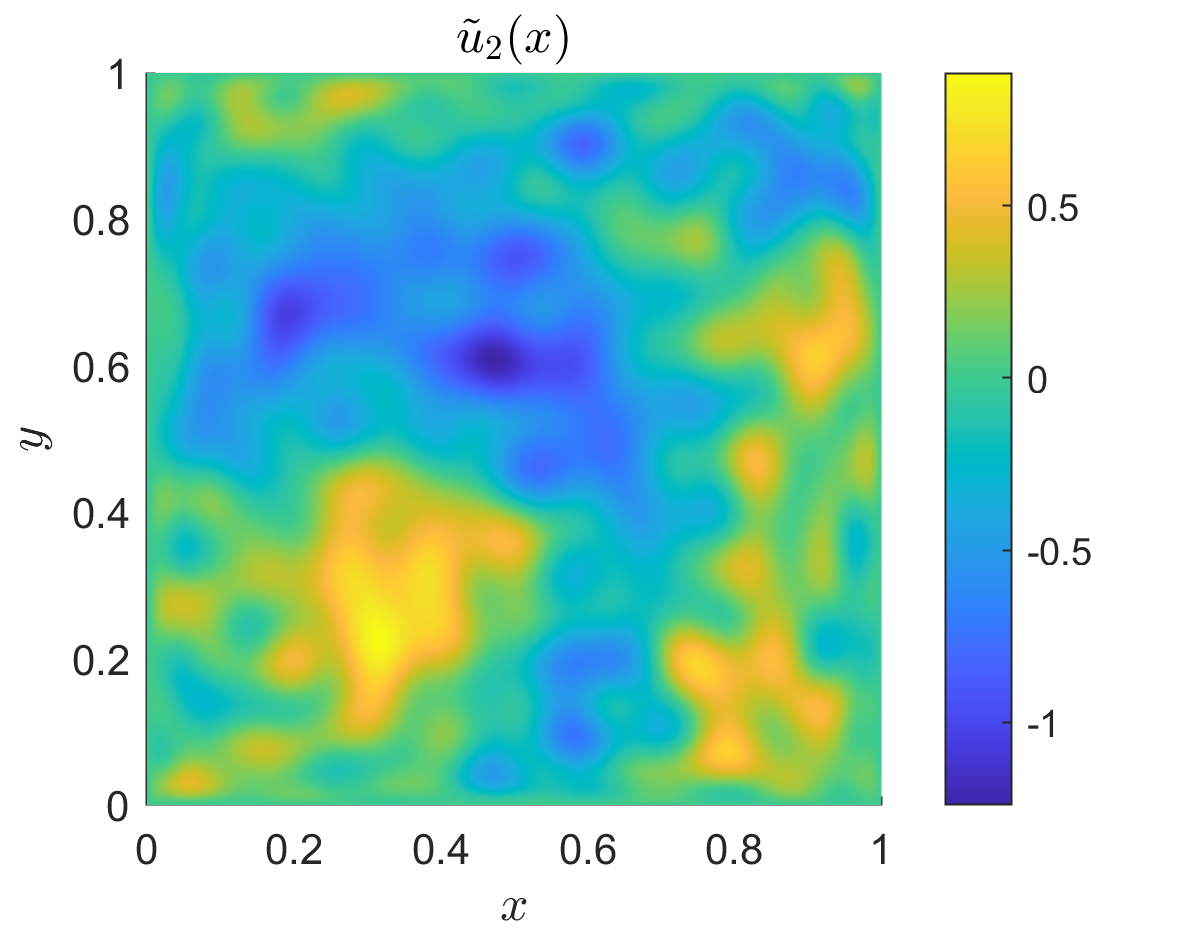}
\includegraphics[width=0.4\linewidth]{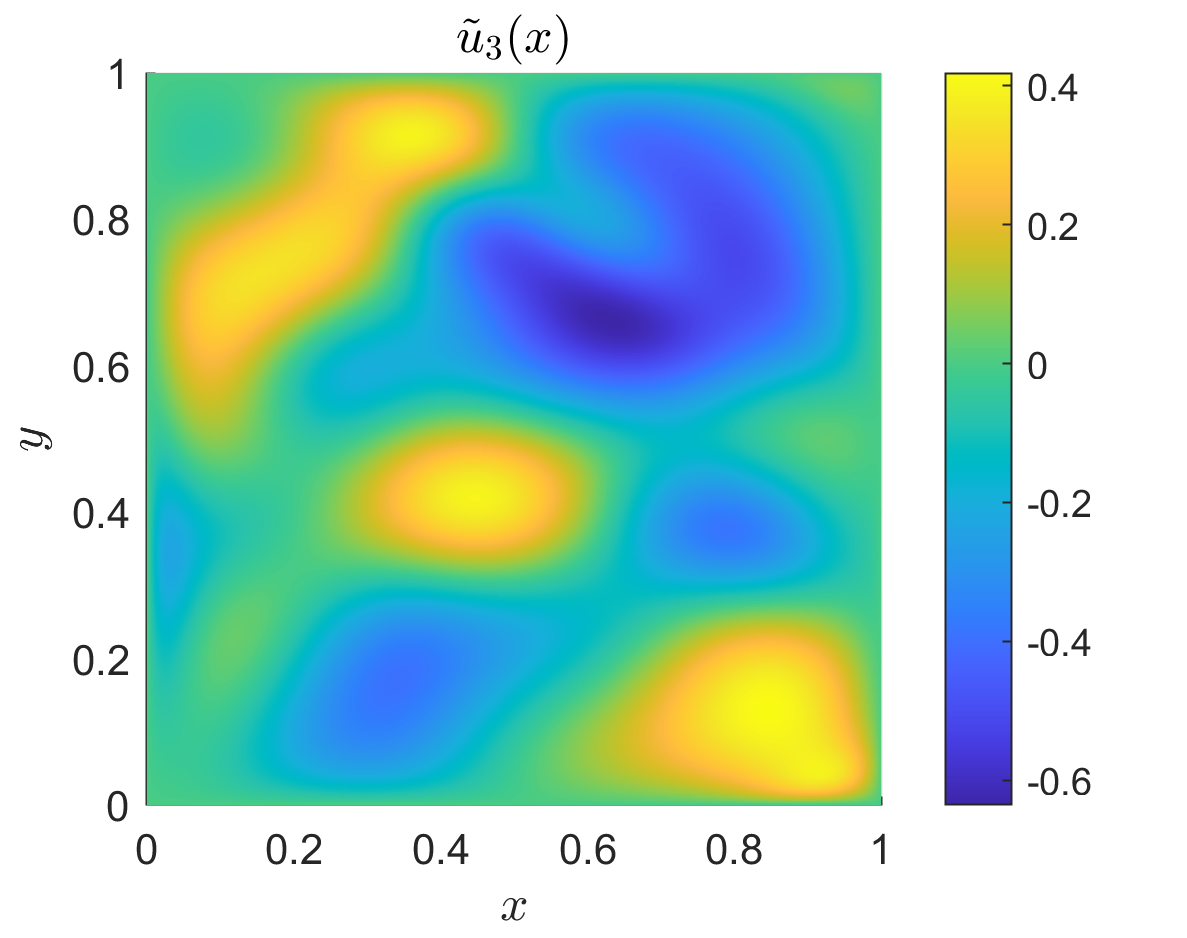}
\caption{Example of Section~\ref{subsec:Poisson2D}: approximation \(\tilde{u}_0(x)\) and corrections \(\tilde{u}_i(x)\), \(i=1,2,3\). }
\label{fig:example_p15_2}
\end{figure}

We show in Figure~\ref{fig:example_p15_1} the evolution of the loss function and of the errors with respect to the number of Adam and L-BFGS iterations. As in the one-dimensional examples, the multi-level neural network approach allows one to reduce the loss function and the errors in the \(L^2\) and \(H^1\) norms down to values around \(10^{-15}\), \(10^{-11}\), and \(10^{-9}\), respectively. The results are in our opinion remarkable since we attain in this 2D example an accuracy comparable to that obtained with classical discretization methods. As indicated in Table~\ref{tbl:example_p15}, the hyper-parameters for the last correction are chosen so that they can capture the low-frequency functions. Figure~\ref{fig:example_p15_1} (right) actually shows that, by the end of the process, we are thereby able to decrease the maximum pointwise error to within \(6\times 10^{-10}\). Finally, we plot in Figure~\ref{fig:example_p15_2} the approximation~\(\tilde{u}_0\) along with  the three corrections~\(\tilde{u}\), \(i=1,2,3\), computed after each level of the multi-level neural networks. One easily observes that all solutions are properly normalized and that, as expected, each corrective function exhibits higher frequencies than in the previous one, except for the approximation~\(\tilde{u}_3\) by the design of the last neural network.

%============================================================
\section{Conclusions}
\label{sec:conclusion}
%============================================================

We have presented in this paper a novel approach to control and reduce errors in deep learning approximations of solutions to linear boundary-value problems. The method has been referred to as the multi-level neural network method in the sense that, at each level of the process, one uses a new neural network, possibly of different sizes, to compute a correction corresponding to the error in the previous approximation. Each successive correction aims at reducing the global error in the resulting approximation of the solution. Although the conceptual idea seems straightforward, the efficiency of the approach relies nonetheless on two key ingredients. Indeed, we have observed that the remaining error at each subsequent level, and equivalently, the resulting residual, have smaller amplitudes and contain higher frequency modes, two circumstances for which we have highlighted the fact that the training of the neural networks usually performs poorly. We have addressed the first issue by normalizing the residual before computing a new correction. To do so, we have developed a normalization approach based on the Extreme Learning Method that allows one to estimate appropriate scaling parameters. The second issue is taken care of by applying a Fourier feature mapping to the input data and the functions used to strongly impose the Dirichlet boundary conditions. We believe that the multi-level neural network method is a versatile approach and can be applied to many deep learning techniques designed to solve boundary-value problems. In this work, we have chosen to present the method in the special case of physics-informal neural networks, which have recently been used for the solution of several classes of initial and boundary-value problems. The efficiency of the multi-level neural network method was demonstrated here on several 1D or 2D numerical examples based on the Poisson equation, the convective-diffusion equation, and the Helmholtz equation. More specifically, the numerical results successfully illustrate the fact that the method can provide highly accurate approximations to the solution of the problems and, in some cases, allows one to reduce the numerical errors in the \(L^2\) and \(H^1\) norms down to the machine precision. 

Even if the preliminary results are very encouraging, additional investigations should be considered to further assess and improve the efficiency of the proposed multi-level neural network method. More specifically, one would like to apply the method to other deep learning approaches, such as the Deep Ritz method~\cite{weinan2018deep} or the weak adversarial networks method~\cite{zang2020weak} to name a few, to time-dependent problems, and to the learning of partial differential operators, e.g.\ DeepONets~\cite{lu2019deeponet} or GreenONets~\cite{aldirany2023operator}, for reduced-order modeling. One could also imagine estimating the correction at each level of the algorithm to control the error in specific quantities of interest following ideas from~\cite{prudhomme-oden-2001,kenan-2017,kenan-2019}. In this work, we have chosen to strongly enforce boundary conditions in order to neglect errors arising from the introduction of penalty parameters in the loss function~\eqref{eq:losspinn}. One should thus assess the efficiency of the method when initial and boundary conditions are weakly enforced. Finally, the multi-level neural network method introduces a sequence of several neural networks whose hyper-parameters are chosen a priori and often need to be adjusted by trial and error. It would hence be very useful to devise a methodology that determines optimal values of the hyper-parameters independently of the user.

%============================================================
\paragraph{Acknowledgements.} 
SP and ML are grateful for the support from the Natural Sciences and Engineering Research Council of Canada (NSERC) Discovery Grants [grant numbers RGPIN-2019-7154, PGPIN-2018-06592]. This research was also partially supported by an NSERC Collaborative Research and Development Grant [grant number RDCPJ 522310-17] with the Institut de Recherche en \'Electricit\'e du Qu\'ebec and Prompt. SP acknowledges the support of the Basque Center for Computational Mathematics to host him in May and June 2023. He also thanks David Pardo for many fruitful discussions on this subject. ZA and SP are thankful to the Laboratoire de M\'ecanique et d’Acoustique UMR 7031, in Marseille, France, for hosting them. This work received support from the French Government under the France 2030 investment plan, as part of the Initiative d'Excellence d'Aix-Marseille Université - A*MIDEX -  AMX-19-IET-010.
%============================================================
\bibliographystyle{abbrv}
\bibliography{biblio}
	
%============================================================
\end{document}